\magnification=1200
\parskip 10 pt

 2
\font\caps = cmcsc10
\font\eighttt=cmtt8

\def\longrightarrow{\relbar\joinrel\rightarrow}
\def\longleftarrow{\leftarrow\joinrel\relbar}
\def\mapright#1.{\buildrel #1 \over \longrightarrow}
\def\mapleft#1.{\buildrel #1 \over \longleftarrow}
\def\mapne#1.{\llap{$\vcenter{\hbox{$\scriptstyle#1$}}$}\nearrow}
\def\mapse#1.{\llap{$\vcenter{\hbox{$\scriptstyle#1$}}$}\searrow\ \ }
\def\mapsw#1.{\llap{$\vcenter{\hbox{$\scriptstyle#1$}}$}\swarrow\ \ }
\def\mapup#1.{\Big\uparrow\rlap{$\vcenter{\hbox{$\scriptstyle#1$}}$}}
\def\mapdown#1.{\Big\downarrow\rlap{$\vcenter{\hbox{$\scriptstyle#1$}}$}}

\def\sc#1{{\cal #1}}

\def\Hom{{\rm Hom}}
\def\Map{\underline{Map}\ }
\def \hoco{\mathop{\rm hocolim}}
\def\hocolim{\hoco}
\def \ho{\mathop{\rm holim}}
\def\holim{\ho}
\def \THH{{\rm THH}}

\def\P{{\bf p}}

\def\longlongrightarrow{\hbox to 30pt{\rightarrowfill}}
\def\longlongleftarrow{\hbox to 60pt{\leftarrowfill}}

\def \u {\underline}

\def \id{{\rm id}}

\def \hoco{\mathop{\rm hocolim}}
\def \ho{\mathop{\rm holim}}
\def \hofib{\mathop{\rm hofib}}
\def \THH{{\rm THH}}

\def\Z{{\bf Z}}

\def \Fp{{\bf F}_p}

\def\Res{{\rm Res}}

\def\oF{{\hat{\otimes}_F}}

\def\U{{\cal U}}

\def\u{\underline}
\def\Map{{\u{\rm Map}}}
\def\oc{{\overline c}}

\def\semiprod{\hbox{$\times$\hskip -6.7pt\vrule height5pt
width 0.4pt depth 0pt \hskip 6.5pt}}

\centerline{\bf On the Taylor Tower of Relative K-theory}
\bigskip
\centerline{Ayelet  Lindenstrauss}
\centerline{Department of Mathematics}
\centerline{Indiana University}
\centerline{Bloomington IN 47405}
\centerline{{\eighttt ayelet@math.indiana.edu}}
\medskip
\centerline{Randy McCarthy\footnote*{Partially supported by NSF grant DMS   03-06429}}
\centerline{Department of Mathematics}
\centerline{University of Illinois at Urbana-Champaign}
\centerline{Urbana IL 61801}
\centerline{{\eighttt randy@math.uiuc.edu}}

\bigskip

\def\obj{{\it obj}}

\noindent{\bf 0. \underbar{Introduction}}

\bigskip
Tom Goodwillie conjectured that one could define a Hochschild homology theory for rings where the tensors are replaced by `tensors' over the sphere spectrum, and get an invariant which would agree with the stable K-theory of the ring.  This theory, now known as topological Hochschild homology and denoted by $\THH$, was defined by Marcel B\"okstedt  [B], and in [D], Bj\o rn Dundas  proved that it agrees, as conjectured, with stable K-theory.

The definition of the stable K-theory of a ring $R$ is a shift by one dimension of the definition of the first derivative at a point, in the sense of Goodwillie's calculus of functors, of the functor which sends a pointed simplicial set $X$ to the spectrum $K(R\semiprod\tilde R[X])$.  (Here $\tilde R[X]$ is the  simplicial abelian group which has the free $R$-module on $X_n$ modulo $R\cdot *$ in dimension $n$; $\semiprod$ denotes the trivial square-zero extension, with zero multiplication on the adjoined ideal.)  The proof  incorporates Waldhausen's S-construction to show the topological Hochschild homology is also a shift by one dimension of the derivative of the functor $K(R\semiprod\tilde R[-])$.

Our goal in this paper is to describe the full Taylor tower of the functor $K(R\semiprod\tilde R[-])$.  As in [DMc1], it turns out that the method in fact applies for studying $K(R\semiprod\tilde M[-])$ for any $R$-bimodule $M$, but that the more natural object to construct a Taylor tower for is algebraic K-theory with coefficients in a bimodule, introduced there and reviewed in Section 9 below.  In Theorem 4.1 of [DMc1], the authors define a natural weak equivalence
$$K(R;B_.M)\mapright\simeq. K(R\semiprod M)$$
which (since all the invariants we discuss for simplicial bi-modules are calculated levelwise) means that
$$K(R;\tilde M[\Sigma X])\simeq K(R\semiprod\tilde M[X])\leqno(0.1)$$ 
for any space $X$.  We therefore study the functor $K(R;\tilde M[-])$, which by functoriality is readily decomposed $K(R;\tilde M[-])\simeq K(R)\times\tilde K(R;\tilde M[-])$ (see, again, the beginning of Section 9).  By [DMc1], the  best `linear' (excisive) approximation of the functor $\tilde K(R;\tilde M[-])$ at the one-point space $*$ is $\THH(\u R;\u{\tilde M}[-])$.  

In [BHM], Marcel B\"okstedt, Wu-Chung Hsiang, and Ib Madsen show that the topological Hochschild homology for a functor with smash product $F$ is a Connes cyclic object, with respect to levelwise cyclic rotation of the coordinates.  By the theory of Connes cyclic objects, this gives an $S^1$-action  on $\THH(F)$, and for any positive $n$, the restriction of this action to the cyclic group $C_ns\subset S^1$ acts simplicially on the $n$'th edgewise subdivision, ${\rm sd}^n \THH(F)$.  Introducing an $F$-bimodule $P$ in one of the coordinates when we look at topological Hochschild homology with coefficients $\THH(F;P)$ spoils the Connes cyclic action.  However, in Section 2 below we define $U^n(F;P)$ to be an analog of ${\rm sd}^n \THH(F)$ with $n$ bimodule coefficient coordinates, and we get the same $C_n$ action as one gets on ${\rm sd}^n \THH(F)$ (what is lost is that  ${\rm sd}^n \THH(F)$
is a simplicial spectrum whose realization is homeomorphic to that of $\THH(F)$ and this is not true for
$U^n(F;P)$ and $U^1(F;P)=\THH(F;P)$).  Note that the construction we generalize is not B\"okstedt's original construction of the topological Hochschild homology of an FSP but rather the construction from [DMc1], [DMc2] of the topological Hochschild homology of an FSP over a category.  This will be important in defining the map from K-theory in Section 9.  

As in [BHM], when $m$ divides $n$ we have restriction maps
$$Res^{n/m}: U^n(F;P)^{C_n}\to U^m(F;P)^{C_m}.$$
If one looks at the category $\bf N^\times$ with objects corresponding to the positive integers and one morphism $n\to m$ whenever $m$ divides $n$, and at its full subcategory 
$\{ \leq n\}$, in Section 4 we define 
$$\eqalign{
W_{n}(F;P) & = \holim_{{\{\leq n\}}}\ U^m(F;P)^{C_m},
\cr
W(F;P) & = \holim_{{\bf N^\times}}\ U^m(F;P)^{C_m}.
}
$$
In [BHM] terminology, $W(F;P)$ is an analog of $TR(F)$ with coefficients in a bimodule.  

Our main theorem, Theorem 9.2, states that  a natural transformation 
$$\beta: \tilde K(R;\tilde M[-])\to W(\u R;\u{\tilde M}[-])$$ 
which we construct induces an equivalence between the two functors on connected $X$.  Moreover,
by Corollary 9.3  $ W_n(\u R;\u{\tilde M}[-])$ is the $n$'th stage of the Goodwillie calculus Taylor tower of
the functor $\tilde K(R;\tilde M[-])$, with the tower structure maps the same as those induced on the homotopy inverse limits by the restriction of categories from $\{\leq n\}$ to $\{\leq n-1\}$.  In other words,
$W(\u R;\u{\tilde M}[-])$ is the Goodwillie calculus Taylor tower for the functor  $\tilde K(R;\tilde M[-])$, and the tower converges for connected $X$, which are the ones we care most about because of equation (0.1) above.

This equivalence is interesting even for finite sets $X$, where we use this method, and Lars Hesselholt and Ib Madsen's calculation of $W(\u {\Fp};\u{\Fp})={\rm TR}(\u{\Fp})$ to completely calculate
$\tilde K(\Fp\semiprod(\oplus_{i=1}^n \Fp))^{\wedge}_p$ in 
a separate paper.  

To obtain this result, we start in Section 5 to use analysis like that done by Goodwillie in [G] for $\THH$ and ${\rm TR}$ to calculate $\hofib[W_n(F;P)\to W_{n-1}(F;P)]$ for a general FSP $F$ and an $F$-bimodule $P$.  Corollary 5.7 shows that
$$U^n(F;P)_{hC_n} \simeq   \hofib[W_n(F;P)\to W_{n-1}(F;P)],\leqno(0.2)$$

In Section 8, this analysis is used to show that each layer  $\hofib[W_n(\u R;\u{\tilde M}[-])\to W_{n-1}(\u R;\u{\tilde M}[-])]$ is a homogenous degree $n$ functor (this is basically Corollary 8.2):
when $F=\u R$ and $P=\u{\tilde M}[-]$, 
$$U^n(\u R;\u{\tilde M}[X] )\simeq U^n(\u R;\u M)\wedge(\bigwedge^n X).$$
This implies already that the $W_n(\u R;\u{\tilde M}[-])$ are $n$-excisive functors which are the $n$'th stage
of the Taylor tower of their homotopy inverse limit $W(\u R;\u{\tilde M}[-])$.  

Section 6 contains more properties of $U^n$ and $W$, specifically results involving the varying of the category over which the invariants are constructed: the category consisting of a single point, the category of finitely generated projective $R$-modules, iterations of Waldhausen's S-construction on the latter, and intermediate categories used to go between these.  This section generalizes results from [DMc2], and the results are used to define the map $\beta$ of Theorem 9.2.

Section 7 discussed the Goodwillie calculus properties of our functors, most importantly: showing that 
$W(\u R;\u{\tilde M}[-])$ is $0$-analytic (Proposition 7.14).  It is already known by Proposition 3.2 of [Mc1] that the functor $ \tilde K(R;\tilde M[-])$ is $0$-analytic as well.  

Thus in Section 9, after constructing the natural transformation $\beta$, we can use a variant of Goodwillie's Theorem 5.3 from [G2]: it states that if there is a natural transformation between two $\rho$-analytic functors $F$ and $G$ which induces an equivalence of the differentials at every space $X$, then for $(\rho+1)$-connected maps $X\to Y$, there is a Cartesian square
$$\matrix{F(X)&\mapright .& G(X)\cr
  \mapdown .& &\mapdown .\cr
 F(Y)&\mapright .&G(Y).\cr}$$
 The variant is simply the observation that the proof in [G2] requires an equivalence of the differentials only on $\rho$-connected $X$.  We will want to apply it for $\rho=0$, $Y=*$, and our two functors above to get that for $0$-connected $X$ (that is: $X$ for which $X\to *$ is $1$-connected), 
$$\beta: \tilde K(R;\tilde M[X])\mapright{\simeq}. W(\u R;\u{\tilde M}[X]).$$ 
It follows from [DMc1] that $\beta$ induces an equivalence of the differentials at $*$, and the remainder of the paper uses that result to show an equivalence of the differentials at arbitrary $0$-connected spaces $X$.  Section 9 concludes by reducing our Main Theorem 9.2 to Technical Lemma 9.4, which states that for a ring $R$ and simplicial $R$-bimodules $M$ and $N$ with $N$ $k$-connected, $\beta$ induces a $2k$-connected map 
$${\rm hofib} \bigl(\tilde K(R;B_.M\oplus B_.N) \to \tilde K(R;B_.M)\bigr)
\to 
{\rm hofib} \bigl(  W(\u{R};\u{B_.M}\oplus\u{B_.N})\to  W(\u{R};\u{B_. M})\bigr).$$
(The lemma is actually stated in terms of $W(\u{\sc P_R};-) $ rather than $W(\u{R};-) $, but the two agree by Proposition 6.13.)

The final two sections of the paper prove Technical Lemma 9.4.  The basic strategy is to write
$$K(R;B_.M\oplus B_.N)\simeq K(R\semiprod(M\oplus N)=K((R\semiprod M)\semiprod N)\simeq K(R\semiprod M; B_.N)$$
and observe that by [DMc1], the homotopy fiber of the map from this to $K(R\semiprod M)\simeq K(R;B_.M)$, this admits a $2k$-connected map to $\THH(\u{R\semiprod M};\u{B_.N})$.  Using the multilinearity of each  $\THH_r(\u{R\semiprod M};\u{B_.N})$ in the $r$ FSP coordinates, and careful analysis of dimensions, in Corollary 10.5 we obtain a decomposition 
$$THH(\u{R\semiprod M};\u{B.N})\simeq\prod_{a=0}^{\infty} U^a(\u R;\u{B.M},\ldots,\u{B.M},\u{B.N}).$$
It is much easier to see that up to order $2k+1$, $ W(\u{R};\u{B_.M}\oplus\u{B_.N})$ also decomposes into the same product; this is Lemma 10.3.

Thus Section 10 establishes the plausibility of Technical Lemma 9.4, by demonstrating that the domain and range there are indeed equivalent up to dimension $2k$.  Section 11 then traces the actual map $\beta$ through, to show that {\it it} induces a $2k$-equivalence, as desired.

\bigskip
\noindent{\bf 1. \underbar{Functors with Smash Products}}

\bigskip
Roughly speaking, we want to define $U^n(R;M)$ to be
a simplicial spectrum
with $C_n = {\bf Z}/n{\bf Z}$ action
obtained by taking various smash products of the
Eilenberg-Mac Lane spectra of $R$ and $M$.
What the homotopy type of a smash product of spectra should be
has long been understood and there are several
different models for these which are equivalent
as spectra. However, the iteration of these constructions
is generally only associative up to homotopy. For simplicial
constructions this is not suitable since
it is not sufficient to have the simplicial identities
only up to homotopy. The solution to this problem
that we will follow is by M. B\"okstedt ([B]).
His construction involves the notion of a functor with smash
product (FSP).
We will actually be needing a straightforward generalization
of B\"okstedt's original construction suitable for
categories which was developed in [DMc2].

\bigskip\noindent
Let $\sc S_*$ denote the category of pointed simplicial sets.

\bigskip
\noindent{\caps Definition 1.1}:
A {\it functor with stabilization
} is a functor $F$ from $\sc S_*$ to $\sc S_*$
together with a natural transformation
$$\lambda_{X,Y}: X\wedge F(Y)\mapright .F(X\wedge Y)$$
such that

\medskip
\item{i)}$\lambda_{S^0, X}: S^0\wedge F(X)\to F(S^0\wedge X) $
is the obvious isomorphism for all $X\in S_*$
\medskip

\item{ii)} $\lambda_{X,Y\wedge Z}\circ (id_X\wedge\lambda_{Y,Z})
  = \lambda_{X\wedge Y,Z}$ for all $X,Y,Z\in S_*$.

\medskip
\item{iii)} If $X$ is $n$--connected, then $F(X)$ is $n$--connected.

\medskip
\item{iv)} Let $\sigma_X:F(X)\mapright .\Omega F(\Sigma X)$ be the
adjoint to $\lambda_{S^1,X}$. Then the following limit system
stabilizes for each $n$:
$$\pi_n |F(X)|\mapright \sigma_X.\pi_n\Omega| F(\Sigma X)|
\mapright \sigma_{\Sigma X}.\pi_n\Omega^2|F(\Sigma^2X)|\mapright
.\cdots$$

\bigskip\noindent
{\it Apology:} Our definition of a functor with stabilization is not as
general as some would  like as they will always be connective. The
above definition is suitable for the applications here and we
apologize to
anyone who will need to make the many technical modifications if
they need to use related results for non-connective FSP's.

\bigskip\noindent{\caps Definition 1.2}:
Let $\sc O$ be a set. A {\it functor with stabilization over
$\sc O$} is a functor $F$
from $\sc S_*\times\sc O\times\sc O$ to $\sc S_*$,
such that for all $A,B\in
\obj(\sc O)$,
$F_{A,B}(\ )$ is a functor
with stabilization, where we use the notations
$F_{A,B}(X)$ and $F(A,B)(X)$ as alternative notations for $F(X,A,B)$.

For $F$ a functor with stabilization, we let
${\bf F}$ be the spectrum with
${\bf F}(m) = F(S^m)$ and structure maps given
by $\sigma_{S^m}$ of (iii). We call
${\bf F}$ the {\it spectrum associated to $F$}.
We let
$\pi_i(F) = \pi_i{\bf F} = \lim_{n\to\infty}\pi_i\Omega^nF(S^n)$.
We say that $F$ is {\it $n$-connected} if
$\pi_i(F) = 0$ for $i\leq n$. Thus by condition (ii),
every functor with stabilization is $-1$-connected and hence
bounded below.
A functor with stabilization over $\sc A$ is $n$-connected if for
every $A,B\in \sc A$, $F_{A,B}$ is $n$-connected.

For $F$ a functor with stabilization over $\sc O$ and $F'$ a functor with stabilization over $\sc O'$,
a morphism $\eta:F\to F'$ is a set map $\sc O\to\sc O'$ and natural transformations of functors with stabilization $F_{A,B}\mapright\eta_{A,B}. F'_{\eta(A), \eta(B)}$ for all $A,B\in\sc O$.

\bigskip\noindent{\caps Definition 1.3}:
A {\it functor with smash product over
$\sc O$} (or just FSP) is a functor $F$
with stabilization over $\sc O$
together with natural transformations for all $A,B,C\in \sc O$:
 $$\eqalign{1_{A;X}:& X\mapright .F_{A,A}(X)\cr
 \mu_{A,B,C;X,Y}&: F_{B,C}(X)\wedge F_{A,B}(Y)\mapright
 .F_{A,C}(X\wedge Y)}$$
such that
$$\eqalign{\mu(\mu\wedge id) =& \mu(id\wedge\mu)\cr
\mu(1_{A;X}\wedge 1_{A;Y}) =& 1_{A;X\wedge Y}\cr
\lambda_{A,B;X,Y} =& \mu_{A,A,B;X,Y}(1_{A;X}\wedge id_{F_{A,B}(Y)})\cr
\rho_{A,B;X,Y} =& \mu_{A,A,B;X,Y}\circ (id_{F_{A,B}(X)}\wedge
1_{A;Y})\cr}$$
For $F$ and $F'$ FSP's over $\sc O$ and $\sc O'$ respectively, a
morphism $\eta$ from $F$ to $F'$ as FSP's is a
set map $\tilde\eta:\sc O\rightarrow \sc O'$ and
morphisms
of functors with stabilizations
$F_{A,B}\mapright\eta_{A,B}.F'_{\tilde\eta(A),\tilde\eta(B)}$
which strictly commute with
the natural transformations $\mu$ and $\mu'$.
We have \underbar{not} assumed a morphism of FSP's preserves
the unit. We will say a morphism is {\it unital} if
it does.

\bigskip\noindent
{\it Examples:}
For the applications in these notes, we will primarily
be interested in the following type of FSP.
Let $\sc A$ be a {\it linear} category
(its Hom sets are abelian groups and composition is
bilinear). For any two objects $A,B\in \obj(\sc A$),
we define the FSP  $\underline{\sc A}$
by
$$(X,A,B)\mapsto \Hom_{\sc A}(A,B)\otimes_{\bf Z}\tilde{\bf Z}[X]$$
where $\tilde{\bf Z}[X] = {\bf Z}[X]/{\bf Z}[*]$.
The
multiplication is given by sending smash to tensor followed
by composition:
$$\eqalign{
\sc A(B,C)(X)\wedge \sc A(A,B)(Y)
&\rightarrow (\sc A(B,C)\otimes_{\bf Z}\sc A(A,B))
\otimes_{\bf Z}\tilde{\bf Z}[X\wedge Y]\cr
&\rightarrow \sc A(A,C)(X\wedge Y)\cr
}$$
and the unit at any $A\in\sc A$ is given by the inclusion
$$X\mapright .\sc A(A,A)(X)$$
$$x\mapsto id_A\otimes 1\cdot x$$

For $F$ an FSP over $\sc O$,
we can form a linear category
with objects the set $\sc O$ and whose
Hom sets are $\pi_0F_{A,B}$. In this
way every small linear category can be
thought of as a special case of an FSP.

For $\sc B$ any category, we can form an
FSP over $\obj(\sc B)$ by
$$(X,A,B)\mapsto Hom_{\sc B}(A,B)_+\wedge X$$
(where $\ _+$ denotes a disjoint basepoint).
Conversely, given $F$ an FSP over $\sc O$ we can
form a category
with objects the set $\sc O$ whose
Hom sets are
the set (obtained by forgetting the topology)
of $F_{A,B}(S^0)$.

\bigskip\noindent{\caps Definition 1.4}: Let $F$ be an FSP over $\sc O$
and
$T$ a functor with stabilization over $\sc O$. A {\it structure
of left $F$--module} on $T$ is a natural transformation
$$l_{A,B,C;X,Y}:F_{B,C}(X)\wedge T_{A,B}(Y)\mapright .T_{A,C}(X\wedge
Y)$$
such that
$$\eqalign{l(\mu\wedge id) =& l(id\wedge l)\cr
\lambda_{A,B;X,Y} =& l_{A,B,B;X,Y}(1_{B;X}\wedge
id_{T_{A,B}(Y)})\cr}$$
The notion of right $F$--module is defined similarly.

\bigskip\noindent
{\it Example}: If we write $id$ for the identity FSP
(sending $X$ to $X$ at any $A\in\sc O$), then the category of functors with stabilization over $\sc O$
is isomorphic to the category of right (left) $id$--modules.

\bigskip\noindent{\caps Definition 1.5}:
A {\it bimodule} over $F$ is a functor
$T$ with stabilization over $\sc O$ together with a structure
of left and right module over $F$ such that
$$l_{A,C,D;X,Y\wedge Z}(id_{F_{C,D}(X)}\wedge r_{A,B,C;Y,Z}) =
r_{A,B,D;X\wedge Y,Z}(l_{B,C,D;X,Y}\wedge id_{F_{A,B}(Z)})$$
where $r$ is the structure of right module over $T$.

For $F$ and $F'$ FSP's  over $\sc O$ and $\sc O'$,  respectively, and
$P$ and $P'$ bimodules of $F$ and $F'$,
a {\it map} $(f;g)$ from $(F;P)$ to $(F';P')$ is
a pair of morphisms such that:

\item{(i)}
$f:F\rightarrow F'$ is a unital map of FSP's

\item{(ii)}
$g:P\rightarrow P'$ is a map of functors with stabilization over the same set map $\sc O\to\sc O'$ as $f$,
which is a map of $F$--bimodules
if we use $f$ to make $P'$  into an $F$--bimodule

\bigskip\noindent
{\it Examples:}
\item{(i)} Let $\sc A$ be a linear category and $T$ any bilinear functor
from $\sc A^{op}\times \sc A$ to abelian groups. We can form the
bimodule $\underline T$ of $\underline \sc A$ by
$$(A,B,X)\mapsto T(A,B)\otimes_{\bf Z}\tilde{\bf Z}[X].$$
\item{(ii)} In particular, if $G_1$ and $G_2$ are two functors of linear
categories $\sc A\to \sc B$ , we
can form the $\underline \sc A$-bimodule
$$(A,B,X)\mapsto \sc B(G_1(A),G_2(B))\otimes_{\bf Z}\tilde{\bf Z}[X].$$
\item{(iii)}
If $F$ is an FSP over $\sc O$, $P$ an $F$-bimodule, and $Y_.$ a finite pointed simplicial
set, we can construct another $F$-bimodule $P\otimes Y$ by letting
$$(P\otimes Y)(A,B)(X)=\vert P(A,B)(X) \wedge Y_.\vert$$
with $F$ acting through its action on $P$; all we are doing here is taking a smash product
with the realization of $Y_.$.  
\item{(iv)} In the same setting as (iii), if we know that for every $A,B\in \sc O$ and $X\in\sc S_*$,
$F(A,B)(X)$ and $P(A,B)(X)$ have abelian group structures compatible with all the FSP and
bimodule structure maps, then we can construct yet another $F$-bimodule $\tilde P[Y]$
by letting
$$\tilde P[Y](A,B)(X)= P(A,B)(X)\otimes_{\bf Z}\tilde{\bf Z}[Y].$$
The $F$-bimodule structure, as before, involves only $P$.  In degree $n$, $(P\otimes Y)(A,B)(X)$
has $P(A,B)(X) \wedge Y_n$, and $\tilde P[Y](A,B)(X)$ has $P(A,B)(X))\otimes_{\bf Z}\tilde{\bf Z}[Y_n]$; thus there is an inclusion map $P\otimes Y\hookrightarrow \tilde P[Y]$.  Note that it
induces a stable equivalence on the associated spectra (stabilizing, of course, in the $X$ coordinate.  The $Y$ coordinate is part of the definition of the bimodule). 

\bigskip\goodbreak\noindent
{\bf  2. \underbar {The construction of $U^n$}}

\bigskip\noindent{\it Notation}:
Let $I$ be the category whose objects are the natural numbers
considered
as ordered sets
($n = \{1<2<\cdots<n\}$)
and whose morphisms are all injective maps. For any $X\in I$ we
denote by $|X|$ the cardinality of $X$ and for any
$\underline{X} = (X_0,\ldots,X_j)\in I^{j+1}$ we let
$\sqcup\underline{X}$ denote
$X_0\sqcup X_1\sqcup\ldots\sqcup X_j$, where $\sqcup$ means {\it
concatenation}.

\bigskip
For $E$ a functor from the small category $\sc C$ to
pointed spaces, we let
$\hoco_{C\in\sc C}E(C)$ be a functorial choice for
constructing the homotopy colimit of the functor $E$.
We recall
the following lemma of B\"okstedt (see [B], [M]).

\bigskip\noindent
{\caps Lemma:} Let {\bf N} be any subcategory of $I$ with the
same set of objects but with exactly one morphism between
any pair of objects $n$ and $m$ where $n\leq m$. Let $G$ be a functor from
$I^{j+1}$ to spaces. If the connectivity of the maps
$G(n_0,\ldots,n_j)\rightarrow G(m_0,\ldots,m_j)$ for
maps in $I^{j+1}$ tends to infinity uniformly with $\Sigma n_i$,
then the inclusion $\hoco_{{\bf N}^{j+1}}G\rightarrow
\hoco_{I^{j+1}}G$ is a (weak) homotopy equivalence.

\bigskip\noindent
{\it Example 2.1:}
We note that there is a functor from $I$ to pointed spaces
given by sending  the ordered finite set
$X$ to $S^X$--the $|X|$--sphere obtained by
smashing together copies of $S^1$ indexed by the elements
of $X$. Given any functor with stabilization $F$, we can
define a functor from $I$ to pointed spaces by
sending the ordered set $X$ to $Map(S^X,F(S^X))$.
Given an injective map $X\mapright\alpha.Y$, let $\beta$
be some isomorphism of $Y$ such that $\alpha = \beta\circ inc$
where $inc$ is the ordered inclusion of $X$ into the
first $|X|$-terms of $Y$. The map $\alpha_*$ from
$Map(S^X,F(S^X))$ to $Map(S^Y,F(S^Y))$ is given by taking
$f\in Map(S^X,F(S^X))$ to the composite
\def\backslash{\setminus}
$$\matrix{
S^Y& &\mapright\alpha_*(f).&  &F(S^Y)  \cr
\mapdown\beta^{-1}.& & & &\mapup F(\beta).\cr
S^X\wedge S^{Y\backslash X}&
\mapright f\wedge id.&F(S^X)\wedge S^{Y\backslash X}&
\mapright \lambda.&F(S^X\wedge S^{Y\backslash X})\cr
}$$

\bigskip
\noindent{\it Notation:}
Let $\vec k = (k_1,\ldots,k_n)$
where the $k_i$ are nonnegative integers.
We set
$$I^{\vec k+n} = I^{k_1+1}\times\cdots\times I^{k_n+1}$$
and for $\underline {X}\in I^{\vec k+n}$
and $(i,j)$ integers such that $0\leq j\leq k_i$
we let
$X_{i,j}\in I$ be
in the $j+1$-st position of the $I^{k_i}$
component of $\underline{X}$.
We let
$$\sc A^{\vec k+n} = \sc A^{k_1+1}\times\cdots\times \sc A^{k_n+1}$$
and for $\underline{A}\in \sc A^{\vec k+n}$
and $(i,j)$ integers such that $0\leq j\leq k_i$
we let
$A_{i,j}$ be in the $j+1$-st position of the $\sc A^{k_i+1}$
component of $\underline{A}$.

\bigskip\noindent
{\caps Definition 2.2:}
Let $F$ be an FSP over a skeletally small category $\sc A$
and  let $P(1),\ldots,P(n)$ be a sequence of
$F$-bimodules.
We define the functor $V$ from $I^{\vec k+n}\times \sc A^{\vec k+n}$
to spaces by setting $V(\underline X;\underline A)$ to be:
$$\matrix{
P(1)(A_{1,1},A_{1,0})(S^{X_{1,0}})\wedge F(A_{1,2},A_{1,1})
(S^{X_{1,1}})
\wedge\ldots\wedge F(A_{2,0},A_{1,k_1})(S^{X_{1,k_1}})\wedge\cr
P(2)(A_{2,1},A_{2,0})(S^{X_{2,0}})\wedge
F(A_{2,2},A_{2,1})(S^{X_{2,1}})
\wedge\ldots\wedge F(A_{3,0},A_{2,k_2})(S^{X_{2,k_2}})\wedge\cr
\vdots\ \ \ \ \ \ \vdots\cr
P(n)(A_{n,1},A_{n,0})(S^{X_{n,0}})\wedge
F(A_{n,2},A_{n,1})(S^{X_{n,1}})
\wedge\cdots\wedge F(A_{1,0},A_{n,k_n})(S^{X_{n,k_n}}).\cr
}$$
Also, we set
$$
V(\underline X, \sc A)=
\bigvee_{\underline A\in \sc A^{\vec k+n}}
V(\underline X;\underline A).$$
If $\sc A$ is not itself small, we use $\sc A$'s skeleton rather than $\sc A$ itself here.  In light of example 6.3 below, which says that for a small category $\sc A$,  $U^n$ taken over a subcategory which is equivalent to all of $\sc A$ is $C_n$-equivariantly homotopy equivalent to $U^n$ taken over $\sc A$, this should not lead to problems when we map between small and large categories.
\bigskip
\noindent{\it Notation:}
Given any two pointed spaces $S$ and $T$, we write
$\underline{Map}(S,T)$ for the functor with stabilization
given by sending $X$ to $Map(S,X\wedge T)$.

\bigskip
\noindent{\caps Definition 2.3:}
Let $F$ be an FSP over a category $\sc A$
and  let $P(1),\ldots,P(n)$ be a sequence of
$F$-bimodules.
We use the above construction of $\underline{Map}$ to define
an $n$-simplicial
functor with stabilization
$U^{n}(F;P(1),\ldots,P(n))$, which is given
in simplicial dimension $\vec k = (k_1,\ldots,k_n)$ by:
$$\hoco _{\underline{X}\in I^{\vec k+n}}
G_{\vec k}(\underline X),$$
where
$$
G_{\vec k}(\underline X)=
\underline{Map}\left( (S^{\sqcup\underline{X}}),
\bigvee_{\underline A\in \sc A^{\vec k+n}}
V(\underline X;\underline A)\right)
=\underline{Map}\left( (S^{\sqcup\underline{X}}),
V(\underline X, \sc A)\right)
.$$
The stabilization maps for the homotopy colimit are given
like those in example 2.1.
The face operators are induced by the natural transformations
(let $d(i)_j$ be the $j$-th face map in the $i$-th simplicial
dimension):
$$d(i)_j:G_{(k_1,\ldots,k_n)} \mapright.
G_{(k_1,\ldots,k_{j}-1,\ldots,k_n)}\circ \partial(i)_j$$
where $\partial(i)_j: I^{(k_1,\ldots,k_n)}\rightarrow
I^{(k_1,\ldots,k_j-1,\ldots,k_n)}$ is the functor
$$\partial(i)_j(\underline X)
=\cases{
\pmatrix{X_{1,0},\ldots,X_{1,k_1}\cr
         \vdots\ \ \ \ \ \vdots\cr
        X_{i,0},\ldots,X_{i,j}\sqcup X_{i,j+1},\ldots,X_{i,k_i}\cr
        \vdots\ \ \ \ \ \vdots\cr
        X_{n,0},\ldots,X_{n,k_n}\cr}&if $0\leq j< n$\cr
\  &\  \cr
\ \ \pmatrix{X_{1,0},\ldots,X_{1,k_1}\cr
         \vdots\ \ \ \ \ \vdots\cr
        X_{i,0},\ldots,X_{i,k_{i-1}}\cr
        X_{i,k_i}\sqcup X_{i+1,0},\ldots,X_{i+1,k_{i+1}}\cr
        \vdots\ \ \ \ \ \vdots\cr
        X_{n,0},\ldots,X_{n,k_n}\cr}&if $j= n$ and $i\not= n$\cr
\ &\  \cr
\ \ \ \ \ \ \pmatrix{X_{n,k_n}\sqcup X_{1,0},\ldots,X_{1,k_1}\cr
        \vdots\ \ \ \ \ \vdots\cr
        X_{n,0},\ldots,X_{n,k_{n-1}}\cr}&if $j=n$ and $i=n$\cr
        }
        $$

and
$$d(i)_j(\underline X) = \cases{(1)&if $j=0$\cr
               (2)&if $0<j<k_i$\cr
               (3)&if $j=k_i$ and $i\not= n$\cr
               (4)&if $j=k_n$ and $i=n$\cr
               }$$
where $(1)$---$(4)$ are determined (in increasing order) by
$\underline{Map}(S^{\underline X},\ \ )$ of

{
\hsize = 6 in

$$\pmatrix{
P(1)(A_{1,1},A_{1,0})(S^{X_{1,0}})&\wedge F(A_{1,2},A_{1,1})
(S^{X_{1,1}})
\wedge\ldots\wedge F(A_{2,0},A_{1,k_1})(S^{X_{1,k_1}})\wedge\cr
         \vdots&\vdots\cr
r(i)&\wedge F(A_{i,3},A_{i,2})(S^{X_{i,2}})\wedge\ldots\wedge
F(A_{i+1,0},A_{i,k_{i}})(S^{X_{i,k_i}})\wedge\cr
         \vdots&\vdots\cr
P(n)(A_{n,1},A_{n,0})(S^{X_{n,0}})&\wedge
F(A_{n,2},A_{n,1})(S^{X_{n,1}})
\wedge\ldots\wedge F(A_{1,0},A_{n,k_n})(S^{X_{n,k_n}})\cr
}
\hfill$$

$$
\pmatrix{
P(1)(A_{1,1},A_{1,0})(S^{X_{1,0}})\wedge
&\ldots\ldots\ldots &\wedge
F(A_{2,0},A_{1,k_1})(S^{X_{1,k_1}})\wedge\cr
     \vdots&\vdots&\vdots\cr
P(i)(A_{i,1},A_{i,0})(S^{X_{i,0}})\wedge&
\ldots\wedge
F(A_{i,j},A_{i,j-1})(S^{X_{i,j-1}})\wedge \mu\wedge\ldots
&\wedge F(A_{i+1,0},A_{i,k_{i}})(S^{X_{i,k_i}})\wedge\cr
     \vdots&\vdots&\vdots\cr
P(n)(A_{n,1},A_{n,0})(S^{X_{n,0}})\wedge
&\ldots\ldots\ldots&\wedge F(A_{1,0},A_{n,k_n})(S^{X_{n,k_n}})\cr
}
$$

$$
\pmatrix{
P(1)(A_{1,1},A_{1,0})(S^{X_{1,0}})&
\!\wedge F(A_{1,2},A_{1,1})(S^{X_{1,1}})
\wedge\ldots\wedge F(A_{2,0},A_{1,k_1})(S^{X_{1,k_1}})\wedge\cr
\vdots&\vdots\cr
P(i)(A_{i,1},A_{i,0})(S^{X_{i,0}})&\!\wedge
F(A_{i,2},A_{i,1})(S^{X_{i,1}})
\wedge\ldots\wedge
F(A_{i,k_{i}},A_{i,k_{i}-1})(S^{X_{i,k_i-1}})\wedge\cr
l(i+1)&\!\wedge
F(A_{i+1,2},A_{i+1,1})(S^{X_{i+1,1}})\wedge\ldots\wedge
F(A_{i+2, 0},A_{i+1,k_{i+1}})(S^{X_{i+1,k_{i+1}}})\wedge\cr
       \vdots&\vdots\cr
P(n)(A_{n,1},A_{n,0})(S^{X_{n,0}})&\!\wedge
F(A_{n,2},A_{n,1})(S^{X_{n,1}})
\wedge\cdots\wedge F(A_{1,0},A_{n,k_n})(S^{X_{n,k_n}})\cr
}
$$

$$\pmatrix{
l(1)&\wedge F(A_{1,2},A_{1,1})
(S^{X_{1,1}})
\wedge\ldots\wedge F(A_{2,0},A_{1,k_1})(S^{X_{1,k_1}})\wedge\cr
     \vdots&\vdots\cr
P(n)(A_{n,1},A_{n,0})(S^{X_{n,0}})&\wedge
F(A_{n,2},A_{n,1})(S^{X_{n,1}})
\wedge\ldots\wedge F(A_{n,k_n},A_{n,k_n-1})(S^{X_{n,k_n-1}})\cr
}
$$

\par}

\noindent
where $r(j)$ and $l(j)$ are the right and left actions of $P(j)$,
$\mu$ is the product of $F$ and both (3) and (4) are after
composing with the necessary shuffle isomorphism of $I^{\vec k+n}$.
The degeneracy maps are defined similarly using the unit of $F$.

\bigskip
\noindent
\underbar{The construction of $\hat\otimes_F$}

We now define a construction which will
be useful but can be ignored by the reader until it is needed.
Given an FSP $F$ over a category $\sc A$, and
a left $F$-module $P$ and a right $F$-module $Q$, we
define a simplicial functor with stabilization over $\sc A$,
$P\hat\otimes_F Q$, as follows: Given $B,C\in\sc A$, we
define for each $k\in I$ a functor
$W_{B,C}$ from $I^{k+2}\times\sc A^{k+1}$ to
spaces by setting $W_{B,C}(\underline X;\underline A)$ to be:

$$
P(A_{1},B)(S^{X_{0}})\wedge F(A_{2},A_{1})
(S^{X_{1}})
\wedge\ldots\wedge F(A_{k+1},A_{k})(S^{X_{k}})
\wedge Q(C,A_{k+1})(S^{X_{k+1}})
$$

\bigskip\noindent
{\caps Definition 2.3:}
We define $P\hat\otimes_FQ$
to be the simplicial
functor with stabilization over $\sc A$ defined
on $B,C\in\sc A$ in simplicial dimension $[k]$ by:
$$\hoco _{\underline{X}\in I^{k+2}}
\underline{Map}\left( (S^{\sqcup\underline{X}}),
\bigvee_{\underline A\in \sc A^{k+1}}
W_{B,C}(\underline X;\underline A)\right)$$
The stabilization maps for the homotopy colimit are given
like those in example 2.1.

The face operators are induced by natural transformations
like those for $U^1(F;P)$. Rather then repeat all the notation,
we will describe these operators in words. We first note
that in each dimension $k$, for any $A_0,A_{k+2}\in\sc A$,
$P\hat\otimes_FF(A_0,A_{k+2})_{[k]}$
is a summand of
$U^1(F;P)_{[k+1]}$.
The restrictions of $d_0, d_1,\ldots, d_k$ and
$s_0, s_1,\ldots s_k$ to
$P\hat\otimes_FF(A_0,A_{k+2})_{[k]}$
send it to the summand $P\hat\otimes_FF(A_0,A_{k+2})_{[k-1]}$
inside $U^1(F;P)_{[k]}$,
and make
$P\hat\otimes_FF$ into a simplicial functor with stabilization over
$\sc A$.
Now we note that
we can replace $F$ in the role of a left module
over itself by an arbitrary left $F$-module $Q$ in
this construction, to obtain a simplicial object
that we define to be $P\hat\otimes_FQ$.

\bigskip\noindent
{\it Remarks 2.4}:
We now collect a few facts about $\hat\otimes_F$ whose
proofs are straightforward and left for the interested
reader. If $P$ and $Q$ are $F$-bimodules, then
$P\hat\otimes_FQ$ is naturally a simplicial $F$-bimodule.
Since the construction $\hat\otimes_F$ is natural, we
can iterate it to form
$P\hat\otimes_FQ\hat\otimes_FR$, a bi-simplicial $F$-bimodule.
If $f:P\rightarrow P'$ is a map of right $F$-bimodules
which is an equivalence, then
$f\hat\otimes_Fid_Q$ is an equivalence also (by the realization
lemma). Right multiplication induces a map of simplicial
right $F$-modules $P\hat\otimes_FF\rightarrow P$ (where
$P$ is the trivial multi-simplicial object with structure
maps all equal to the identity) which is an equivalence
(a simplicial contraction is given by the unused degeneracy
operator in $U^1(F;P)$).

\bigskip\noindent
{\caps Definition 2.5:} We define $\hat\otimes_F^n$ to
be the functor from $F$-bimodules to $n$--fold simplicial
$F$--bimodules given by:
$$\overbrace{P\hat\otimes_F\cdots\hat\otimes_FP}
^{n\rm\;times}.$$

\bigskip\noindent
{\caps Lemma 2.6:}
There is a
homotopy equivalence
 of
$n$--fold simplicial functors with stabilization
$$U^n(F;P(1),\ldots,P(n))\mapright.
U^{n-1}(F;P(1)\hat\otimes_FP(2),P(3),\ldots,P(n))
$$
By iteration of this equivalence, there is for
each $m|n$ an equivalence $\alpha_m$ of
$n$--fold simplicial functors with stabilization
$$U^n(F;P(1),\ldots,P(n))\mapright\alpha_m.
U^{n/m}(F;Q(1),\ldots,Q(n/m))
$$
where $Q(i) =
P(m(i-1)+1)\hat\otimes_F\cdots\hat\otimes_FP(m(i-1)+m).$

\bigskip
\noindent
{\it Proof}: This is a formal consequence of the definitions,
using
the fact that homotopy colimits commute with one another
and the fact that suspension
$$
\displaylines{\qquad
\hoco _{\underline{X}\in I^{\vec k+n}}
\underline{Map}\left( S^{\sqcup\underline{X}},
\bigvee_{\underline A\in \sc A^{\vec k+n}}
V(\underline X;\underline A)\right)
\hfill\cr
\hfill
\mapright ~.
\hoco _{\underline{X}\in I^{\vec k+n}}
\underline{Map}\left( S^{\sqcup\underline{X}}\wedge S^Y,
\bigvee_{\underline A\in \sc A^{\vec k+n}}
V(\underline X;\underline A)\wedge S^Y\right)
\qquad}
$$
induces a homotopy equivalence when evaluated at any space.

\bigskip
\noindent{\bf 3. \underbar{First properties of $U^n$}}

\bigskip
We now establish some elementary first properties of $U^n$.
Since for a fixed FSP $F$ over a category $\sc O$, a morphism of $F$-bimodules is just a map
of functors with stabilization at each $(A,B)\in\sc O\times\sc O$ which strictly commutes with the left and right actions of $F$, it follows immediately from the definitions that
$U^n(F;\ )$
is a functor from the $n$-fold product category of
$F$-bimodules to spectra. The $U^n$ are also functorial
in the FSP variable but the statement is slightly messy and
left to the interested reader.

\bigskip
\noindent
{\caps Lemma 3.1}: If
 $f(i):P(i)\rightarrow P'(i)$ is an $m$--connected
map of bimodules,
then the induced map of functors with stabilization
$$U^n(f):
U^n(F;\ldots,P(i),\ldots)\mapright.U^n(F;\ldots,P'(i),\ldots)$$
is (at least) $m$--connected.

\bigskip\noindent
{\it Proof:}
Since
the associated spectrum of
$U^n$ is an $\Omega$-spectra (the structure
maps $\Omega E_n\rightarrow E_{n+1}$ are homotopy equivalences)
it suffices to show $U^n(f)(S^0)$ is $m$--connected. A map
of simplicial spaces which is $m$--connected in each simplicial
dimension is $m$--connected upon realization (essentially
because homotopy colimits preserve connectivity) and hence it
suffices to show that $U^n(f)(0)_{[k]}$ is $m$--connected
for all $k\geq 0$. If $f:X\rightarrow X'$ is a $q$--connected
map of (pointed) spaces, then $id\wedge f: Y\wedge X\rightarrow
Y\wedge X'$
is $p+q+1$--connected for any $p$--connected space $Y$.
Assume for convenience that $i=1$.
Consider the homotopy equivalent rewriting
$$\hoco _{\underline{X}\in I^{\vec k+n}}
Map\ \left( S^{\sqcup\underline{X}},
\bigvee_{\underline A\in \sc A^{\vec k+n}}
V(\underline X;\underline A)\right)\simeq$$
$$
\hoco _{\underline{X'}\in
I^{k_1}\times I^{k_2+1}\times\cdots\times I^{k_n+1}}
\hoco _{X_{1,0}\in {I}}
Map\ \left( S^{X_{1,0}\sqcup\underline{X}'},
\bigvee_{\underline A\in \sc A^{\vec k+n}}
V(X_{1,0}\sqcup\underline X;\underline A)\right).$$
By our previous observations and the fact that our
functors with stabilization
take $n$--connected spaces to $n$--connected spaces,
we see that the map induced by $f(i)$ on
$$\hoco _{X_{1,0}\in {I}}
Map\ \left( S^{X_{1,0}\sqcup\underline {X}'};
\bigvee_{\underline A\in \sc A^{\vec k+n}}
V(X_{1,0}\sqcup\underline X;\underline A)\right)$$
will be $m$--connected. The result follows from the fact
that homotopy colimits preserve connectivity.

\bigskip
\noindent
{\caps Corollary 3.2}: The functor $U^n$ is a reduced homotopy functor
in each variable $i$ separately.

\bigskip\noindent
{\it Constructing $\Omega$-FSP's}:
It will often be useful to replace a given functor with stabilization
$F$ by an equivalent one (that is: an FSP whose associated spectrum is stably equivalent to that of the original FSP) whose associated spectrum is an
$\Omega$-spectrum. Define $\Omega^{\infty}F$ to be the
new functor with stabilization defined by
$$\Omega^{\infty}F = \hocolim_{X\in I} \underline{Map} (S^X,F(S^X))$$
with $\lambda_{Z,Y}$ defined by the natural composite
$$\eqalign{
Z\wedge \Omega^{\infty}F(Y) &=\ Z\wedge
                \hocolim_{X\in I}\ Map\ (S^X,F(Y\wedge S^X))\cr
  &\ \cong \hocolim_{X\in I}\ Z\wedge Map\ (S^X,F(Y\wedge S^X))\cr
  &\mapright.\hocolim_{X\in I}\ Map\ (S^X,Z\wedge F(Y\wedge S^X))\cr
  &\mapright \lambda.\hocolim_{X\in I}\ Map(S^X,F(Z\wedge Y\wedge S^X))\cr
  &=\ \Omega^{\infty}F(Z\wedge Y). \cr
  }$$
The natural map $F\rightarrow \Omega^{\infty}F$ of functors
with stabilization gives a stable
 equivalence of the associated spectra by condition 1.1 (iii).
If $F$ was an FSP, then we can make $\Omega^{\infty}F$
an FSP by defining $\mu$ to be the composite
$$\eqalign{\Omega^{\infty}F(Z) & \wedge \Omega^{\infty}F(Y)\cr
&= \hoco _{X\in I}Map\ (S^X,F(Z\wedge S^X)) \wedge
   \hoco _{X'\in I}Map\ (S^{X'},F(Y\wedge S^{X'}))\cr
&\mapright\alpha.\hocolim_{X,X'\in I\times I}\
Map\ (S^{X\sqcup X'},F(Z\wedge S^X)\wedge F(Y\wedge S^{X'}))\cr
&\mapright{\mu_F}.
\hocolim_{X,X'\in I\times I}\
Map\ (S^{X\sqcup X'},F(Z\wedge S^X\wedge Y\wedge S^{X'}))\cr
&\mapright\beta.
\hocolim_{X,X'\in I\times I}\
Map\ (S^{X\sqcup X'},F(Z \wedge Y\wedge S^{X\sqcup X'}))\cr
&\mapright\gamma.
\hocolim_{\tilde X\in I}\
Map\ (S^{\tilde X},F(Z \wedge Y\wedge S^{\tilde X}))\cr
&= \Omega^{\infty}F(Z\wedge Y)\cr
}$$
where $\alpha$ is obtained by smashing maps, $\beta$ by
switching factors and $\gamma$ is induced by the concatenation functor
$\sqcup:I\times I\rightarrow I$.
Checking that $\mu$ is strictly associative is straightforward
and left to the reader.
The natural map from $F$ to $\Omega^{\infty}F$ is
an equivalence of (unital) FSP's (which determines the
units in $\Omega^{\infty}F$).

If $P$ is a right/left/bi-module of $F$, then $\Omega^{\infty}P$
is again a right/left/bi-module of $\Omega^{\infty}F$
(defined as we did for FSP's above)
and the natural map $P\rightarrow \Omega^{\infty}P$
is a map of right/left/bi-modules.

\bigskip\noindent
{\caps Corollary 3.3}: The natural map of functors with
stabilization
$$U^n(F;P(1),\ldots,P(n))\rightarrow U^n(\Omega^{\infty}F;
\Omega^{\infty}P(1),\ldots,\Omega^{\infty}P(n))$$
is an equivalence. Thus, we can always replace the
FSP and its associated bimodules with equivalent ones
whose associated spectra are $\Omega$-spectra.

\bigskip\noindent
{\caps Lemma 3.4}:
If each bimodule $P(i)$ is $m_i$--connected,
then $U^n(F;P(1),\ldots,P(n))$ is (at least)
$\Sigma_{i=1}^nm_i + (n-1)$
connected.

\bigskip\noindent
{\it Proof:}
It suffices to prove the result after taking $\Omega^{\infty}$
of $F$ and the $P(i)$'s. Thus, we may assume that
$P(i)(S^X)$ is $m_i+|X|$--connected for all $i$ and $X\in I$.
Recall that if $X_i$ is $x_i$--connected for $1\leq i\leq n$,
then $X_1\wedge\cdots\wedge X_n$ is $\Sigma_{i=1}^n x_i + (n-1)$
connected.
Thus, $V(\underline{X},\underline A)$ will
be $|\underline{X}|+\Sigma_{i=1}^nm_i+(n-1)$
connected for all $\underline{X}$.
Since homotopy colimits preserve connectivity, we see
that
$U^n(F;P(1),\ldots,P(n))$ is (at least)
$\Sigma_{i=1}^nm_i + (n-1)$ connected.

\bigskip
For a continuous map of (pointed)
spaces $f:X\rightarrow Y$, we set:
$${\rm hofib\;}(f) = {\rm pullback\;}\pmatrix{ & &X\cr
                   & &\mapdown f.\cr
                   PY&\mapright ev_1.&Y\cr
                   }
$$
where $PY$ is the pointed path space of $Y$
($Map_*([0,1],Y)$)
and $ev_1$ is
evaluation of the path at $1$. Using the natural
homeomorphism $\Omega PY\cong P\Omega Y$, we
see that there is a natural homeomorphism
$\Omega {\rm hofib\;}(f)\cong {\rm hofib\;}(\Omega f)$.
Given a map $f:{\bf E}\rightarrow {\bf E}'$ of spectra, we let
${\rm hofib\;}(f)$ be the spectrum defined by
${\rm hofib\;}(f)(n) = {\rm hofib\;}(f(n))$ with structure maps
given by the composite
$$\Omega {\rm hofib\;}(f)(n)\cong {\rm hofib\;}(\Omega f(n))
\mapright. {\rm hofib\;} (f(n+1)).$$
We call a sequence of maps (of spectra)
${\bf E}''\mapright f''.{\bf E}\mapright f'.{\bf E}'$
a {\it fibration} if
$f'\circ f'' = *$ and the map induced by $f''$ from
${\bf E}''$ to ${\rm hofib\;}(f')$ is an equivalence.

\bigskip
\noindent
{\caps Definition 3.5}: We call a sequence of maps of functors with
stabilization
$$F''\mapright\alpha.F\mapright\beta.F'$$
a {\it fibration}
if the associated map of underlying spectra is a fibration. We will
call a sequence of maps of bimodules {\it short exact} if
it is a fibration as a sequence of functors with stabilization.

\bigskip\noindent
{\it Example}: A sequence of abelian groups
$A\mapright f.B\mapright g.C$ is short exact if and only if
the associated sequence of functors with stabilization
(given by $A\mapsto A\otimes_{\bf Z}\tilde{\bf Z}[\ ]$ {\it etc.})
is short exact. If one considers a sequence of simplicial
abelian groups
$A_*\mapright f_*.B_*\mapright g_*.C_*$, then the associated
sequence of functors with stabilization is short exact if
and only if the sequence of simplicial abelian groups
is quasi-exact (the associated chain complexes
and maps form
a triple in the associated triangulated category).

\bigskip\noindent
{\caps Lemma 3.6}: If
$$P(i)''\mapright f''.P(i)\mapright f'.P(i)' $$
is a short exact sequence of bimodules, then
the induced sequence of
functors with stabilization
$$U^n(F;\ldots,P(i)'',\ldots)
\mapright.U^n(F;\ldots,P(i),\ldots)
\mapright.
U^n(F;\ldots,P(i)',\ldots)$$
is a fibration.

\bigskip\noindent
{\it Proof:}
First note that it suffices to show that
it suffices to
know that the map to the hofib computed degreewise
(fixing the simplicial directions)
is an equivalence.
For
$$\underline {X}'\in
I^{k_1+1}\times\cdots\times I^{k_i}\times\cdots\times I^{k_n+1}$$
define
$$
G''(X') =
\hoco _{X\in I}
Map\ \left( (S^{X\sqcup\underline{X}'}),
\bigvee_{\underline A\in \sc A^{\vec k+n}}
S^2\wedge V(X\sqcup\underline X';\underline A)\right)$$
where we are using the bimodule $P''(i)$. Similarly
define $G(X')$ (respectively $G'(X')$)
using the bimodule $P(i)$ (respectively $P'(i)$).
Since $\Omega$ commutes with taking hofib and
homotopy colimits commute with finite homotopy
pull-backs (in spectra), it suffices to
show that the map
$$G''(\underline {X}')\mapright.{\rm hofib}[
G(\underline {X}')\mapright.
G'(\underline {X}')]$$
is a homotopy equivalence
and this follows by an argument similar to
that used in lemma 3.1.

\bigskip
Let $(f;g):(F;P)\rightarrow (F';P')$ be a morphism of
pairs of an FSP and a bimodule over it (1.5). Thus, for all $n\geq 1$ we
similarly
obtain a morphism of pairs (techincally of
simplicial $F$--bimodules):
$$(f;\hat\otimes^n_Fg):(F;\hat\otimes_F^nP)\rightarrow
(F';\hat\otimes_{F'}^nP')$$
(see 2.5).

\bigskip\noindent
{\caps Lemma 3.7:} Let $(f;g)$ be a morphism of pairs such that
the morphism
$U^1(f;\hat\otimes_F^ng)$ is a weak equivalence, then
the map $U^n(f;g)$ is also. In particular, if $P=F$ then whenever
$U^1(f;g)$ an equivalence,
$U^n(f;g)$ must be an equivalence too for all $n\geq 1$
(this last statement also follows using edge-wise subdivision,
see [BHM]).

\bigskip\noindent
{\it Proof}: The first statement follows from Lemma 2.6. The second
part follows from the first and the fact that
$\otimes_F^nF\mapright\simeq.F$ as $F$--bimodules (see Remarks 2.4).

\bigskip\noindent
{\bf 4.\underbar{The construction of $W_n$}}

\bigskip
In this section, we will only be considering
$U^n(F;P(1),\ldots,P(n))$ as a simplicial
functor with stabilization
by taking the diagonal of the $n$--dimensional
simplicial functor with stabilization defined in 2.3. We first note
that there is an isomorphism
$$t: U^n(F;P(1),\ldots,P(n))\cong U^n(F;P(n),P(1),\ldots,P(n-1))$$
induced in each simplicial dimension $[k]$ by precomposing each map
with $\tau^{-1}$,
%
where $\tau: I^{(k,\ldots,k)}\rightarrow
I^{(k,\ldots,k)}$ is the functor
$$\tau(\underline X)
=
\pmatrix{X_{n,0},\ldots,X_{n,k}\cr
        X_{1,0},\ldots,X_{1,k}\cr
        \vdots\ \ \ \ \ \vdots\cr
        X_{n-1,0},\ldots,X_{n-1,k}\cr},
$$
and
post-composing with the analogous permutation

{
\hsize = 6 in

$$ V(\underline X;\underline A) \to
\pmatrix{
P(n)(A_{n,1},A_{n,0})(S^{X_{n,0}})
\wedge\ldots\wedge F(A_{1,0},A_{n,k})(S^{X_{n,k}})\cr
P(1)(A_{1,1},A_{1,0})(S^{X_{1,0}})
\wedge\ldots\wedge F(A_{2,0},A_{1,k})(S^{X_{1,k}})\cr
         \vdots\ \ \ \ \ \ \vdots\cr
P(n-1)(A_{n-1,1},A_{n-1,0})(S^{X_{n-1,0}})
\wedge\ldots\wedge F(A_{n,0},A_{n-1,k_n})(S^{X_{n-1,k}})\cr
}
\hfill$$
\par}

\bigskip
It is sometimes convenient to work with a slightly modified
version of $U^n$ we write as $\tilde U^n$.
The only difference between $\tilde U^n$ and $U^n$ is
that for $\tilde U^n$ one takes the homotopy colimit over
the {\it diagonal} of $(I^{k+1})^{\times n}$ in simplicial
dimension $[k]$. The face and degeneracy maps are easily
seen to restrict to this sub-limit system as well as
the simplicial isomorphism $t$. The natural map
of simplicial functors with stabilization from $\tilde U^n$ to $U^n$
is
determined by the inclusion of subcategory
is always an equivalence (by finality and the realization lemma).

\bigskip\noindent
{\caps Definition 4.1}: For $F$ an FSP and $P$ a bimodule, we
define $U^n(F;P)$ to be
the simplicial
functor with stabilization
with $C_n$--action (given by $t$)
$\tilde U^n(F;P,\ldots,P)$. Thus, $U^n$ is a functor from
the category of
pairs $(F;P)$ (an FSP $F$ and $F$--bimodule $P$) with
morphisms of pairs (1.5) to simplicial functors
with stabilization with $C_n$--action.

\bigskip\noindent
{\it Remark:} $U^1(F;P)$ is just $THH(F;P)$: the topological
Hochschild homology of $F$ with coefficients in $P$ as defined
in [DMc2] which was a straightforward generalization of the
definition found in [PW] to FSP's with several objects.
The spectrum $\tilde U^n(F;F,\ldots,F)$ is isomorphic (as simplicial
functors with stabilization
with $C_n$--action) to the $n$-th edgewise
subdivision of $THH(F;F)$.

\bigskip\noindent
{\caps Lemma 4.2:} For every $m|n$, the
$n$--fold
simplicial isomorphism
$$\alpha_m:U^n(F;P,\ldots,P)\rightarrow U^{n/m}
(F;P^{\hat\otimes m},\ldots,P^{\hat\otimes m})$$
of 2.6 is $C_{n/m}$--equivariant if we take
diagonals. The composite of
$\tilde U^n(F;P,\ldots,P)\mapright\simeq.U^n(F;P)$
with the diagonal of $\alpha_m$ factors to produce
a $C_{n/m}$--equivariant map:

$$\matrix{\tilde U^n(F;P,\ldots,P)&\mapright\tilde\alpha_m.&
   \tilde U^{m/n}(F;P^{\hat\otimes^m_F},\ldots,P^{\hat\otimes_F^m})\cr
   \mapdown =.& &\mapdown =.\cr
   U^n(F;P)&\mapright\alpha_m.&U^{n/m}(F;P^{\hat\otimes_F^m})\cr
}$$

\noindent
which is an equivalence and natural with respect to morphism pairs.

\bigskip\noindent
{\it Proof}: This is formally true from the definitions and
left to the reader.

\bigskip
Given a simplicial $C_n$--set $X_*$, the $C_n$ fixed
point space of the realization of $X_*$ is homeomorphic
to the realization of the simplicial set
$X_*^{C_n}$ obtained by taking fixed points degreewise.
It follows that the same is true for a simplicial
$C_n$--CW complex (or anything $C_n$-equivariantly homotopy
equivalent to one in each degree) and hence one can compute
$U^n(F;P)(X)^{C_m}$ degreewise for every $m|n$.
Similarly, since our model for  homotopy colimits
is given by simplicial spaces, we see that
if $E$ is a functor from the small category $\sc C$
to $C_n$--equivariant CW--complexes (or spaces $C_n$-equivariantly equivalent
to them) then
$(\hoco _{C\in\sc C}E(C))^{C_n}\cong \hoco _{C\in\sc C}E(C)^{C_n}$.
Thus, $U^n(F;P)^{C_m}$ for $m|n$ is naturally homeomorphic
to the realization of
$$[k]\mapsto U^n(F;P)^{C_m}_k\cong
\hoco _{\underline {X}\in I^{k+1}}\underline{Map\ }
\left((S^{\sqcup \underline {X}})^{\wedge n},\bigvee_{\underline
{A}\in
(\sc A^{k+1})^{\times n}}V(\underline{X}^{\times n};
\underline{A})\right)^{C_m}.$$
We note further, that with the specified $C_n$-action,
we have natural $C_n/C_m\cong C_{{n\over m}}$ equivariant
homeomorphisms

$$((S^{\sqcup \underline{X}})^{\wedge n})^{C_m}
\cong (S^{\sqcup \underline{X}})^{\wedge ({n\over m})}\leqno(4.2)$$
$$\left(\bigvee_{\underline {A}\in
(\sc A^{k+1})^{\times n}}V(\underline{X}^{\times n};
\underline{A})\right)^{C_m}
\cong
\left(\bigvee_{\underline {A}\in
(\sc A^{k+1})^{\times ({n\over m})}}V(\underline{X}^{\times ({n\over
m})};
\underline{A})\right).$$

We recall that if $G$ is a group with normal subgroup $H$,
$X$ and $Y$ $G$--spaces, then there
is a continuous map from $Hom(X,Y)^{G}$ to $Hom(X^H,Y^H)^{G/H}$
given by the restriction to fixed point subspaces. In fact, it
is the composite:
$$Hom(X,Y)^G\mapright res^H.Hom(X^H,Y)^G = Hom(X^H,Y^H)^{G/H}\leqno(4.3)$$
By the sequence $X^H\mapright .X\mapright .X/X^H$, which is a cofibration since
$H$ acts simplicially.
We see that the fiber the map of $res^H$
(and hence of the composite) is just $Hom(X/X^H,Y)^G$.

\bigskip\noindent{\caps Definition 4.3}:
For $r$, $s$ and $t$ integers greater
than $0$, we let
$$Res^r:
U^{rst}(F;P)^{C_{rs}} \to U^{st}(F;P)^{C_s}$$ be the
map of simplicial functors with stabilization
with $C_t$--action
defined
degreewise by applying $res^{C_r}$
and making the appropriate identifications by (4.2).
Thus,
$$Res^1 = id$$
and
$$Res^r Res^{s} = Res^{r\cdot s} = Res^{s}Res^r.$$

\bigskip\noindent
{\caps Lemma 4.4.}
For $r$, $s$ and $t$ integers greater
than $0$ the following diagram commutes
$$\matrix{
U^{rst}(F;P)^{C_{rs}}&\mapright Res^r.&U^{st}(F;P)^{C_{s}}\cr
\mapdown \alpha_t.&&\mapdown\alpha_t.\cr
U^{rs}(F;P^{\hat\otimes t})^{C_{rs}}&\mapright Res^r.&
U^s(F,P^{\hat\otimes t})^{C_s}.\cr}
$$

\bigskip\noindent
{\it Proof.}
As in the previous discussion, we check the claim simplicially.
$$ U^{rst} (F;P)^{C_{rs}}_k\cong
\hoco _{\underline {X}\in I^{k+1}}\underline{Map\ }
\left((S^{\sqcup \underline {X}})^{\wedge {rst} },\bigvee_{\underline
{A}\in
(\sc A^{k+1})^{\times {rst} }}V(\underline{X}^{\times {rst} };
\underline{A})\right)^{C_{rs}}.$$
We consider the model of $P^{\hat\otimes t}$ arising from the diagonal of
the $(t-1)$-simplicial construction $(((P\hat\otimes_F P)\hat\otimes_F P)
\hat\otimes_F
\cdots)\hat\otimes_F P$.  Also, we consider the simplicial
model of $U^{rs}(F;P^{\hat\otimes t})$ which is the diagonal on the simplicial
structure of $U$ and all the simplicial structures of the $P^{\hat\otimes t}$
simultaneously.
Then $\alpha_t$ is induced by grouping together the first $(t-1)(k+1)+1$
coordinates in each $t(k+1)$-tuple of coordinates, and sending them to
the corresponding coordinate in $P^{\hat\otimes t}$.

Therefore, if instead of looking at  $C_{rs}$-equivariant maps on the full 
$ (S^{\sqcup \underline {X}})^{\wedge {rst}} $
we look at their restrictions to the $C_r$-fixedpoints in the domain, which
must land in that part of the range whose coordinates repeat themselves in
$r$ blocks of $ts(k+1)$, the effect of grouping together before or after the
$r$-fold repetition is the same.

\bigskip\noindent
{\caps Definition 4.5.}
For $r$ and $s$ positive integers greater than $0$, we
let
$$I^s: U^{rs}(F;P)^{C_{rs}}\mapright. U^{rs}(F,P)^{C_r}$$
be the obviouis inclusion of the fixed set of a group in the fixed set of
its subgroup.

\bigskip\noindent
{\caps Lemma 4.6.} For $r$, $s$ and $t$ integers greater
than $0$ the following diagram commutes
$$\matrix{
U^{rst}(F;P)^{C_{rs}}&\mapright Res^r.&U^{st}(F;P)^{C_{s}}\cr
\mapdown I^s.&&\mapdown I^s.\cr
U^{rst}(F;P)^{C_{r}}&\mapright Res^r.&
U^{st}(F,P).\cr}
$$

\bigskip\noindent
{\it Proof.}
This is immediate from the definition of $Res^r$: by equation (4.3), the
restriction  sends
$$Hom(X,Y)^{C_{rs}}
\mapright res^{C_r}.Hom(X^{C_r},Y^{C_r})^{C_{rs}/C_r} = Hom(X^{C_r},Y^{C_r})^{C_s}
.$$

\bigskip\noindent
{\caps Definition 4.7.} For $s$ a positive integer greater than $0$
we define
$$P^s = \alpha_s\circ I^s: U^{rst}(F;P)^{C_{rs}}\mapright.
U^{rt}(F;P^{\hat\otimes s})^{C_{r}}.$$

By lemmas 4.4 and 4.6 we see that whenever they are both defined,
$$Res^r P^s = P^s Res^r.\leqno(4.4)$$

\bigskip
Let ${\bf N}^{\times}$ be the natural numbers $\{1,2,\ldots\}$
as a partially ordered set with $n<m \Leftrightarrow m|n$.
For $F$ an FSP and $P$ an $F$ bimodule,
we have a   functor from ${\bf N}^{\times}$ to
functors with stabilization sending every natural number $n$ to 
$ U^n(F;P)^{C_n}$ and every morphism $n<m$ to  $Res^{{n\over m}}.$

\bigskip\noindent
{\caps Definition 4.8}: Let $F$ be an FSP and $P$ an $F$ bimodule.
For ${\sc M}$ a subcategory of ${\bf N}^{\times}$, we set
$$W_{\sc M}(F;P) = \holim_{{n\in \sc M}}\ U^n(F;P)^{C_n}.$$

We will use simplified notation
for various distinguished subcategories of ${\bf N}^{\times}$
as follows. First, we set
$$W(F;P) = W_{{\bf N}^{\times}}(F;P).$$
We let $\{\leq n\}$ be the  full subcategory of ${\bf
N}^{\times}$
generated by $\{1,2,\ldots,n\}$ and write
$$W_n(F;P) = W_{\{\leq n\}}(F;P).$$
We let $(p)$ be the full subcategory of ${\bf N}^{\times}$
generated by the powers of $p$, $(p) = \{1,p,p^2,\ldots\}$ and
write
$$W^{(p)}(F;P) = W_{(p)}(F;P)$$
(we use this notation so we can have room for a subscript).
Let $(\leq p^n)$ be the  full subcategory of
$(p)$ generated by $\{1,p,\ldots,p^n\}$ and write
$$W^{(p)}_n(F;P) = W_{(\leq p^n)}(F;P).$$

\bigskip\noindent
{\caps Definition 4.9}: For $n$ a natural number, we let
$$W_n(F;P)\mapright R_n. W_{n-1}(F;P)
\ \ \ \ \ \ W_n^{(p)}(F;P)\mapright R^{(p)}_n. W^{(p)}_{n-1}(F;P)$$
be the natural maps obtained by restriction to subcategories.
Thus,
$$W(F;P) = \holim_{\infty\leftarrow n}W_n(F;P)\ \ \ \ \ \
W^{(p)}(F;P) = \holim_{\infty\leftarrow n}W^{(p)}_n(F;P)$$
with structure maps given by the $R_n$'s.

\bigskip
\noindent
{\caps Lemma 4.10:} For $s$ a positive integer, $P^s$ induces
a natural map
$$P^s: W(F,P)\mapright. W(F,P^{\hat\otimes s}).$$
For $p$ a prime, we also have
$$P^p: W_n^{(p)}(F,P)\mapright. W_{n-1}^{(p)}(F,P^{\hat\otimes p})$$
which commutes with $R^{(p)}$.

\bigskip\noindent
{\it Proof.}
Since $P^s$ commutes with $Res^r$ whenever both are defined, $P^s$ induces a
map between the inverse system $U(F;P)$  on the subcategory 
$s\cdot {\bf N}^\times$ of multiples of $s$ and the inverse system
$U(F;P^{\hat\otimes s})$ on ${\bf N}^\times$.  Therefore we get a map
$$W_{s\cdot \sc M}(F;P) 
\to
W_{\sc M}(F;P^{\hat\otimes s})$$ 
for any subcategory $\sc M$ of 
${\bf N}^\times$.
But the subcategory $s\cdot {\bf N}^\times$ is coinitial in
${\bf N}^\times$, so the obvious inclusion of categories induces a
natural  isomorphism
$$W_{s\cdot {\bf N}^\times}(F;P)\cong W(F;P)$$
through which we can define the first map claimed in the lemma.

Since there are initial objects, we know that 
$W_n^{(p)}(F;M)=U^{p^n}(F;M)^{C_{p^n}}$ and
$W_{n-1}^{(p)}(F;M)=U^{p^{n-1}}(F;M)^{C_{p^{n-1}}}$
for $M=P$ and $M=P^{\hat\otimes p}$, and that $R^{(p)}_n$ in both systems is
simply $Res^p$; it commutes with $P^p$ by equation (4.4).

\bigskip
\noindent{\bf 5. \underbar{The fiber of $R_n$}}

\bigskip
Our goal in this section is to identify the fiber of
the maps $R_n$ and $R^{(p)}_n$ up to natural equivalence
with $U^n_{hC_n}$ and $U^{p^n}_{hC_{p^n}}$. This was
essentially done by T.~Goodwillie
in the appendix to
his MSRI notes [G]. Since these MSRI notes are not published,
in this section we reproduce what is needed from them
(5.3, 5.4 and 5.5 below) to establish
the result. We have modify some of the constructions
found in [G] to make the proofs more transparent.

\bigskip
Let $G$ be a group.
Recall that for $X$ a (pointed) space with $G$--action,
we define the {\it homotopy orbit} space of $X$ to
be $X_{hG} = X\wedge_G EG_+
= (X\wedge EG_+)_G$. We recall that if
$f:X\rightarrow Y$ is an $n$--connected $G$--equivariant
map then $f_{hG}$ is also $n$--connected. We note that
if $F$ is a functor with stabilization with $G$--action, then
$X\mapsto F(X)_{hG}$ is again naturally a functor
with stabilization.

We define the {\it homotopy fixed-point} space of $X$
to be $X^{hG} = Map_G(EG_+,X) = Map_*(EG_+,X)^G$. If
$f:X\rightarrow Y$ is a $G$--equivariant map which is
also an equivalence, then $f^{hG}$ is also an
equivalence but $(\ )^{hG}$ does \underbar{not}
preserve connectivity in general. Thus, if $F$ is
a functor with stabilization then
$X\mapsto F(X)^{hG}$ satisfies 1.1(i) but not necessarily
(ii), (iii), or (iv) and hence is not again a functor with stabilization.

\bigskip\noindent
{\caps Definition 5.1}: A {\it functor with structure} will
be a functor $F$ from
$\sc S_*$ to $\sc S_*$ together with a natural transformation
$$\lambda_{X,Y}: X\wedge F(Y) \mapright. F(X\wedge Y)$$
which satisfies 1.1 (i) but not necessarily 1.1(ii), 1.1(iii), or 1.1(iv).
A {\it functor with structure over $\sc O$} for a set $\sc O$ is
a functor from $\sc S_*\times\sc O\times\sc O$ to $\sc S_*$
such that for all $A,B\in\sc O$, $F_{A,B}(\ )=F(A,B)(\ )$ is
a functor with structure.

\bigskip\noindent
{\caps Definition 5.2}: Let $G$ be a group and $F$ a functor
with stabilization with $G$-action. We define
the {\it homotopy orbits} of $F$ to be the
functor with stabilization
$$F_{hG} :
X\mapsto { \hocolim_m}  \ \Omega^m [F(\Sigma^m X)]_{hG}
$$
and the {\it homotopy fixed--points} of $F$ to be
the functor with structure
$$F^{hG} :
X\mapsto { \hocolim_m }\ \Omega^m [{ \hocolim_\ell}\ \Omega^\ell
F(\Sigma^{m+\ell} X)
]^{hG}).$$

\goodbreak

\bigskip\noindent
{\it Important Remark:} If $F_*$ is a simplicial functor
with stabilization with $G$--action, then  the commuting diagram
$$
\matrix{
|(F_*)_{hG}|&
=&{\rm \hocolim}_{\Delta^{op}}[(F_*)_{hG}]\cr
\mapdown.& &\mapdown\cong.\cr
|(F_*)|_{hG}&
=& ({\rm \hocolim}_{\Delta^{op}}[F_*])_{hG}\cr
}$$
shows that homotopy orbits commute with realizations. However,
since $[q]\mapsto (F_q)^{hG}$ is only a simplicial functor with
structure, we do \underbar{not} get that
the natural map
$$|(F_*)^{hG}|\mapright.|F_*|^{hG}$$
is an equivalence. Hence, homotopy fixed points do not
in general commute with realizations.
However, if each $(F_q)^{hG}$ is again a functor with
stabilization (in particular, connective) then
this map does induce  an equivalence on the associated spectra.

\bigskip
\noindent
\underbar{The Tate Map}

For $G$ a \underbar{finite} group,
the {\it Tate map} is a chain of natural maps of functors
with structure from $F_{hG}$ to
$F^{hG}$ which we now wish to define.
But first, we establish a sequence of natural equivalences
$$(G_+\wedge F)_{hG}\simeq \Omega^{\infty}F\simeq
(G_+\wedge F)^{hG}$$

\noindent
For $X$ a $G$--space, we let $\gamma$ be the $G$--equivariant
map
$$ G_+\wedge X\cong\bigvee_G X
\mapright inc.
\prod_G X \cong
Map\ (G_+,X)$$
that is:
$$\gamma(g\wedge x)(u) = \cases{x&if $g=u$\cr
                         *&otherwise\cr}
                         $$
Thus, if $X$ is $k$--connected, then $\gamma$ is $(2k-1)$--connected
by Blakers-Massey and we obtain the diagram:
$$\matrix{(G_+\wedge
X)_G&\mapleft\cong.&X&\mapright\cong.&Map(G_+,X)^G\cr
    \mapup\simeq.&&& &\mapdown\simeq.\cr
    (G_+\wedge X)_{hG}&&& &Map(G_+,X)^{hG}\cr
    \mapdown \gamma_{hG}.&&& &\mapup\gamma^{hG}.\cr
    Map(G_+,X)_{hG}&&& &(G_+\wedge X)^{hG}\cr
    }\leqno{(1)}$$
 Since $\gamma$ is $(2k-1)$--connected, so is
$\gamma_{hG}$;
but we do not
know anything about the connectivity of $\gamma^{hG}$. However,
if $F$ is functor with stabilization with $G$--action then
for all $X$
we know that the composite
map
$$
\matrix{
[{\hocolim}\ \Omega^n(G_+\wedge F(\Sigma^n X))]^{hG}\cr
\mapdown\gamma^{hG}.\cr
[{\hocolim}\ \Omega^nMap(G_+,F(\Sigma^n X)]^{hG}\cr
\mapdown\cong.\cr
[{\hocolim}\ Map(G_+,\Omega^n F(\Sigma^n X)]^{hG}
}\leqno{(2)}$$
is an equivalence. We also note that since $G$ is finite,
the natural $G$--equivariant map
$${\hocolim}\ Map(G_+,\Omega^n F(\Sigma^n X))
\mapright\simeq.
Map(G_+,{\hocolim}\  \Omega^n F(\Sigma^n X))
\leqno{(3)}$$
is an equivalence.
Thus, we can assemble all these remarks to obtain the following
sequence of natural equivalences
of spaces
(a symbol in ``$(\ )$'' indicates the previous statement which
implies the map is an equivalence)
$$\matrix{
\hfill (G_+\wedge F)_{hG}(X) = &
{\hocolim}\ \Omega^n [(G_+\wedge F(\Sigma^n X)_{hG}]\cr
   &  \mapdown\simeq\ (1).\cr
 \hfill\Omega^{\infty}{F}(X) =&
 {\hocolim}\ \Omega^n F(\Sigma^n X)\cr
     &\mapdown\simeq\ (1).\cr
     &[{\rm Map\ }(G_+,{\hocolim} \Omega^n F(\Sigma^n X)]^{hG}\cr
     &\mapup\simeq\ (3).\cr
     &[{\hocolim}\ \Omega^nMap(G_+,F(\Sigma^n X))]^{hG}\cr
     &\mapup\simeq\ (2).\cr
\hfill (G_+\wedge F)^{hG}(X) =&
  [{\hocolim}\ \Omega^n (G_+\wedge F(\Sigma^n X))]^{hG}\cr
     }\leqno{(A)}$$
The map on the last line is an equivalence since the previous lines show we
already
have an omega-spectrum.
The maps in (A) assemble (as $X$ varies)
into a natural sequence of equivalences of functors
with stabilization with $G$--action.

Using the $G$--equivariant equivalence $EG_+\wedge
X\mapright\simeq.X$,
we obtain
natural equivalences on $G$--homotopy orbits and homotopy fixed points.
We
can obtain $EG_+$ as the realization of a simplicial $G$--set
$\ [q]\mapsto \bigwedge^{q+1}G_+\ $
(the simplicial path space of the bar construction for $G$, with $G$ acting on
the $0$'th coordinate).
Thus,
$$EG_+\wedge F\  \cong_G\
|[q]\mapsto \bigwedge^{q+1}G_+\wedge F|$$

\bigskip\noindent
{\caps Definition 5.3:} (T. Goodwillie) The {\it Tate ``map''}
is the following natural diagram:
$$\matrix{
F_{hG}\cr
\mapup\simeq.\cr
(EG_+\wedge F)_{hG}\cr
\mapdown\cong.\cr
|[q]\mapsto (
\bigwedge^{q+1}G_+\wedge F
)|_{hG}\cr
\mapup\simeq.\cr
|[q]\mapsto (\bigwedge^{q+1}G_+\wedge F)_{hG}|\cr
\ \ \ \ \simeq\ (A)\cr
|[q]\mapsto (\bigwedge^{q+1}G_+\wedge F)^{hG}|\cr
\mapdown.\cr
|[q]\mapsto (\bigwedge^{q+1}G_+\wedge F)|^{hG}\cr
\mapdown\cong.\cr
|EG_+\wedge F|^{hG}\cr
\mapdown\simeq.\cr
F^{hG}\cr
}$$
(In the middle step,
$ \bigwedge^{q+1}G_+\wedge F$
should be viewed as
$G_+\wedge \bigwedge^{q}G_+\wedge F$.)

\bigskip
There is one case where the Tate map is easily seen to be
an equivalence: when
$E$ is a functor with stabilization with
$G$-action and $F = G_+\wedge E$. This follows from
the remarks following 5.2 or
the commuting
diagram (using the projection maps $\pi$)
$$\matrix{
|[q]\mapsto [\bigwedge^{q+1}G_+\wedge (G_+\wedge
E)]_{hG}|&\mapright\pi.&
(G_+\wedge E)_{hG}\cr
\ \ \ \ \simeq\ (A)& &\cr
|[q]\mapsto [\bigwedge^{q+1}G_+\wedge (G_+\wedge E)]^{hG}|
& &\ \ \ \ \simeq\ (A)\cr
\mapdown.& &\cr
|[q]\mapsto [\bigwedge^{q+1}G_+\wedge(G_+\wedge E)]|^{hG}
&\mapright\pi.&
(G_+\wedge E)^{hG}\cr
     }\leqno{(B)}$$
Another case which follows from this one
is that of functors with stabilization of the
form ${\rm Map}\ (G_+, E)$, where we have the stable equivalences of homotopy
orbits from (1) and of homotopy fixedpoints from (2).

\bigskip\noindent
{\caps Proposition 5.4}: (T. Goodwillie)
The Tate map for $F$ is an equivalence
if $F$ is either $U\wedge E$ or ${\rm Map\ }(U, E)$, where
$ E$ is a functor with stabilization with $G$--action and
$U$ is a pointed finite
free $G$--space (i.e., a simplicial $G$--set with finitely many
nondegenerate non-basepoint simplicies permuted freely by $G$).

\bigskip
The proof is by induction over
skeleta; the
cells attached at stage $n$ are dealt with by applying
the above discussion
 to the case $(E\wedge \bigvee^{t} S^n)\wedge G_+$ and ${\rm Map}\
(\bigvee^t S^n\wedge G_+, E)$, using the fact that after we apply
$\Omega^\infty$,
cofibrations become fibrations.

\bigskip
\noindent
{\caps Theorem 5.5}: (T. Goodwillie)
Let $U$ be a free finite based $G$--complex
of dimension $n$ and $W$ a $(n-1)$--connected based $G$--complex. Then
the spectrum associated to the functor with stabilization
$\Map(U,W)^G$ is naturally equivalent to
that associated to $\Map(U,W)_{hG}$.

\bigskip\noindent
{\it Proof}: Consider the diagram
$$\eqalign{\Map(U,W)^G&\mapright\alpha.\Map(U,\hocolim\Omega^k(S^k\wedge
W))^G\cr
   &\mapright\beta.\Map(U,\hocolim\Omega^k(S^k\wedge W))^{hG}\cr
   &\mapleft\gamma. (\hocolim\Omega^k(\Map(U,(S^k\wedge W))^{hG})\cr
   &\mapleft\delta.\hocolim\Omega^k(\Map(U,(S^k\wedge W))_{hG})\cr
   }$$
The first map, $\alpha$, is induced by the inclusion
$W\rightarrow \hocolim\Omega^k(S^k\wedge W)$. Since $W$ is
$(n-1)$--connected
this map is $(2n-1)$--connected. Since
$U$ is a free $G$--space of dimension $n$, the map
itself is $(n-1)$--connected.

The second map, $\beta$, is the canonical map from fixed points
to homotopy fixed points. It is an equivalence because this is always
so for function spaces $Map(U,?)$ where $U$ and $?$ are $G$-spaces
and $U$ is free.

The third map, $\gamma$, is an equivalence because $U$ is a finite
complex; this has nothing to do with the $G$--action.

The fourth map, $\delta$, is the Tate map, and is an
equivalence by 5.4.

\bigskip
\noindent{\it Notation}:
Let ${\sc M}\subset {\bf N}^\times$ be a full subcategory, and let
$M\in{\sc M}$.
Let $p_1,,\ldots,p_t$ be the distinct prime divisors of $M$.
For $U\subseteq\{1,\ldots, t\}$, we let
$<U> =
\prod_{u\in U}p_u$ ($<\emptyset> = 1$).  Assume  ${M\over <U>}\in {\sc M}$ for all $U\subseteq\{1,\ldots,t\}$, and
let $\vec M$ denote the full subcategory of $\sc M$
with objects $\{{M\over <U>}|U\subseteq\{1,\ldots,t\}\}$.
We define
$\vec M - M$ to be the full subcategory of $\vec M$
generated by all the objects {\it except} $M$.

\noindent
{\it Assumption}: $\sc M$ is a subcategory of ${\bf N}^{\times}$
and $M\in\sc M$ is such that
$\sc M$ is covered by the object of, and compositions of the morphisms of,
 $\vec M$ and $\sc M-M$
(that is,
there does not exist an $M'\in\sc M-M$ such that $M|M'$).

\bigskip\noindent
{\it Remark}:
Since $\vec M-M$ is the intersection of
$\vec M$ and $\sc M-M$, for any functor $F$ from $\sc M$
to functors with stabilization
the following natural
diagram is
homotopy cartesian:
$$\matrix{\holim_{\sc M}F&\mapright .&\holim_{\sc M - M}F\cr
  \mapdown .& &\mapdown .\cr
  \holim_{\vec M}F&\mapright .&\holim_{\vec M-M}F\cr}$$
Since $M$ is initial in $\vec M$, the natural map
$F(M)\mapright .\holim_{\vec M}F$ is an equivalence and
the homotopy fiber of
the composite $F(M)\mapright .\holim_{\vec M-M}F$ is
naturally equivalent to the total fiber of the $t$--dimensional
cube
determined by $F$ on $\vec M$ (see  [G2], 1.1b).

If we consider the functor $ U(F;P)$
defined after Definition 4.3, we see that the homotopy fiber of
the restriction map from $W_{\sc M}(F;P)$
to $W_{\sc M-M}(F;P)$ is naturally equivalent to the total
fiber of the $t$--dimensional cube determined by $U(F;P)$
on $\vec M$. Since $U(F;P)$ takes values in simplicial
functors with stabilization, this total fiber is naturally
equivalent to the realization of the total fibers computed
in each simplicial dimension separately.

\bigskip
\noindent
{\caps Proposition 5.6}:
If $M\in\sc M$ is such that $\sc M$ is covered by
the objects in, and compositions of the morphisms in,
$\vec M$ and $\sc M-M$ then there is a natural
(in $F$ and $P$) chain of equivalences of functors
with stabilization
$$U^M(F;P)_{hC_M} \simeq
{\rm hofib}\ [W_{\sc M}(F;P)\mapright.W_{\sc M-M}(F;P)]
.$$

\bigskip\noindent
{\it Proof}:
For $\underline{X}\in I^{\vec k+1}$, set
$$Z = (S^{\underline{X}})^{\wedge M}$$
$$Y =  V(\underline{X},\sc O)^{\wedge M}$$
(recall Definition 2.2),
and $x = \Sigma_{i=0}^k|X_i|$.

By the above remark,
$$
{\rm hofib}\ [W_{\sc M}(F;P)\mapright.W_{\sc M-M}(F;P)]
\simeq
{\rm hofib}\ [W_{\vec M}(F;P)\mapright.W_{\vec M-M}(F;P)]
.$$
This is the total fiber of
the $t$--dimensional cube
$$\sc X(U) = U^{M/<U>}(F;P)^{C_{M/<U>}}
=\Map(Z^{C_{<U>}},Y)^{C_M}$$
with  maps by the restriction to fixed subsets.
By applying the functor $\Map(\ ,Y)^{C_M}$, we see that
this is the
same as $\Map({\sc U},Y)^{C_M}$ where ${\sc U}$ is
the total cofiber of the $t$--dimensional cube
with
$\sc Y(U) = Z^{C_{<U>}}$
and structure maps given by
inclusion.
Thus, ${\sc U} = Z/Z'$ where $Z'$ is the push--out
of the diagram
$$\sc Z(U) = (S^{\underline{X}})^{\wedge (M / <U>)}
\ \ \ \ \ \ \ \ U\subseteq
\{1,\ldots,n\}
; U\not=\emptyset$$
with the maps given
by inclusions. 

Now ${\sc U}$ is a finite free based
$C_M$--space of dimension  $Mx$.
Since $Y$ is $(Mx-1)$--connected,
$\Map({\sc U},Y)^{C_M}$ is naturally equivalent to $\Map({\sc U},Y)_{hC_M}$
by Theorem 5.5.
The quotient map $Z\rightarrow {\sc U}$ produces a natural map
$$\Map({\sc U},Y)_{hC_M}\mapright\epsilon.\Map(Z,Y)_{hC_M}.$$
Since ${\rm dim}(Z')= {\rm max} \{{Mx\over p_i}|1\leq i\leq n\}= Mx/p$ (for $p$ the smallest prime divisor 
 of $M$)
and $S^\ell\wedge Y$ is $(\ell+Mx-1)$--connected, the
map
$$\Map ({\sc U},(S^\ell\wedge Y))\rightarrow \Map(Z,(S^\ell\wedge Y))$$
is
$(\ell+Mx(1-1/p)-1)$--connected and
hence $\epsilon$ is
$(Mx(1-1/p)-1)$--connected.
Since the quotient map from $Z$ to ${\sc U}$ is functorial in
$I^{\vec k+1}$, we can take the homotopy colimit
with respect to $I^{\vec k+1}$ and hence obtain an equivalence
$$\hocolim_{I^{\vec k+1}}
\Map ({\sc U},(S^\ell\wedge Y))\rightarrow
\hocolim_{I^{\vec k+1}}\Map(Z,(S^\ell\wedge Y)).$$

\noindent
Thus, we have obtained a natural sequence of equivalences
$$\bigl(U^M(F;P)_{[k]}\bigr)_{hM}\simeq
{\rm hofib}[W_{\sc M}(F;P)_{[k]}\rightarrow W_{\sc M-M}(F;P)_{[k]}].$$
It remains to check that these maps respect the simplicial
operators. This straightforward but messy detail is left
for the interested reader.

\bigskip\noindent
{\caps Corollary 5.7:}
For any FSP $F$ and bimodule $P$
there is a natural chain of equivalences of
functors
$$\eqalign{
U^n(F;P)_{hC_n} &\simeq {\rm hofib}[W_n(F;P)\mapright R.W_{n-1}(F;P)]\cr
U^{p^n}(F;P)_{hC_{p^n}} &\simeq
{\rm hofib}[W^{(p)}_n(F;P)\mapright R.W^{(p)}_{n-1}(F;P)]\cr
}$$

\bigskip\noindent
{\caps Corollary 5.8 :}
Let $P$ be a $(k-1)$--connected $F$-bimodule.
The natural map
$$W(F;P)\mapright R. W_n(F;P)$$ is $((n+1)k-2)$--connected
and the natural map
$$W^{(p)}(F;P)\mapright R^{(p)}.W^{(p)}_n(F;P)$$
is $(p^{n+1}k-2)$--connected.

\bigskip\noindent
{\it Proof}: By lemma 3.4, $U^n(F;P)$ is $kn+(n-1)$
connected for all $n\geq 1$. Since homotopy orbits
preserve connectivity the result follows
from corollary 5.7.

\bigskip

\noindent{\bf 6. \underbar{More Properties of $U^n$ and $W$}}
\bigskip
In this section we will develop several useful properties of $U^n$ and $W$  which will be
needed to define the map from algebraic K-theory to $W$, and which can also be useful
in calculations.  These results are all analogous to results from [DMc2], where they were
proved for the $\THH$ (i.e. $U^n$ for $n=1$) case, with or without coefficients in a
bimodule.   We will restict our attention to FSP's over small categories  (note that in the [DMc2] nomenclature, what we here call an FSP is a unital ring functor;
they reserve the name FSP for the case where the category in question consists of a single point, as in
B\"okstedt's original definition).  Later in the section we will restrict ourselves further to small 
linear categories (where the homomorphism set between any two objects is an abelian group and composition is bilinear), and to a particular FSP on them (see the first example after Definition 1.3
above).

The proofs from [DMc2] can be adapted as explained below
to give homotopy equivalences between $U^n( F;P_1,\ldots,P_n)$ for  different FSP's $F$ over 
categories $\sc C$ and $F$-bimodules $P_1,\ldots, P_n$ (the proofs that the maps are indeed homotopy equivalences work when the bimodules are distinct; we mention only the case $P_1=\cdots=P_n$, which is what we will mostly use, in order to simplify notation).  In all the cases, there 
are maps inducing these homotopy equivalences which are naturally $C_n$-equivariant
when $P_1=\cdots=P_n$.    This turns out to be enough to show that they are $C_n$-equivalences:  Using the fundamental sequence of Proposition 5.6 (for $\sc M$ equal
to all of $n$'s divisors), by induction on $n$ (since a $C_n$-equivariant map which is a homotopy
equivalence induces a homotopy equivalence on the $C_n$ homotopy quotient) we can see that
they induce an equivalence on the $C_n$-fixedpoints of $U^n$.  Using groupings
$U^n(F;P)\simeq U^m(F;P^{\oF {n\over m}})$ as in Lemma 2.6 we can get an equivalence of
the $C_m$-fixedpoints for any $m\vert n$ in a similar way.
Once we know that the maps in each of these claims are $C_n$ equivalences, it follows 
(except in Lemma 6.4 where the map is from a  direct limits of $U^n$'s) that they induce equivalences on $W$ and all its variants as well.

\medskip
If F is an FSP over a small category 
$\sc C$ and $P$ an $F$-bimodule, then for any small category $\sc D$
and functor $\phi:\sc D\to \sc C$ we can get an FSP $\phi^*F$ over $\sc D$ by letting
$\phi^*F_{d,d'}(X)=F_{\phi(d), \phi(d')}(X)$ for all $d,d'\in\sc D$, $X\in \sc S_*$.  We can similarly
define the $\phi^*F$-bimodule $\phi^*P$.

\bigskip
\noindent
{\caps Lemma 6.1:}  Let $\phi_1, \phi_2:\sc D\to\sc C$ be two naturally isomorphic functors between small categories, and let
$F$ be an FSP over $\sc C$ and $P$ an $F$-bimodule.  Then the natural isomorphism induces a $C_n$-homeomorphism
$$U^n(\phi_1^*(F);\phi_1^*(P))  \mapright\cong. U^n(\phi_2^*(F);\phi_2^*(P)) $$
for all $n$ and therefore a homeomorphism on $W$.

\bigskip
\noindent
{\it Proof}:   Analogous to that of Lemma 1.6.2 of [DMc2]: the natural isomorphism $\eta$
induces an
equivalence  $F(\eta(a)^{-1},\eta(b))(id_X): F(\phi_1(a), \phi_1(b))(X) \mapright\sim. F(\phi_2(a), \phi_2(b))(X) $
for all $X$, and similarly for $P$.   These are compatible with the multiplicative structure. 
\bigskip
\noindent
{\caps Proposition 6.2:}  Let $\phi:\sc D\to \sc C$ be an equivalence of categories, and let $F$
be an FSP over $\sc C$ and $P$ an $F$-bimodule.  Then $\phi $ induces a $C_n$-equivalence
$$U^n(\phi^*(F);\phi^*(P))   \mapright\sim. U^n(F;P)$$
and therefore an equivalence on $W$.

\bigskip
\noindent
{\it Proof}:  Using the natural transformation
 $\psi:\sc C\to \sc D$ such that both compositions are naturally  isomorphic to the identity and Lemma 6.1, as in the proof of Lemma 1.6.6 in [DMc2].

\bigskip
\noindent
{\it Example} 6.3:  Let $\sc C$ and $\sc D$ be small linear categories, with an
equivalence of categories $\phi:\sc D\to \sc C$.  Then we can look at the FSPs $\u{\sc C}$
and $\u{\sc D}$, as defined after Definition 1.3.  There is a natural isomorphism
between the FSPs $\phi^* \u{\sc C}$ and $\u{\sc D}$ on $\sc D$: Whenever the identity is 
naturally isomorphic to a functor $F:\sc E\to\sc E$ then $F^*:\Hom_\sc E(e_1,e_2)
\to  \Hom_\sc E(F(e_1),F(e_2))$ is one-to-one and onto for all $e_1, e_2\in \sc E$ (because
if $H:id_\sc E\to F$ is a natural isomorphism, for every $\alpha:e_1\to e_2$, $F(\alpha)=H(e_2)\alpha
H(e_1)^{-1}$).  So we can use this on   $\Hom_\sc D(d_1, d_2)\mapright\phi^*.
\Hom_\sc C(\phi(d_1),\phi( d_2))$ and 
$\Hom_\sc C(c_1,c_2)\mapright\psi^*.
\Hom_\sc D(\psi( c_1),\psi(c_2))$ first to establish that $\phi^*$ is one-to-one on morphism sets and
$\psi^*$ onto, and then in the opposite order to establish that $\psi^*$ is one-to-one on morphism sets and
$\phi^*$ onto.

Now say we have two functors
$F,G:\sc C\to\sc B$ and a $\u{\sc C}$-bimodule of the form
$P(a,b)(X)= \sc B(F(a),G(b))\otimes_{\bf Z}\tilde{\bf Z}[X]$  (see the
second example after Definition 1.5 above).  Then $\phi^*P$ is a bimodule on the same form on $\sc D$,  corresponding to the functors $F\circ \phi$ and $G\circ \phi$, and we get a $C_n$-equivalence
$$U^n(\u{\sc D}; \phi^*P)  \mapright\sim. U^n(\u{\sc C};P).$$

\bigskip
 We now want to show that our construction of $U^n$ commutes with direct limits.  Call the category we are working on $\sc C$, and assume that there is a directed set of subcategories
 $\sc C_j$  of $\sc C$, $j\in J$, such that for any object $c\in \sc C$ there is $j\in J$ with
 $c\in \sc C_j$.  If $J$ satisfies this condition, we say that it is a {\it saturated} directed set in 
 $\sc C$.  Note that this condition is really a condition on the underlying sets of the small 
 categories involved. 
 
 \bigskip
\noindent
{\caps Lemma 6.4:}  If $J$ is a saturated directed set in $\sc C$, and $F$ is an FSP on $\sc C$
with $P$ an $F$-bimodule, then we have a $C_n$-equivalence
$$\lim_{j\in J}\U^n(F\vert_{C_j};P\vert_{C_j}) \mapright\sim. U^n(F;P).$$

\bigskip
\noindent
{\it Proof}: As in the proof of Lemma 1.6.9 in [DMc2], this follows from the fact that any map from
a sphere $S^{\sqcup\underline{X}}$ (which is compact) to $V(\underline{X},\sc C)$ (see Definition 2.3 above)
has a compact image, and therefore can intersect only finitely many summands which each involve 
only finitely many elements of $\sc C$, by commuting colimits and homotopy colimits.
The $C_n$ equivalence of $U^n$ is proved as usual, but note that this lemma does not imply a similar result for $W$ because of the problem of commuting direct and inverse limits.

 \bigskip
\noindent
{\caps Lemma 6.5:}  If we have FSPs $F$ and $F'$ over $\sc C$ and bimodules $P$ and $P'$
over $F$ and $F'$, respectively, and a map $(f,g):(F;P)\to(F';P')$ (see Definition 1.5) so that
$f$ and $g$ induce a stable equivalence on the associated spectra ${\bf F} \mapright\sim. {\bf F'}$,
${\bf P} \mapright\sim. {\bf P'}$ (see Definition 1.2), then we have a $C_n$-equivalence
$$U^n(F;P)   \mapright{(f;g)_*}. U^n(F';P')$$
and therefore an equivalence on $W$.

\bigskip
\noindent
{\it Proof}: Since the associated spectra are all equivalent, we get an equivalence of the $k$-simplices in $U^n$ for all $n$.

\bigskip
Given an FSP $F_1$ over a small category $\sc C_1$ with an $F_1$-bimodule $P_1$, and
an FSP $F_2$ over a small category $\sc C_2$ with an $F_2$-bimodule $P_2$, one can define
functors with stabilization over $\sc C_1\times \sc C_2$
$$\eqalign{
(F_1\times F_2)((a_1, a_2), (b_1, b_2))(X)  & =F_1(a_1, b_1)(X) \times F_2(a_2, b_2)(X)
\cr 
(F_1\vee F_2)((a_1, a_2), (b_1, b_2))(X) & =F_1(a_1, b_1)(X) \vee F_2(a_2, b_2)(X)
.}$$
Similar definitions can be made for bimodules.  Note that $F_1\times F_2$ is an FSP;  
$F_1\vee F_2$ has no unit, so is not an FSP.  However the inclusion of
the latter in the former induces a stable equivalence of the associated spectra, so if one used
the definition of $U^n$ on $F_1\vee F_2$ with coefficients in $P_1\vee P_2$, as in Lemma 6.5 above one would get
the same thing as $U^n( F_1\times F_2;P_1\times P_2)$.
One can also define an FSP over $\sc C_1\amalg \sc C_2$
$${ 
(F_1\amalg F_2)(a,b)(X) \cases{F_1(a,b)(X) &  if  $a,b\in \sc C_1$\cr
F_2(a,b)(X) &  if  $a,b\in \sc C_2$ \cr
* & if $a,b$ lie in different $\sc C_i$,
} }$$
and similarly define an $F_1\amalg F_2$-bimodule $P_1\amalg P_2$.  

 \bigskip
\noindent
{\caps Lemma 6.6:} 
For FSPs $F_1$ and $F_2$ over small categories $C_1$ and $C_2$, respectively, with bimodules
$P_1$ and $P_2$ we have a $C_n$-equivalence
$$U^n( F_1\amalg F_2; P_1\amalg P_2)  \mapright\sim.U^n(F_1;P_1)\times U^n(F_2;P_2)$$
inducing an equivalence on $W$.

\bigskip
\noindent
{\it Proof}:  Following the proof of Lemma 1.6.13 in [DMc2], 
if we pick $\u{X}$ and look at $V(\u{X}, \sc C_1\amalg \sc C_2)$ calculated with respect to the
FSP 
$F_1\amalg F_2$ and the bimodule $P_1\amalg P_2$, we can see that only summands
which correspond to sequences of elements which are all in $\sc C_1$ or all in $\sc C_2$ are non-trivial.  Thus  for every $\u X$ we have  
$$V(\u{X}, \sc C_1\amalg \sc C_2)= V_{(F_1;P_1)}(\u{X},\sc C_1)\vee V_{(F_2;P_2)}(\u{X},\sc C_2),$$
and so the limit of $\u{Map}(S^{\sqcup\u X}, V(\u{X}, \sc C_1\amalg \sc C_2))$
will be weakly equivalent to the limit of 
$\u{Map}(S^{\sqcup\u X}, V_{(F_1;P_1)}(\u{X},\sc C_1))\times V_{(F_2;P_2)}(\u{X},\sc C_2)).$

 \bigskip
\noindent
{\caps Lemma 6.7:} 
For FSPs $F_1$ and $F_2$ over small categories $C_1$ and $C_2$, respectively, with bimodules
$P_1$ and $P_2$.  Then the projections to both coordinates define a $C_n$-equivalence
$$U^n( F_1\times F_2; P_1\times P_2)  \mapright\sim.U^n(F_1;P_1)\times U^n(F_2;P_2)$$
inducing an equivalence on $W$.

\bigskip
\noindent
{\it Proof}:  Call the map induced by the projections 
(which is a $C_n$-equivariant map)  $f$; we need to show that it
is an equivalence.  Like in the proof of Lemma 1.6.15 in [DMc2], one can construct 
a commutative diagram 
$$\matrix{U^n( F_1\vee F_2; P_1\vee P_2)  &\mapright g .&U^n( F_1\amalg F_2; P_1\amalg P_2) \cr
  \mapdown{{\rm incl}_*} .& &\mapdown{{\rm Lemma} \ 6.6}  .\cr
 U^n( F_1\times F_2; P_1\times P_2) &\mapright f.&U^n(F_1;P_1)\times U^n(F_2;P_2)\cr}$$
with the vertical maps weak equivalences.  The map $g$ is obtained by restricting the maps
$$V_{(F_1\times F_2;P_1\times P_2)}(\u{X}, \sc C_1\times \sc C_2)\to  V_{(F_1;P_1)}(\u{X},\sc C_1)\times V_{(F_2;P_2)}(\u{X},\sc C_2),$$
used to define $f$ (induced by the product of the projections) to $V_{(F_1\vee F_2;P_1\vee P_2)}(\u{X}, \sc C_1\times \sc C_2)$, and observing that they give maps 
$$V_{(F_1\vee F_2;P_1\vee P_2)}(\u{X}, \sc C_1\times \sc C_2)\to  V_{(F_1;P_1)}(\u{X},\sc C_1)\vee V_{(F_2;P_2)}(\u{X},\sc C_2)=V_{( F_1\amalg F_2; P_1\amalg P_2)}(\u{X}, \sc C_1\amalg \sc C_2).$$
In the opposite direction, one can define a  map 
$U^n( F_1\amalg F_2; P_1\amalg P_2)\mapright i.  U^n( F_1\vee F_2; P_1\vee P_2)$ 
by first mapping $\sc C_1\amalg\sc C_2\to \sc C_1\times\sc C_2$
using some fixed $a\in \sc C_2$ as a `filler'  second coordinate for $\sc C_1$ and some fixed
$b\in \sc C_1$ as a `filler' first coordinate for $\sc C_2$, and then mapping $F_i$ or $P_i$ into
$F_1\vee F_2$ or $P_1\vee P_2$, respectively.  Now $g\circ i=id$,  and the proof in 1.6.15 [DMc2] works to show that $i\circ g\cong id$ (if one ignores the cyclic action, it does not matter whether the coordinates are the FSP or the bimodule).  Since the vertical maps are known to be equivalences, 
the fact that $g$ is an equivalence implies that $f$ is one, too.

 \bigskip
\noindent
{\caps Definition  6.8:} 
Given a functor with stabilization $A$ over $\sc C$, we can define its $\ell\times \ell$ matrices as a functor with
stabilization over $\sc C^\ell$ in two ways:
$$M_\ell(A)((c_1,\ldots,c_\ell),(c_1',\ldots,c_\ell'))(X)=\prod_{r=1}^\ell\bigvee_{s=1}^\ell A(c_r, c_s')(X),$$
$$M_\ell(A)_{\vee}((c_1,\ldots,c_\ell),(c_1',\ldots,c_\ell'))(X)=\bigvee_{r=1}^\ell\bigvee_{s=1}^\ell A(c_r, c_s')(X).$$
Now if $F$ is an $FSP$, $M_\ell(F)$ is an FSP too, using matrix multiplication (see 1.2.6 in  [DMc2];
since the FSP multiplication sends $F(b,c)(X)\wedge F(a,b)(Y)\to F(a,c)(X\wedge Y)$, one should
think of $A(c_r, c_s')(X)$ as the $(s,r)$'th entry in the matrix);
$M_\ell(F)_{\vee}$ using the same multiplication does not have a unit but of course the associated
spectra are stably equivalent.  If $F$  is an FSP and $P$ is an $F$-bimodule, $M_\ell(P)$
is a $M_\ell(F)$-bimodule.

One can also define upper-triangular matrices (see above comment about indexing)
$$T_\ell(A)((c_1,\ldots,c_\ell),(c_1',\ldots,c_\ell'))(X)=\prod_{r=1}^\ell\bigvee_{s=1}^r A(c_r, c_s')(X),$$
$$T_\ell(A)_{\vee}((c_1,\ldots,c_\ell),(c_1',\ldots,c_\ell'))(X)=\bigvee_{r=1}^\ell\bigvee_{s=1}^r A(c_r, c_s')(X),$$
 and again if $F$  is an FSP then  $T_\ell(F)$ is one too,  and $T_\ell(A), \ T_\ell(A)_\vee$ are stably equivalent for any $A$.

  \bigskip
\noindent
{\caps Proposition  6.9} (Morita Equivalence): Let $F$ be an FSP on $\sc C$ and let $P$ be
an $F$-bimodule.  Then there is a $C_n$-equivalence
$$U^n( F;P)  \mapright\sim.U^n(M_\ell(F);M_\ell(P))$$
inducing an equivalence on $W$. 

 \bigskip
\noindent
{\it Proof}: We will define the map which induces this equivalence; it will as usual be $C_n$-equivariant.  The proof that it is a homotopy equivalence is completely analogous to the proof of Proposition 1.6.18 in [DMc2].  The map is defined by picking some element $c_0\in\sc C$,
 and using it to embed $\sc C \mapright{{ i}}. \sc C^\ell$ by $a\mapsto (a,c_0,\ldots, c_0)$ on objects and $f\mapsto(f,id_{c_0},\ldots, id_{c_0})$ on morphisms.  Then on $\sc C$
 we map $F$ to ${{i }}^*M_n (F)_\vee$ by including, for every $a, b\in\sc C$ and any finite simplicial $X$, $F(a,b)(X)$ as the $(1,1)$'st coordinate in $M_n (F)_\vee( (a,c_0,\ldots, c_0),
  (b,c_0,\ldots, c_0))$, which in turn includes into $M_n (F)$.
 
   \bigskip
\noindent
{\caps Proposition  6.10:} Let $F$ be an FSP on $\sc C$ and let $P$ be
an $F$-bimodule.  Then the inclusion of the diagonal matrices in the upper-triangular ones
induces a $C_n$-equivalence
$$U^n( \prod_{i=1}^\ell F; \prod_{i=1}^\ell P)  \mapright\sim.U^n(T_\ell(F);T_\ell(P))$$
inducing an equivalence on $W$. 
 \bigskip
\noindent
{\it Proof}: There is an obvious map from the upper-triangular matrices to the diagonal ones---collapsing all the off-diagonal terms---which shows that the above map must be the inclusion of
 a retract.  The proof that it is in fact an equivalence is analogous to that of Proposition 1.6.20 in [DMc2], and is done
 by replacing $U^n( \prod_{i=1}^\ell F; \prod_{i=1}^\ell P) $ with the equivalent $U^n( \bigvee_{i=1}^\ell F; \bigvee_{i=1}^\ell P)$  and $U^n(T_\ell(F);T_\ell(P))$ with the equivalent
$ U^n(T_\ell(F)_\vee;T_\ell(P)_\vee)$.
 
  \bigskip
\noindent
{\caps Definition  6.11:}
 We now restrict ourselves to small linear categories $\sc C$,  to  FSPs of the
 form $\u{\sc C}(a,b)=\sc C(a,b)\otimes_{\bf Z}\tilde{\bf Z}[X]$, as described in the example after Definition 1.3, and to bimodules over
 them of the form 
 $$P(a,b)(X)= \sc B(G_1(a),G_2(b))\otimes_{\bf Z}\tilde{\bf Z}[X]$$
 for some other linear category $\sc B$ and two functors $G_1, G_2:\sc C\to\sc B$ which respect the linear structure of the morphism sets.
 In this case we can define FSPs
 $$m_\ell(\u\sc C)((c_1,\ldots,c_\ell),(c_1',\ldots,c_\ell'))(X)=\bigl(\bigoplus_{r=1}^\ell\bigoplus_{s=1}^\ell \sc C(c_r, c_s')\bigr)\otimes_{\bf Z}\tilde{\bf Z}[X],$$
 $$M_\ell(\u\sc C)_\oplus((c_1,\ldots,c_\ell),(c_1',\ldots,c_\ell'))(X)=\bigoplus_{r=1}^\ell\bigoplus_{s=1}^\ell \sc C(c_r, c_s')\otimes_{\bf Z}\tilde{\bf Z}[X].$$
 Observe that the two are, in fact, homeomorphic on any $X$.  One can do the same construction for matrices over $P$ of the above form.  Observe also that the obvious
 inclusions $M_\ell(\u\sc C)\to M_\ell(\u\sc C)_\oplus$, $M_\ell( P)\to M_\ell(P)_\oplus$ are stable equivalences, and therefore by Lemma 6.5 above we have $C_n$-equivalences
 $$U^n(M_\ell(\u\sc C);M_\ell( P)) \mapright\sim.U^n(M_\ell(\u\sc C)_\oplus;M_\ell( P)_\oplus)=
 U^n(m_\ell(\u\sc C);m_\ell( P)).$$
 One can similarly construct
 upper-triangular versions
 $t_\ell(\u\sc C)$, $t_\ell(P)$, $T_\ell(\u\sc C)_\oplus$, $T_\ell(P)_\oplus$
and get 
 $$U^n(T_\ell(\u\sc C);T_\ell( P)) \mapright\sim.U^n(T_\ell(\u\sc C)_\oplus;T_\ell( P)_\oplus)=
 U^n(t_\ell(\u\sc C);t_\ell( P)).$$

  \bigskip
\noindent
{\caps Proposition  6.12} (Cofinality): Let $\sc D$ be a  a small additive category (that is, a small linear category
with a notion of $\oplus$ on the objects which corresponds to taking the direct sum of abelian groups on the morphisms, and with a zero object).  Let $\sc C$ be a full subcategory of it which is cofinal, that is: for
any $d\in\sc D$ there is $d'\in\sc D$ such that $d\oplus d'\in\sc C$.  Let $P$ be an FSP of the
form  $P(a,b)(X)= \sc B(G_1(a),G_2(b))\otimes_{\bf Z}\tilde{\bf Z}[X]$ on $\sc D$ for some
additive category functors $G_1, G_2:\sc D\to \sc B$.  Then the inclusion induces a $C_n$-equivalence
$$U^n(\u\sc C;P\vert_{\sc C}) \mapright\sim.U^n(\u\sc D;P)$$
and an equivalence on $W$. 
 \bigskip
\noindent
{\it Proof}: Since by construction $P$ respects the direct sum structure, the proof of Lemma 2.1.1 of
 [DMc2] works if we use it in some of the coordinates. 
\bigskip

The following proposition connects $U^n(\u{\sc P_R};P)$  with the definition of $U$ of the FSP associated to a ring with 
coefficients in  a bimodule which is analogous to B\"okstedt's original definition of ${\rm THH}$ in [B].  The former will be needed to construct the map from K-theory in section 9 below; the latter is more compact  and easier to use for calculations.
\bigskip
\noindent
{\caps Proposition  6.13} (Another Morita Equivalence): Let $R$ be an associative ring with unit.
We can view $R$ as the full subcategory on the rank $1$ free module inside $\sc P_R$, the category of finitely generated projective right $R$-modules.  Let $P(a,b)(X)= \sc B(G_1(a),G_2(b))\otimes_{\bf Z}\tilde{\bf Z}[X]$ for some  additive category functors $G_1, G_2:\sc P_R\to \sc B$.  Then the inclusion $R\hookrightarrow
\sc P_R$ induces a $C_n$-equivalence
$$U^n(\u R;P\vert_R) \mapright\sim.U^n(\u{\sc P_R};P)$$
and an equivalence on $W$. 
 \bigskip
\noindent
{\it Proof}: Following the proof of Proposition 2.1.5 in [DMc2], we let $\sc F_R$ be the category of
 finitely generated free right $R$-modules, and $\sc F_R^k$ its full subcategory on the modules of
 rank less than or equal to $k$.  Then the inclusion $m_k(R)\hookrightarrow \sc F_R^k$, where we regard $m_k(R)$ as the full subcategory on a rank $k$ free module, is an equivalence of categories.  Note that on the category with one object, the FSP $m_\ell(\u R)$ of Definition 6.11
is the same FSP as what we would call $\u{m_\ell(R)}$, the one associated to the full subcategory described above.

So we get a commutative diagram of $C_n$-equivariant maps
$$\matrix{U^n(\u R;P\vert_R)  &\!\!\!\!\! \mapright.&\!\!\!\!\! U^n(\u{\sc P_R};P) &\!\!\!\!\!  \mapleft{{\rm Prop \ 6.12}} .
& U^n(\u{\sc F_R};P\vert_{\sc F_R})\cr
  \mapdown{{\rm Prop\ 6.9}} .& &&&\mapup{{\rm Lemma \ 6.4}}  .\cr
 \lim\limits_{k\to\infty}U^n(M_k(\u R); M_k(P\vert_R)) &\!\!\!\!\! \mapright{{\rm Def \ 6.11}}.&\!\!\!\!\! 
 \lim\limits_{k\to\infty}U^n(m_k(\u R); m_k(P\vert_R)) &\!\!\!\!\! \mapright{{\rm Ex \ 6.3}}.&\!\!\!\!\!
 \lim\limits_{k\to\infty}U^n(\u{\sc F_R^k};P\vert_{\sc F_R^k}) 
\cr}$$
where the labels on the arrows indicate from where it follows that those maps are homotopy
equivalences.  Thus the unlabeled map must be a homotopy equivalence as well.

\bigskip
In  [DMc1] the authors construct a map from $K(R;M)$ to  topological Hochschild homology,
and then show that it is the map from $K(R;M)$ to the first layer of its Goodwillie Taylor tower.
They start with another functor which maps very naturally to  topological Hochschild homology
and then look at the functor induced on the Waldhausen S-constructions of domain and range.
In Section 9 below, we will follow the same method.  So we briefly recall the S-construction from [W] and
[DMc1].  

Given an exact category (an additive category with a compatible notion of exact sequences)
$\sc C$, one can define for any $n\geq 0$ another exact category $S_n\sc C$
whose objects are sequences of admissible monomorphisms  $0=c_0\hookrightarrow c_1
\hookrightarrow\cdots\hookrightarrow c_n$ with particular identifications of $c_j/c_i$
for all $i\leq j$, and whose morphisms are commuting diagrams.
Assembled over all $n$, with the obvious composition maps for $\partial_i$, $0<i<n$, omission of
and quotienting by $c_1$ for $\partial_0$, and omission of $c_n$ for $\partial_n$, these form a simplicial exact category. 

If the original category $\sc C$ is split exact, that is: all exact sequences split in it, then $S_.\sc C$ 
is a split simplicial exact category.  One can take  functors from small categories to  spaces
or spectra and define them on 
a simplicial exact category levelwise, and then realize.  One can also define the iterated
S-construction $S_.^{(k)}\sc C$ 
to be the simplicial exact category obtained by taking the diagonal of the $k$-simplicial exact
category one would get by iterating the process $k$ times.

Note that if $\sc C$ is a small exact category, and we have a bimodule over the FSP $\u\sc C$
of the form  $P(a,b)(X)= \sc B(G_1(a),G_2(b))\otimes_{\bf Z}\tilde{\bf Z}[X]$ for some  exact
 functors $G_1, G_2$ from $\sc C$ to an exact category $\sc B$, the $G_i$ induce 
 simplicial exact functors $S_.G_i:S_.\sc C\to S_.\sc B$ and so we can define a $\u{S_.\sc C}$
 bimodule 
 $$S_nP(\bar a,\bar b)(X)= S_n \sc B(S_nG_1(\bar a),S_nG_2(\bar b))\otimes_{\bf Z}\tilde{\bf Z}[X]$$
for all $\bar a, \bar b\in S_n\sc C$ for all $n\geq 0$.  

\bigskip
\noindent
{\caps Proposition  6.14:} Let $\sc C$ be a split small exact category and let $P$ be a bimodule
over the FSP  $\u\sc C$
of the form  $P(a,b)(X)= \sc B(G_1(a),G_2(b))\otimes_{\bf Z}\tilde{\bf Z}[X]$ for some  exact
 functors $G_1, G_2:\sc C\to\sc B$.  Then there is a 
 $C_n$-equivalence 
 $$U^n(\u\sc C;P) \mapright\sim.\Omega \vert U^n(\u{S_.\sc C};S_.P)\vert $$
 induced by the adjoint to the map $\Sigma U^n(\u\sc C;P) \mapright. \vert U^n(\u{S_.\sc C};S_.P)\vert $ coming from the identification of $S_0\sc C$ with the trivial category and
 $S_1\sc C$ with $\sc C$ via $\{0\hookrightarrow c\}\leftrightarrow c$.
 This equivalence yields an equivalence
 $$W(\u\sc C;P) \mapright\sim.\Omega \vert W(\u{S_.\sc C};S_.P)\vert. $$
  \bigskip
\noindent
{\it Proof}:  As in the proof of Proposition 2.1.3 in [DMc2], the key point is that because $\sc C$
 is split exact, for any $k\geq 0$, if we look at the functor $\sc C^k\mapright f. S_k\sc C$ defined by
 $(c_1,\ldots, c_k)\mapsto\{0\hookrightarrow c_1\hookrightarrow c_1\oplus c_2 \hookrightarrow
 \cdots\hookrightarrow c_1\oplus\cdots\oplus c_k\}$, it is an equivalence of categories.
 The morphisms of $S_k\sc C$
 pull back exactly to the upper-triangular matrices so $f^*(\u{   S_k\sc C})=t_k(\u \sc C)$ of
 Definition 6.11 above.  Similarly, since the $G_i$ are exact functors, they send direct sums to direct sums, and so $P$ preserves direct sums and $f^*(   S_k  P)=t_k(P)$.  We also want to look at 
 the functor $S_k\sc C\mapright g. \sc C^k$ sending $\{0\hookrightarrow c_1 \hookrightarrow c_2\hookrightarrow\cdots\hookrightarrow c_k)\mapsto(c_1,c_2/c_1,\ldots,c_k/c_{k-1})$.
 We get a commutative diagram of $C_n$-equivariant maps
$$\matrix{U^n(\u{   S_k\sc C};S_k P)  & \mapright g_*.& U^n(\prod_{i=1}^k\u{\sc C};\prod_{i=1}^kP) &  \mapright{{\rm Lemma \ 6.7}} .
& U^n(\u{\sc C};P)^k\cr
  \mapdown{f_*\atop{\rm Example\ 6.3}} .& &\mapup{{\rm Prop \ 6.10}}  .&&\cr
U^n(t_k(\u \sc C);t_k (P))& \mapright{{\rm Def \ 6.11}}.&
U^n(T_k(\u \sc C);T_k (P))&&
\cr}$$
where the labels on the arrows indicate from where it follows that those maps are homotopy
equivalences.  We deduce that $U^n(\u{   S_k\sc C};S_k P)   \mapright \sim.U^n(\u{\sc C};P)^k $
for every $k$.  

We will show that if we apply the maps $g_*$, followed by the projections of Lemma 6.7, levelwise, we get an equivalence $\vert U^n(\u{   S_.\sc C};S_. P) \vert \mapright \sim. 
\vert B_. U^n(\u{\sc C};P)\vert$ compatible with the identifications $U^n(\u{   S_1\sc C};S_1 P) =U^n(\u{\sc C};P)=B_1U^n(\u{\sc C};P)$, which will complete our proof.  (The classifying space
$B_.$ is taken with respect to the operation induced by the abelian group structure on the morphism sets of $\sc C$.)  To show this, we consider that for any simplicial object $X.$ one 
can look at its simplicial path space $(\P X).$ defined by $(\P X)_k=X_{k+1}$  with the original
$\partial_0,\ldots\partial_k$ and $s_0,\ldots,s_k$ as structure maps.  The `extra' degeneracy map
$s_{k+1}$ from $(\P X)_k$ to $(\P X)_{k+1}$ allows us to embed the cone on $\vert (\P X)_.\vert$
in $\vert (\P X)_.\vert$, showing that $\vert (\P X)_.\vert$ is contractible.  The `extra'  boundary map 
$\partial_{k+1}:(\P X)_n\to X_n$ is a simplicial map.  Then we have a commutative diagram
for each $k$
$$\matrix{   
U^n(\u{\sc C};P) 
& \mapright {(s_0)^k_*}. 
& U^n(\u{\P S_k \sc C};\P S_kP) 
& \mapright {(\partial_{k+1})_*}. 
& U^n(\u{S_k\sc C};S_kP) \cr  
 \mapdown.&&\mapdown.&&\mapdown.\cr
U^n(\u{\sc C};P) 
& \mapright .& (\P B)_k(U^n(\u{\sc C};P)=U^n(\u{\sc C};P) ^{k+1}  
& \mapright .& B_k(U^n(\u{\sc C};P)=U^n(\u{\sc C};P) ^k 
\cr}$$
where the vertical arrows are the maps $g_*$ followed by the projections of Lemma 6.7, so we know that they are all equivalences, and the bottom row is the trivial product fibration, inserting
$U^n(\u{\sc C};P) $ in the last coordinate.  The two fibrations are therefore homotopy equivalent.
\bigskip

\noindent
{\caps Proposition  6.15:} Let $\sc C$ be a split small exact category and let $P$ be a bimodule
over the FSP  $\u\sc C$
of the form  $P(a,b)(X)= \sc B(G_1(a),G_2(b))\otimes_{\bf Z}\tilde{\bf Z}[X]$ for some  exact
 functors $G_1, G_2:\sc C\to\sc B$.  Then there is a 
 $C_n$-equivalence 
 $$\lim_{k\to\infty}\Omega^k U^n_0(\u{S^{(k)}\sc C};S^{(k)}P) \mapright\sim.
 \lim_{k\to\infty}\Omega^k U^n(\u{S^{(k)}\sc C};S^{(k)}P).$$
 
  \bigskip
\noindent
{\it Proof}: This is analogous to the proof of Proposition 2.2.3 in [DMc2], but is done simultaneously
in all $n$ blocks.  As usual, we will prove that for all $n$ the given map, which respects the $C_n$ action, is a homotopy equivalence, and the $C_n$-equivalence will follow by induction.

As in [DMc2],  sections 2.0.7 and 2.2.1, we can replace $U^n(\sc C, P)$ by the simplicial abelian group $R_.(\sc C)$ with
$$\eqalign{
R_p(\sc C)=
&\hocolim_{\u X\in I^{n(p+1)}}
\Omega^{\sqcup\u X}
\bigoplus_{(c_{1,0},\ldots,c_{1,p},c_{2,0},\ldots,c_{n,p})\in\sc C^{n(p+1)}}\cr
&
s\sc B(G_1(a_{1,0}),G_2(a_{1, -1}))\otimes
\tilde \Z[s\sc C(a_{1,1}, a_{1,0})]\otimes
\cdots\otimes
\tilde \Z[s\sc C(a_{1,p}, a_{1,p-1})]\otimes
\cr &
\tilde\Z[s\sc B(G_1(a_{2,0}),G_2(a_{2, -1}))]\otimes
\cdots\otimes
\tilde \Z[s\sc C(a_{2,p}, a_{2,p-1})]\otimes
\cdots\otimes
\tilde \Z[s\sc C(a_{n,p}, a_{n,p-1})]
}$$
where $s\sc B$, $s\sc C$ denote the categories of simplicial objects in $\sc B$, $\sc C$,
$$a_{1,-1}=c_{n,p}\otimes\tilde \Z[S^{\sqcup\u X }],$$
for $1< i\leq n$,
$$a_{i,-1}=c_{i-1,p}\otimes\tilde \Z[S^{\sqcup_{(k,l)>(i-1,p)} X_{k,l}}],$$
and for $1\leq i\leq n$, $0\leq j\leq p$
$$a_{i,j}=c_{i,j}\otimes\tilde \Z[S^{\sqcup_{(k,l)>(i,j)} X_{k,l}}].$$
The ordering on pairs of indices is lexicographic.

For every $p$, then, $R_p$ clearly satisfies the requirements of Lemma 2.2.2 in [DMc2],
yielding
$$d_0+d_2\simeq d_1:
 \lim_{k\to\infty}\Omega^k R_p(S^{(k)}S_2\sc C)
 \to
  \lim_{k\to\infty}\Omega^k R_p(S^{(k)}\sc C).$$

For every $p$ we have maps
$$R_0(\sc C)\mapright{s_0^p}. R_p(\sc C)\mapright{d_0^p}. R_0(\sc C)$$
with $d_0^p\circ s_0^p=\id_{R_0(\sc C)}$.  We will show that 
$s_0^p\circ d_0^p\simeq\id_{R_p(\sc C)}$ as well, so that the $R_p(\sc C)$ are all homotopy equivalent to $R_0(\sc C)$.

Since $d_0^{p-1}\circ d_i=d_0^p$, $0\leq i<p$, these $d_i:R_p(\sc C)\to R_{p-1}(\sc C)$ are homotopy equivalences compatible with the equivalence $R_p(\sc C)\to R_{0}(\sc C)$.  
Since  $d_0^{p+1}\circ s_i=d_0^p$ for all  $0\leq i\leq p$, $s_i:R_p(\sc C)\to R_{p+1}(\sc C)$ are
similarly compatible.
One could, similarly to the calculation we are about to make, show that $d_1\circ d_{2}\circ \cdots\circ d_p
:R_p(\sc C)\to R_0(\sc C)$ is a homotopy inverse of $s_0^p$, which implies that 
$d_p:R_p(\sc C)\to R_{p-1}(\sc C)$ is also a homotopy equivalence compatible with the
 equivalence $R_p(\sc C)\to R_{0}(\sc C)$.  Thus
 $$\vert\{ p\mapsto R_p(\sc C)\}\vert
 =\hocolim R_0(\sc C)\simeq R_0(\sc C),$$
 where the homotopy colimit is taken over the co-simplicial category, and the last equivalence is because all the maps involved are  homotopic to the identity and the homotopy colimit of a point over the category is contractible.

To define the homotopy $s_0^p\circ d_0^p\simeq\id_{R_p(\sc C)}$ for a fixed $p$, note that $R_p(\sc C)$ is generated by products of  blocks of maps
$$\eqalign{
& \alpha_{i,0}:G_1(a_{i,0})\to G_2(a_{i-1,p}) 
\cr &
a_{i,0}\mapleft {\alpha_{i,1}}.a_{i,1}\mapleft {\alpha_{i,2}}.a_{i,2}\mapleft {\alpha_{i,3}}.\cdots 
\mapleft {\alpha_{i,p}}.a_{i,p}
}$$
for all $1\leq i\leq n$ (defining $a_{-1,p}=a_{n,p}$).  We can write such a generator as
$\u\alpha=(\alpha_{1,0},\alpha_{1,1},\ldots,\alpha_{1,p},
\alpha_{2,0},\ldots,\alpha_{n,0},\alpha_{n,1},\ldots,\alpha_{n,p}).$  Define 
$$\beta_{i,j}=
\alpha_{i,j}\circ\alpha_{i,j+1}\circ\cdots\circ\alpha_{i,p}:a_{i,p}\to a_{i,j-1}$$
for $1\leq j\leq p$ and
 $$\beta_i=\alpha_{i,0}\circ G_1(\beta_{i,1}):G_1(a_{i,p})\to G_2(a_{i-1,p}).$$
 Then
$$s_0^p\circ d_0^p(\u\alpha)=(\beta_1,\id_{a_{1,p}},\ldots,\id_{a_{1,p}},
\beta_2,\id_{a_{2,p}},\ldots,\id_{a_{2,p}},\ldots,
\beta_n,\id_{a_{n,p}},\ldots,\id_{a_{n,p}}).$$
Now define $t^i_\id (\u\alpha)$ for $1\leq i\leq n$ to consist of
$$\matrix{
&G_2(a_{i-1,p}) 
& \mapleft 0. 
&G_1(a_{i,p})
\cr
&\mapdown{G_2(i_1)}.
&&\mapdown{G_1(i_1)}.
\cr
&G_2(a_{i-1,p}\oplus a_{i-1,p}) 
& \mapleft {G_2(\Delta)\alpha_{i,0}G_1(\pi_2)}. 
&G_1(a_{i,p}\oplus a_{i,0})
\cr
&\mapdown{G_2(\pi_2)}.
&&\mapdown{G_1(\pi_2)}.
\cr
&G_2(a_{i-1,p}) 
& \mapleft {\alpha_{i,0}}. 
&G_1(a_{i,0})
}$$
and
$$\matrix{
&a_{i,p}
&=
&a_{i,p}
&=
&\cdots
&=
&a_{i,p}
\cr
&\mapdown{i_1}.
&&\mapdown{i_1}.
&&&\mapdown{i_1}.
\cr
&a_{i,p}\oplus a_{i,0}
&\mapleft{\id\oplus\alpha_{i,1}}.
&a_{i,p}\oplus a_{i,1}
&\mapleft{\id\oplus\alpha_{i,2}}.
&\cdots
&\mapleft{\id\oplus\alpha_{i,p}}.
&a_{i,p}\oplus a_{i,p}
\cr
&\mapdown{\pi_2}.
&&\mapdown{\pi_2}.
&&&\mapdown{\pi_2}.
\cr
&a_{i,0}
&\mapleft{\alpha_{i,1}}.
&a_{i,1}
&\mapleft{\alpha_{i,2}}.
&\cdots
&\mapleft{\alpha_{i,p}}.
&a_{i,p},}$$
where $i_1$ is the inclusion into the first factor in the direct sum and $\pi_2$ the projection into
the second.  
Then the map
$$\u\alpha\mapsto(t^1_\id(\u\alpha), t^2_\id(\u\alpha),\ldots,t^n_\id(\u\alpha))$$
induces a map
$$T_\id:R_p(\sc C)\to R_p(S_2\sc C).$$
Also, 
define $t^i_{\beta_i} (\u\alpha)$ for $1\leq i\leq n$ to consist of
$$\matrix{
&G_2(a_{i-1,p}) 
& \mapleft \beta_i. 
&G_1(a_{i,p})
\cr
&\mapdown{G_2(\Delta)}.
&&\mapdown{G_1((\id\oplus\beta_{i,1})\Delta)}.
\cr
&G_2(a_{i-1,p}\oplus a_{i-1,p}) 
& \mapleft {G_2(\Delta)\alpha_{i,0}G_1(\pi_2)}. 
&G_1(a_{i,p}\oplus a_{i,0})
\cr
&\mapdown{G_2(\id,-\id)}.
&&\mapdown{G_1(\beta_{i,1},-\id)}.
\cr
&G_2(a_{i-1,p}) 
& \mapleft 0. 
&G_1(a_{i,0})
}$$
and
$$\matrix{
&a_{i,p}
&=
&a_{i,p}
&=
&\cdots
&=
&a_{i,p}
\cr
&\mapdown{(\id\oplus\beta_{i,1})\Delta}.
&&\mapdown{(\id\oplus\beta_{i,2})\Delta}.
&&&\mapdown{\Delta}.
\cr
&a_{i,p}\oplus a_{i,0}
&\mapleft{\id\oplus\alpha_{i,1}}.
&a_{i,p}\oplus a_{i,1}
&\mapleft{\id\oplus\alpha_{i,2}}.
&\cdots
&\mapleft{\id\oplus\alpha_{i,p}}.
&a_{i,p}\oplus a_{i,p}
\cr
&\mapdown{(\beta_{i,1},-\id)}.
&&\mapdown{(\beta_{i,2},-\id)}.
&&&\mapdown{(\id,-\id)}.
\cr
&a_{i,0}
&\mapleft{\alpha_{i,1}}.
&a_{i,1}
&\mapleft{\alpha_{i,2}}.
&\cdots
&\mapleft{\alpha_{i,p}}.
&a_{i,p}.}$$
The map
$$\u\alpha\mapsto(t^1_{\beta_1}(\u\alpha), t^2_{\beta_2}(\u\alpha),\ldots,t^n_{\beta_n}(\u\alpha))$$
induces a map
$$T_\beta:R_p(\sc C)\to R_p(S_2\sc C).$$
Recall that for $c_1\hookrightarrow c_2$ with an identification of $c_2/c_1$ in $S_2\sc C$,
$d_0$ is $c_2/c_1$, $d_1$ is $c_2$, and $d_2$ is $c_1$.  Thus
$$\matrix{
&d_0T_\id=\id && d_0T_\beta=0\cr
&d_1T_\id=d_1T_\beta &&\cr
&d_2T_\id=0 &&d_2T_\beta=s^p_0\circ d_0^p.
}$$
Thus the induced maps
$$\lim_{k\to\infty}\Omega^k R_p(S^{(k)}\sc C) \to \lim_{k\to\infty}\Omega^k R_p(S^{(k)}S_2\sc C)
 \mapright{d_j}. \lim_{k\to\infty}\Omega^k R_p(S^{(k)}\sc C)$$
satisfy
$$\id=d_0T_\id\simeq d_1T_\id=d_1T_\beta\simeq d_2 T_\beta=s^p_0\circ d_0^p,$$
where the homotopies use the homotopy
$d_0+d_2\simeq d_1:
 \lim_{k\to\infty}\Omega^k R_p(S^{(k)}S_2\sc C)
 \to
  \lim_{k\to\infty}\Omega^k R_p(S^{(k)}\sc C)$ mentioned above, along with the vanishing of $d_2T_\id$, $d_0T_\beta$.

\bigskip

\noindent
{\caps Proposition  6.16:} Let $\sc C$ be a split small exact category and let $P$ be a bimodule
over the FSP  $\u\sc C$
of the form  $P(a,b)(X)= \sc B(G_1(a),G_2(b))\otimes_{\bf Z}\tilde{\bf Z}[X]$ for some  exact
 functors $G_1, G_2:\sc C\to\sc B$.  Then there is a 
 $C_n$-equivalence 
 $$\lim_{k\to\infty}\Omega^k \bigvee_{\u c\in S^{(k)}\sc C} U^n_0(\u{S^{(k)}\sc C}\vert_{\u c} ;S^{(k)}P\vert_{\u c}) \mapright\sim.
 \lim_{k\to\infty}\Omega^k U^n_0(\u{S^{(k)}\sc C};S^{(k)}P).$$
 
  \bigskip
\noindent
{\it Proof}: 
As in the proof of Proposition 6.45 above, we replace $U^n_0$ by the abelian group
$R_0$; the corresponding replacement for $\bigvee_{c\in\sc C}U^n_0(\u{\sc C}\vert_c,
P\vert_c)$ is
$$\eqalign{
R'(\sc C)= &
\bigoplus_{c\in\sc C}
\hocolim_{\u X\in I^{n(p+1)}}
\Omega^{\sqcup\u X}
\Bigl(s\sc B(G_1(a'_1),G_2(a'_0))\otimes
\cr &\qquad\qquad\qquad 
\tilde \Z[s\sc B(G_1(a'_2),G_2( a'_1)]\otimes
\cdots\otimes
\tilde \Z[s\sc B(G_1(a'_n),G_2( a'_{n-1})]\Bigr)
\cr 
\simeq &\hocolim_{\u X\in I^{n(p+1)}}
\Omega^{\sqcup\u X}
\bigoplus_{c\in\sc C}
\Bigl(s\sc B(G_1(a'_1),G_2(a'_0))\otimes
\cr &\qquad\qquad\qquad 
\tilde \Z[s\sc B(G_1(a'_2),G_2( a'_1)]\otimes
\cdots\otimes
\tilde \Z[s\sc B(G_1(a'_n),G_2( a'_{n-1})]\Bigr)
}$$
where $a'_0=c\otimes\tilde\Z[S^{\sqcup\u X}],\ \ a'_1=c\otimes\tilde\Z[S^{X_2\sqcup \cdots\sqcup X_n}], \ \ldots\ , a'_{n-1}=c\otimes\tilde\Z[S^{ X_n}],\ \ a'_n=c$.

We have the obvious inclusion 
$$i:R'(\sc C)\mapright i. R_0(\sc C)$$
and a map
$$R_0(\sc C)\mapright j. R'(\sc C)$$
sending each summand corresponding to $(c_1,\ldots, c_n)\in\sc C^n$ in $R_0(\sc C)$ 
to the summand corresponding to $c=c_1\oplus\cdots\oplus  c_n$ in $R'(\sc C)$ via the map induced by the obvious inclusions and projections $c_j\mapright{i_j}.c_1\oplus\cdots\oplus  c_n
\mapright{\pi_j}. c_j$,
$$
s\sc B(G_1(\pi_1),G_2(i_n))\otimes
\tilde \Z[s\sc B(G_1(\pi_2),G_2( i_1)]\otimes
\cdots\otimes
\tilde \Z[s\sc B(G_1(\pi_n),G_2( i_{n-1})].
$$
Since $R_0$ satisfies the requirements of Lemma 2.2.2 in [DMc2],
$$d_0+d_2\simeq d_1:
 \lim_{k\to\infty}\Omega^k R_0(S^{(k)}S_2\sc C)
 \to
  \lim_{k\to\infty}\Omega^k R_0(S^{(k)}\sc C).$$
Define a map $T:R_0(\sc C)\to R_0(S_2\sc C)$ by sending
$$\matrix{
&G_2(a_0)\mapleft\alpha_1. G_1(a_1)
&G_2(a_1)\mapleft\alpha_2. G_1(a_2) 
&\cdots
&G_2(a_{n-1})\mapleft\alpha_n. G_1(a_n)}$$
to 
$$\matrix{
&G_2(a_0)&\mapleft\alpha_1.& G_1(a_1)
&G_2(a_1)&\mapleft\alpha_2. &G_1(a_2) 
&\cdots
&G_2(a_{n-1})&\mapleft\alpha_n.& G_1(a_n)
\cr
&\mapdown{G_2(i_n)}.
&&\mapdown{G_1(i_1)}.
&\mapdown{G_2(i_1)}.
&&\mapdown{G_1(i_2)}.
&&\mapdown{G_2(i_{n-1})}.
&&\mapdown{G_1(i_n)}.
\cr
&G_2(a'_0)&\mapleft{\gamma_1}.& G_1(a'_1)
&G_2(a'_1)&\mapleft{\gamma_2}. &G_1(a'_2) 
&\cdots
&G_2(a'_{n-1})&\mapleft{\gamma_n}.& G_1(a'_n)
\cr
&\mapdown{G_2(p_n)}.
&&\mapdown{G_1(p_1)}.
&\mapdown{G_2(p_1)}.
&&\mapdown{G_1(p_2)}.
&&\mapdown{G_2(p_{n-1})}.
&&\mapdown{G_1(p_n)}.
\cr
&G_2(a''_0)&\mapleft0.& G_1(a''_1)
&G_2(a''_1)&\mapleft0. &G_1(a''_2) 
&\cdots
&G_2(a''_{n-1})&\mapleft0.& G_1(a''_n)
}$$
where $\gamma_j=G_2(i_{j-1})\alpha_jG_1(\pi_j)$ (with $i_0$ defined as $i_n$), and for all 
$0\leq i\leq n$ (with $0$ cyclically identified with $n$ as an index for the $c_i$ as usual),
$$\eqalign{
a_i &= c_i\otimes\tilde\Z[S^{\sqcup_{j>i} X_j}] \cr
a'_i &= c\otimes\tilde\Z[S^{\sqcup_{j>i} X_j}] \cr
a''_i &= \bigoplus_{k\neq i} c_k\otimes\tilde\Z[S^{\sqcup_{j>i} X_j}] .}$$
By construction we have short exact sequences $0\mapright. a_j\mapright {i_j}. a'_j\mapright {p_j}. a''_j\mapright. 0$.

We have $d_0T=0$, $d_2T=\id$, and $d_1T=i\circ j$, so on 
$  \lim_{k\to\infty}\Omega^k R_0(S^{(k)}\sc C)$,
$$\id=d_2T=d_0T+d_2T\simeq d_1T=i\circ j.$$
Now $j\circ i$ is not equal to the identity, but note that the map $T$ above sends $R'(\sc C)$ to
$R'(S_2\sc C)$ and $R' $ also satisfies the conditions of Lemma 2.2.2 in [DMc2], so the formula for $T$ also gives a homotopy $j\circ i\simeq\id_{ \lim_{k\to\infty}\Omega^k R'(S^{(k)}\sc C)}$.
\bigskip




\noindent
{\caps Corollary  6.17:} Let $\sc C$ be a split small exact category and let $P$ be a bimodule
over the FSP  $\u\sc C$
of the form  $P(a,b)(X)= \sc B(G_1(a),G_2(b))\otimes_{\bf Z}\tilde{\bf Z}[X]$ for some  exact
 functors $G_1, G_2:\sc C\to\sc B$.  Then there is a 
 $C_n$-equivalence 
 $$\lim_{k\to\infty}\Omega^k \bigvee_{\u c\in S^{(k)}\sc C} U^n_0(\u{S^{(k)}\sc C}\vert_{\u c} ;S^{(k)}P\vert_{\u c}) \mapright\sim.
 \lim_{k\to\infty}\Omega^k U^n(\u{S^{(k)}\sc C};S^{(k)}P).$$

\bigskip
\noindent{\bf 7. \underbar{Cubical diagrams and analyticity}}

\bigskip
Our goal in this section is to
recall some definitions and terminology from
[G2] which are suitable for our applications
and then to show that the functor $W(F;P\otimes -)$
is $0$-analytic.
Instead of developing the tools to apply the Taylor tower
to the category of all $F$--bimodules, we will
instead define this functor from spaces to spectra, which
we will be able to apply [G2] and [G3] to without
additional modifications. This is done to avoid complications
that are not necessary for our applications.

For $S$ a finite set, let $\sc P(S)$ be the poset
of all subsets of $S$. An $S$--cube (or $n$--cube
where $n=|S|$) in a category $\sc C$ is a functor
$\sc X$ from $\sc P(S)$ to $\sc C$. Thus, a $0$--cube is
an object of $\sc C$, a $1$--cube is a morphism
and a $2$--cube is a commuting square in $\sc C$.

Let $\sc X$ be an $S$--cube of spaces or spectra.
Since $\emptyset$ is initial in $S$,
$\sc X(\emptyset)\rightarrow{\rm \holim}_S\sc X$
is an equivalence. Let $\sc P_0(S)$
be the poset of {\it nonempty} subsets of $S$.
The {\it homotopy fiber}
of $\sc X$ is:
$${\rm hofib}(\sc X) =
{\rm hofib}\ (\sc X(\emptyset)\mapright\alpha.{\rm \holim}_{\sc P_0(S)}
\sc X).$$
We say that $\sc X$ is {\it $k$--Cartesian} if $\alpha$ is
$k$--connected and {\it Cartesian} if $\alpha$ is an equivalence.

Since $S$ is final in $\sc P(S)$,
${\rm colim}_{\sc P(S)}\sc X\rightarrow \sc X(S)$
is an equivalence. Let $\sc P_1(S)$
be the poset of $T\subseteq S$ such that
$T\ne S$.
The {\it homotopy cofiber} of $\sc X$ is:
$${\rm hcofib} (\sc X) = {\rm hcofib}\ ({\rm hocolim}_{\sc P_1(S)}\sc X
\mapright\beta.\sc X(S)).$$
We say that $\sc X$ is {\it $k$--co-Cartesian} if $\beta$
is $k$--connected and {\it co-Cartesian} if $\beta$ is
an equivalence.

\noindent
\underbar{Terminology for $F$-bimodules}:
Let $F$ be a fixed FSP. We now want to modify the above
definitions for functors from spectra to spectra so
we may apply them to the category of $F$--bimodules.
Recall that a morphism of $F$--bimodules is $k$--connected
if
for all $A,B\in\sc O$, the map of associated spectra
$f_{A,B}$ is $k$--connected.
We will call a cubical diagram $\sc X$ of $F$--bimodules
(or functors with stabilization)
$k$--(co)Cartesian if
for all $A,B\in\sc O$
its associated diagram of
spectra $\sc X_{A,B}$ is $k$--(co)Cartesian.

\bigskip
\noindent
{\caps Definition 7.1:} An $S$--cube $\sc X$ of spaces or spectra
is {\it strongly (co-) Cartesian} if each face of
dimension $\geq 2$ is (co-) Cartesian.

\bigskip
\noindent
{\caps Definition 7.2:} Let $F$ be a homotopy functor from
$\sc C$ to $\sc D$. Then $F$ is {\it $n$--excisive},
or satisfies {\it n-th order excision}, if
for every strongly co-Cartesian $(n+1)$--cubical diagram
$\sc X:\sc P(S)\rightarrow\sc C$ the diagram
$F(\sc X)$ is Cartesian.

\noindent
{\it Example} 7.3:
By Proposition 3.4 of [G2], if $M:\sc C^{\times r}\rightarrow\sc D$
is $n_i$--excisive in the $i$-th variable for all $1\leq i\leq r$,
then the composition with the diagonal inclusion
$\sc C\rightarrow \sc C^{\times r}$ is $n$--excisive with
$n=n_1+\cdots+n_r$.
Thus, by Lemma 3.3 of [G2], $U^n(F;\ )$ is an $n$--excisive functor.

\bigskip
\noindent
{\it Example} 7.4:
Spectra
have the nice property that
for $n$--cubes
of spectra,
`$k$--Cartesian' = `(k+n-1)--co-Cartesian'
(1.19 of [G2]) and hence, by definition,
the category of $F$--bimodules
(or functors with stabilization) also shares this property.
In particular, every co-Cartesian diagram of $F$--bimodules
is Cartesian and
thus for any FSP $F$,
the identity functor  from $F$--bimodules to
itself is 1-excisive.
More generally,
given a unital map $f:F\rightarrow F'$ of FSP's,
the functor $f^*$ from $F'$--bimodules to $F$--bimodules
is 1-excisive.

\bigskip
\noindent
{\it Example} 7.5:
The functor $\Omega^{\infty}$ from $F$--bimodules
to $\Omega^{\infty}F$--bimodules
is 1-excisive (since $id\mapright\simeq.\Omega^{\infty}$).  Since $F$ maps to  $\Omega^{\infty}F$,
one can regard the resulting $\Omega^{\infty}F$--bimodules as $F$--bimodules again.
Thus one can replace every strong cofibration
square of $F$--bimodules by an equivalent strong
co-Cartesian square of $F$--bimodules for which the associated spectra sending $n$ to the
value of the functor on $S^n$ are $\Omega$-spectra.

\bigskip\noindent
{\caps Definition/Example 7.6}:
Given a functor with stabilization $Q$, we
let $Q_*$ denote the homotopy functor
$X\mapsto Q(X)$ which is necessarily reduced
($Q(*)\simeq *$). Then
$\Omega^{\infty}Q_*$ is a 1-excisive functor from spaces to spaces
which satisfies the {\it limit axiom}.
We recall from 1.7 of [G1] that a
homotopy functor $G$ satisfies the limit
axiom if for all spaces $X$,
$${\hocolim}_{Y\subseteq X}G(Y)\mapright\simeq.G(X)$$
where the limit system runs over all compact CW-subspaces
and maps the inclusions.

\bigskip
\noindent
{\caps Definition 7.7:} (Goodwillie [G2]) Let $F$ be a homotopy functor from
$\sc C$ to $\sc D$. Then $F$ is {\it stably $n$--excisive},
or satisfies {\it stable n-th order excision}, if
the following is true for some numbers $c$ and $\kappa$:

\noindent
\hang
$E_n(c,\kappa)$: If $\sc X:\sc P(S)\rightarrow \sc C$ is any
strongly co-Cartesian $(n+1)$--cube such that for all $s\in S$
the map $\sc X(\emptyset)\rightarrow \sc X(s)$ is $k_s$--connected
and $k_s\geq \kappa$, then the diagram
$F(\sc X)$ is $(-c+\Sigma k_s)$--Cartesian.

\noindent
If $E_n(c,\kappa)$ holds for all $\kappa$ then we simply say that
$F$ satisfies $E_n(c)$.

\bigskip\noindent
{\it Example} 7.8: Every $n$--excisive functor satisfies $E_n(c)$ for
all $c$.

\bigskip
\noindent
{\it Example} 7.9: If $F_*$ is a simplicial object
of functors from spaces or spectra to spectra satisfying
$E_n(c,\kappa)$, then the realized functor $|F_*|$
(to spectra)
also
satisfies $E_n(c,\kappa)$. This is because
the realization functor is equivalent to a homotopy colimit
which preserves connectivity and commutes with finite
homotopy inverse limits (because we are in spectra).

\bigskip\noindent
{\it Example} 7.10: For any functor with stabilization $P$,
the functor from spaces to functors with stabilization 
sending a space $X$ to the functor with stabilization
$P\otimes X:Y\mapsto P(Y)\wedge X$ is 1-excisive,
and the functor from functors with stabilization to
spaces defined by $P\mapsto {\hocolim}_{X\in I}
Map(S^X,P(S^X))$ is 1-excisive.

\bigskip
\noindent
{\caps Definition 7.11:} (Goodwillie [G2]) The functor $F$ is {\it $\rho$--analytic} if
there is some number $q$ such that $F$ satisfies
$E_n(n\rho-q,\rho+1)$ for all $n\geq 1$.

\bigskip\noindent
{\caps Definition 7.12:} Given an $F$--bimodule $P$, we
write $U^n(F;P\otimes- )$ and $W_{\sc M}(F;P\otimes-)$
for the homotopy functors from spaces to
simplicial functors with stabilization
defined by composition with the functor $P\otimes -$ (of Example 7.10)
from spaces
to $F$--bimodules.
We note that $U^n(F;P\otimes-)$ satisfies the limit axiom
but $W_{\sc M}(F;P\otimes-)$ only does if $\sc M$ is finite
or $P$ is $0$--connective.

\bigskip
\noindent
{\it Example 7.13}: For any $m|n$
and $F$-bimodule $P$,
the functor $U^n(F;P\otimes-)_{hC_m}$ from
space
to functors with stabilization
is $n$--excisive and
satisfies
$E_k(0)$ for all $k\geq 1$ and hence is
$-1$--analytic.

\bigskip\noindent
{\it Proof}: (after 4.4 of [G2])
Since homotopy orbits preserve connectivity and commute
with finite homotopy inverse limits, it suffices to
show the result for $U^n(F;P\otimes- )$.  

By Example 7.5, we may assume $P$ is an $\Omega$ F--bimodule.
Now $U^n(F;P\otimes-)$ is
the diagonal on an $n$--multi-excisive functor so
by Example 7.3 it is $n$--excisive.  By Proposition 3.2 in [G2], being $n$-excisive implies
being $k$-excisive for all $k\geq n$, which implies (as observed in Example 7.8 above) 
that  $U^n(F;P\otimes-)$ satisfies $E_k(c)$ for all $k\geq n$ and all $c$.  Below we will show
that for all $k\geq 1$,  $U^n(F;P\otimes-)$ satisfies $E_k(0)$.  Once we know that, we can take $q=-n$ and get that for all $k\geq 1$,  $U^n(F;P\otimes-)$ satisfies $E_k(-k+n, 0)$ and so it is
$-1$--analytic.

By Example 7.9, to show that $U^n(F;P\otimes-)$ satisfies $E_k(0)$ it suffices to show that
$U^n(F;P\otimes-)_{[m]}$ satisfies $E_k(0)$ for
all $k\geq 1$ and $m\geq 0$, so we fix $m$ and consider $U^n(F;P\otimes-)_{[m]}$.

Let $\sc X$ be a
strongly co-Cartesian
$S$--cube of $F$--bimodules
such that $\sc X(\emptyset)\rightarrow \sc X(s)$
is $k_s$--connected for $s\in S$.
Now for any space $Z$, the spectra
$$U^n(F;P)_{[m]}\wedge \bigwedge^n Z\simeq  U^n(F;P\otimes Z)_{[m]} $$ 
are equivalent.
By example 4.4 of [G2], $\wedge^n\sc X$
is $\Sigma k_s+ n$--co-Cartesian (reduce to
a CW case and consider cells), thus
$U^n(F;P\otimes \sc X)_{[m]} $
is a $\Sigma k_s + n$--co-Cartesian diagram of spectra and so  a $\Sigma k_s$--Cartesian
diagram of spectra.

\bigskip
\noindent
{\caps Proposition 7.14:}
Let $\sc M\subseteq {\bf N}^{\times}$
be such that for all $M\in\sc M$, $\vec M\subseteq\sc M$
(notation as in Proposition 5.6). Then for $F$-bimodules $P$,
$W_{\sc M}(F;P\otimes-)$ satisfies
$E_n(1)$ for all $n\geq 1$ and hence is $0$--analytic.

\bigskip\noindent
{\it Proof}:
Note, if $\sc M\subseteq {\bf N}^{\times}$ is finite, then
there exists an $M\in\sc M$ such that $\sc M$ is
covered by $\vec M$ and $\sc M-M$ (as in 5.6). Thus, by
5.6 we have a natural fibration:
$$U^{M}(F;P\otimes-)_{hC_M}
\mapright.{\holim}
_{\sc M}U(F;P\otimes-)|_{\sc M}
\mapright Res.{\holim}
_{\sc M-M}U(F;P\otimes-)|_{\sc M-M}
.$$
By Example 7.13, $U^{M}_{hC_M}$ satisfies $E_n(0)$.
Now $\sc M-M$ is again such that if $M'\in\sc M-M$
then $\vec M'\subseteq\sc M-M$ and so by
induction on the number of objects of $\sc M$,
${\holim}_{\sc M-M}U(F;P\otimes-)|_{\sc M - M}$
satisfies $E_n(0)$
and hence
${\holim}_{\sc M}U(F;P\otimes-)$ satisfies $E_n(0)$ too.

In general,
$W_{\sc M}(F;P\otimes-)$ can be written as
$$W_{\sc M}(F;P\otimes-)\simeq \holim_{\infty\leftarrow n}
W_{\sc M_n}(F;P\otimes-)$$
for $\sc M_n$ finite and as above.  Now since homotopy inverse limits commute,
the homotopy fiber of $W_{\sc M}(F;P\otimes-)$ on a strongly co-Cartesian $(n+1)$-cube
is the homotopy inverse limit of the homotopy fibers of the $W_{\sc M_n}(F;P\otimes-)$ 
on that cube.
By [BK], XI.7.4, there is a natural short exact sequence
$$0\rightarrow {\rm lim}_{n\to\infty}^1\pi_{i+1}({\rm the\ homotopy\ fiber\ of\ }W_{\sc M_n})
\mapright.\pi_{i}({\rm the\ homotopy\ fiber\ of\ }W_{\sc M})$$
$$\mapright.
{\rm lim}_{n\to\infty}\pi_{i}({\rm the\ homotopy\ fiber\ of\ }W_{\sc M_n})\rightarrow 0$$
and since the homotopy fibers of the $W_{\sc M_n}$ are  $\Sigma k_s$-connected,
the homotopy fiber of  $W_{\sc M}$ must be $-1+\Sigma k_s$-connected.  So
$W_{\sc M}(F;P\otimes-)$ always satisfies $E_n(1)$.

\bigskip
\noindent
{\it Remark}: Certainly Proposition 7.14 holds for more general categories
$\sc M$ but what is proven suffices for our purposes. For example,
using $\vec M$ and the equivalence
$$U^{M}(\ )^{C_M}\mapright\simeq.{\holim}_{\vec M}U(F;P\otimes-)|_{\sc M}$$
($M$ is initial in $\vec M$)
we see that $U^M(F;P\otimes-)^{C_M}$ also satisfies $E_n(0)$.

\bigskip
\noindent{\bf 8. \underbar{Taylor towers and $W_{\sc M}(F;P\otimes-)$}}

\bigskip
We now want to identify the Taylor tower  for
$W(F;P\otimes-)$ and $W^{(p)}(F;P\otimes-)$ when $P$ is
an $F$-bimodule. Technically, what we will write down
is what one might call the Maclaurin  series since
it is  the Taylor tower at the basepoint.
We will be
using universal properties and the results of section 6.

Recall that for a $\rho$--analytic functor $F$ there
is an $n$--excisive homotopy functor called the
$n$-th degree Taylor polynomial which we write as $P_nF$.
We also define:
$$D_nF = {\rm hofib} (P_nF\mapright q_nF.P_{n-1}F).$$

\noindent
{\it Terminology}: An {\it admissible sequence}
$\{x_1,\ldots,x_n\}$ is a (possibly infinite) strictly
increasing sequence of positive integers such that

\item{(i)} $x_1 = 1$

\item{(ii)} If $m|x_j$ then ${x_j\over m}\in
\{x_1,\ldots,x_{j-1}\}$

Note that if $\{x_1,\ldots,x_n\}$ is an admissible sequence
then so is $\{x_1,\ldots,x_{n-1}\}$. Some examples of admissible
sequences are:
$\{1,2,3,\ldots\}$, $\{1,p,p^2,\ldots\}$,
$\{1,p,q,pq\}$ where $p$ and $q$ are prime.
Any multiplicatively closed subset of $N^{\times}$ which contains
all its prime divisors  and $1$ determines an admissible sequence
and every admissible sequence is a subsequence of one obtained
in this manner.

Given an admissible sequence $\{x_1,\ldots,x_n\ldots\}$,
let $\vec\sc X_j\subseteq {\bf N}^{\times}$ be the
full subcategory generated by $\{x_1,\ldots,x_j\}$.
Thus, $\sc X_j\subseteq \sc X_{j+1}$ and
$\sc X_{j+1}$ is covered (as in Proposition 5.6) by $\sc X_{j}$ and
$\vec x_{j+1}$. Let $\sc X$ be the
full subcategory generated by all $x_i$'s.

\bigskip\noindent
{\caps Proposition 8.1}: Given an admissible sequence
$\{x_1,x_2,\ldots\}$ and an $F$--bimodule $P$, then
$$\eqalign{
D_n(W_{\sc X_j}(F;P\otimes-)) \simeq&\
\cases{U^n(F;P\otimes-)_{hC_n}& if $n\in
\{x_1,\ldots,x_j\}$\cr
*&otherwise\cr}
\cr
P_n(W_{\sc X_j}(F;P\otimes-)) \simeq&\
\cases{W_{\sc X_k}(F;P\otimes-) & where $x_k\leq n < x_{k+1}$, $k<j$\cr
W_{\sc X_j}(F;P\otimes-) &  where $x_j\leq n$\cr}
}$$
$$\eqalign{
D_n(W_{\sc X}(F;P\otimes-)) \simeq&\
\cases{U^n(F;P\otimes-)_{hC_n}& if $n\in
\{x_1,\ldots\}$\cr
*&otherwise\cr}
\cr
P_n(W_{\sc X}(F;P\otimes-)) \simeq&\
W_{\sc X_k}(F;P\otimes-)\ \ \ \ {\rm where}\ x_k\leq n < x_{k+1}\cr
}$$
with structure maps $q_n:P_n\rightarrow P_{n-1}$
given by restriction to subcategories.

\bigskip\noindent
{\it Proof}:
In order to ease notation, we will drop $F$ and $P$
from our notation so that $W_{\sc X}$ will
represent $ W_{\sc X}(F;P\otimes-)$.

We recall from [G3] the following facts about
the functor $P_n$ (at the basepoint):
$P_n$ preserves
equivalences and fibrations of functors.  For $F$  $\rho$--analytic, the
natural transformation $F\mapright p_nF.P_nF$
is universal with respect to mapping $F$ to an $n$--excisive functor which
agrees with it to order $n$, that is:
 so that for some constant $q$, 
$p_nF$ is at least $nk+q+k$--connected on $k$--connected
spaces for $k\geq\rho$ (see [G3], just after the proof of Proposition 1.6).

We will use the fibration 
$$U^{x_j}_{hC_{x_j}}\simeq
{\rm hofib}[W_{\sc X_j}\mapright res.W_{\sc X_{j-1}}]$$
from Proposition 5.6.
By Example 3.5 of [G2], $U^n_{hC_n}$ is an $n$--excisive
functor.  If $F''\rightarrow F\rightarrow F'$ is a fibration
of homotopy functors and both $F''$ and $F'$
satisfy $E_n(c)$, then  $F$ satisfies $E_n(c)$ also.
So by induction and Proposition 7.14, $W_{\sc X_j}$ is $0$--analytic and
$x_j$--excisive.

By the universality of $P_nF$, if $F$ is $0$-analytic and $m$-excisive,
the natural transformation $F\to P_nF$ is an equivalence for all $n\geq m$ on all $0$-connected
spaces and so a stable equivalence for those $n$.  Therefore $P_n(W_{\sc X_j})\simeq W_{\sc X_j}$
for all $n\geq x_j$.

Now the connectivity condition in the universal property of $P_n$ shows that $U^n_{hC_n}$ is 
in fact a {\it homogenous} $n$--excisive functor, {\it i.e.}
$$P_k(U^n_{hC_n}) = \cases{U^n_{hC_n}&if $k\geq n$\cr
              *&if $k<n$.\cr
              }$$
Using the fibration of Proposition 5.6 again, this implies that 
$P_k(W_{\sc X_j})\mapright P_k(res).P_k(W_{\sc X_{j-1}})$
is an equivalence for $k < x_j$, so the first two results follow by induction on $j$. 

As in Corollary 5.8, we see that
$W_{\sc X}\mapright res.W_{\sc X_j}$ is
$(x_{j+1}(k+1)-2)$ connected for $X$ $k$--connected.  Since $x_{j+1}\geq x_j+1$, 
this implies $((x_j+1)( k +1)-2)$--connectedness, that is: $(x_jk+x_j-1+k)$--connectedness.
Thus, by the universal property of $p_{x_j}$,
$P_{x_j}W_{\sc X}\mapright res.P_{x_j}W_{\sc X_j}$
is an equivalence. The formula for the layers $D_nW_{\sc X}$ now
follows by induction since we know the layers
$D_nW_{\sc X_j}$ for all $j$ and $n$.

\bigskip
\noindent
{\caps Corollary 8.2:}
For any FSP $F$ and linear bimodule $P$, one has
$$\eqalign{
D_n(W(F;P\otimes-)) \simeq&\ U^n(F;P\otimes-)_{hC_n}\cr
P_n(W(F;P\otimes-)) \simeq&\ W_n(F;P\otimes-)\cr
\ &\ \cr
D_n(W^{(p)}(F;P\otimes-))\simeq&\ \cases{U^{p^k}(F;P\otimes-)_{hC_{p^k}}&if
$n=p^k$\cr
                                *&otherwise\cr}\cr
P_n(W^{(p)}(F;P\otimes-))\simeq&\ W^{(p)}_k(F;P\otimes-)\ \ \ \ {\rm where}\ p^k\leq n <
p^{k+1}\cr
}$$
with structure maps $q_n:P_n\rightarrow P_{n-1}$
given by restriction to subcategories.

\bigskip
\noindent{\bf 9. \underbar{
Relating $K(R\semiprod\tilde M[X])$ to  $W(\u R;\u{M}\otimes \Sigma X)$}}

When $R$ is an associative ring with unit and $M$ is a simplicial $R$-bimodule, Dundas
and McCarthy defined in [DMc1] the invariant $K(R;M)$.  For $M$ discrete, this is the algebraic K-theory
of the exact category whose objects are $(P,\alpha)$, where $P$ is a finitely generated projective $R$-module and $\alpha:P\to P\otimes_R M$ is a right $R$-module map, and whose morphisms
from $(P,\alpha)$ to $(Q,\beta)$ are right $R$-module homomorphisms $f:P\to Q$ such that
$\beta\circ f=(f\otimes {\bf 1}_M)\circ \alpha$.   For $M$ a simplicial $R$-bimodule, this definition
is applied degreewise.  In Theorem 4.1 there, Dundas and McCarthy show that there is a natural homotopy equivalence 
$$K(R\semiprod M)\simeq K(R;B.M)$$
where $B.M\simeq \tilde M[S^1]$ is the bar construction on $M$.   The functoriality of this identification means that the direct summand $K(R)=K(R;0)$ maps compatibly to both sides,
so we have a natural  equivalence on the reduced theories $\tilde K(R\semiprod M)\simeq \tilde K(R;B.M)$ as well.  The argument in [DMc1]
goes on to show that topological Hochschild homology $\THH(R;M)$ is the first derivative at the one-point space $*$ of the functor $X\mapsto\tilde K(R; \tilde M[X]).$

Now Proposition 3.2 in [Mc1] says that the functor $X\mapsto
\tilde K(R;\tilde M[ X])$ is $0$-analytic.  We want to show that for connected $X$, this functor is 
naturally homotopy equivalent to the functor $X\mapsto W(\u R;\u{\tilde M[X] })\simeq W(\u R;\u M\otimes X)$, which is also $0$-analytic by Proposition 6.14 above.  Here $\u R$ is the FSP associated (as in the example after Definition 1.3) to the category having a single object and
$R$ as its morphism set, and $\u M$ is the $\u R$-bimodule sending $X\mapsto \tilde M[X]$.
The equivalence $W(\u R;\u M\otimes X) \simeq W(\u R;\u{\tilde M[X] })$ is by Lemma 6.5 above.

It will turn out to be more convenient to map into $W(\u{\sc P_R};\u M)$ instead of $W(\u R;\u M)$,
where $\sc P_R$ is the category of finitely generated projective right $R$-modules, and $\u M$
is the $\u{\sc P_R}$-bimodule given by
$$\u M(A,B)(X)=\Hom_{\sc M_R}(A;B\otimes_R M)\otimes_{{\bf Z}}\tilde{\bf Z}[X]$$
for all $A,B\in\sc P_R$, $X\in\sc S_*$, where $\sc M_R$ is the category of all right $R$-modules.
The restriction of this $\u M$ on $\sc P_R$ to the full subcategory on the rank $1$ free module (which is isomorphic to the category $R$) is the $\u M$ of the previous paragraph.

\bigskip
\noindent
{\caps Lemma 9.1:} Let $R$ be an associative ring with unit and let $M$ be a discrete $R$-bimodule.  Then we can define maps
$$\tilde K(R;M)\mapright{\beta_n}. U^n(\u{\sc P_R};\u M)$$
for all $n\geq 1$ such that $\Res^{n\over m}\circ\beta_n=\beta_m$ for all $m\vert n$ and such
that $\beta_1$ is the map of Theorem 3.4 in [DMc1].
  \bigskip
\noindent
{\it Proof}:  
 As defined in section 3 of [DMc1],
 $$K(R;M)=\Omega\vert\coprod_{\oc\in S_.\sc P_R}\Hom_{S_.\sc M_R}(\oc, \oc\otimes_R M)\vert.$$
 We will let 
 $$\sc K(R;M)=\coprod_{c\in\sc P_R}\Hom_{\sc M_R}(c,c\otimes_R M).$$
 It is easy to map
 $$\sc K(R;M)\mapright{b_n} .U^n(\u{\sc P_R};\u M)$$
 for all $n\geq 1$ such that $\Res^{n\over m}\circ b_n=b_m$ for all $m\vert n$ and such that 
 $b_1$ is the map [DMc1] used at this level: send
$$\eqalign{\Hom_{\sc M_R}(c,c\otimes_R M)
&\to
\underbrace{\Hom_{\sc M_R}(c,c\otimes_R M)\wedge\cdots\wedge\Hom_{\sc M_R}(c,c\otimes_R M)}_{n \rm\ times}
\cr
\alpha &\mapsto \alpha^{\wedge n},}$$
that is: map $\Hom_{\sc M_R}(c,c\otimes_R M)$ into the multi-simplicial degree $(0,0,\ldots 0)$
part of $U^n(\u{\sc P_R};\u M)$ (where there are no $\u{\sc P_R}$ coordinates, only bimodule
coordinates), into the term corresponding to $\u X=(0,0,\ldots, 0)$ (since $\otimes_\Z\tilde\Z[S^0]=
\otimes_\Z \Z$ does not do anything to an abelian group) by sending $\alpha $ to $n$ copies 
of itself.

For each $k$, on $S_k$ of Waldhausen's S-construction (as reviewed before Proposition 6.14 above) $b_n$ induce maps
$$\eqalign{
 \coprod_{\oc\in S_k\sc P_R}\Hom_{S_k\sc M_R}(\oc, \oc\otimes_R M)
&\to 
U^n(S_k\u{\sc P_R};S_k\u M) \cr
\alpha& \mapsto \alpha^{\wedge n}}$$
and so we get maps
$$
K(R;M)=\Omega\vert\coprod_{\oc\in S_.\sc P_R}\Hom_{S_.\sc M_R}(\oc, \oc\otimes_R M)\vert
\mapright.
\Omega\vert U^n(S_.\u{\sc P_R};S_.\u M)\vert\simeq U^n(\u{\sc P_R};\u M)$$
by the equivalence of Proposition 6.14.  These satisfy $\Res^{n\over m}\circ\beta_n=\beta_m$ for all $m\vert n$ and generalize the construction of [DMc1].

Note that by naturality, these $\beta_n$ send
$K(R)=K(R;0)\to U^n(\u{\sc P_R};\u 0)\simeq *$ and 
and so factor through reduced K-theory maps
$$\tilde K(R;M)\mapright{\beta_n}. U^n(\u{\sc P_R};\u M).$$
\bigskip

Recall that for simplicial $R$-bimodules $N$, $K(R;N)$ is defined by geometrically 
realizing the $K(R;N_n)$'s with respect to the maps induced by $N$'s simplicial structure.
Since the functors $U^n(F;-)$ commute with realizations, we could do the same for
$U^n(\u{\sc P_R};\u N)$.  Therefore, Lemma 9.1 allows us to define maps
$$\tilde K(R;N)\mapright{\beta_n}. U^n(\u{\sc P_R};\u N)$$
for any simplicial $R$-bimodule $N$, satisfying $\Res^{n\over m}\circ\beta_n=\beta_m$ for all $m\vert n$ and generalizing the construction of [DMc1].
Since the $\beta_n$ are compatible with the restriction maps, they define a map
$$\tilde K(R;N)\mapright{\beta}. W(\u{\sc P_R};\u N).$$

We will want to apply this for $N=\tilde M[X]$ for $M$ a discrete $R$-bimodule and $X$ a finite
pointed simplicial set.  Note that the FSP $\u{\tilde M[X]}$ associated to the simplicial $R$-bimodule $\tilde M[X]$ is the same as the FSP $\u{\tilde M}[X]$ of Example (iv) after Definition 1.5,
and as such has an associated spectrum stably equivalent to that of $\u M\otimes X$ of Example (iii), which has been studied in the previous section. The following Theorem will be proved in
several steps. 
 
 \bigskip
\noindent
{\caps Main Theorem 9.2:} Let $R$ be an associative ring with unit and let $M$ be a discrete $R$-bimodule.  Then the natural transformation
$$\tilde K(R;\tilde M[X])\mapright{\beta}. W(\u{\sc P_R};\u{\tilde M}[X])$$
induces an equivalence when $X$ is connected.
Since by Proposition 6.13,
$$W(\u R;\u{\tilde M}[X]) \simeq W(\u{\sc P_R};\u{\tilde M}[X])$$
is a homotopy equivalence (note that in the proof there are maps given in both directions), 
this gives a homotopy equivalence 
$$\tilde K(R;\tilde M[X])\mapright \simeq. W(\u R;\u{\tilde M}[X])$$
for connected $X$.

 \bigskip
\noindent
{\caps Corollary 9.3:} The Taylor tower of the functor $X\mapsto  \tilde K(R;\tilde M[X])$ at $*$
has  $W_n(\u R;\u{\tilde M}[X]) $ as its $n$'th stage, with the tower maps given by category restriction.

  \bigskip
\noindent
{\it Proof of the Corollary}:  In Corollary 8.2 above, we have seen that the $W_n(R;M\otimes-)$ are the finite stages in the Taylor tower of $W(R;M\otimes -)$ at $*$, related by the category restriction maps.  By Lemma 6.5, this means
  that the $W_n(\u R;\u{\tilde M}[X]) $ are the finite stages in the Taylor tower of 
  $W(\u R;\u{\tilde M}[X]) $ at $*$, related by the same maps.
The coefficients in the Taylor tower can be computed by looking arbitrarily close to the space at which we are working.  See Remark 1.1 in [G3] for a discussion of this.  So to find the Taylor
tower of $X\mapsto  \tilde K(R;\tilde M[X])$ it is enough to look at connected $X$, where Theorem
9.2 tells us it agrees with $W(\u R;\u{\tilde M}[X]) $.

  \bigskip
\noindent
{\it Proof of the Theorem}:  
  We will  use a variant of Theorem 5.3 of [G2].  Theorem 5.3 says
that if two  $\rho$-analytic functors $F$ and $G$ have a natural 
transformation between them which induces an equivalence  of the differentials
$D_XF\to D_XG$ at every space $X$, then for $(\rho+1)$-connected maps $X\to Y$ there is a Cartesian 
diagram
$$\matrix{F(X)&\mapright .& G(X)\cr
  \mapdown .& &\mapdown .\cr
 F(Y)&\mapright .&G(Y).\cr}$$
The same proof shows that  if two  $\rho$-analytic functors $F$ and $G$ have a natural 
transformation between them which induces an equivalence  of the differentials
$D_XF\to D_XG$ at every $\rho$-connected space $X$ (i.e. every $X$ for which the map $X\to *$ is $(\rho+1)$-connected---see the comment just after Definition 1.3  in [G2]), then
for every $\rho$-connected $X$ there is a 
Cartesian 
diagram
$$\matrix{F(X)&\mapright .& G(X)\cr
  \mapdown .& &\mapdown .\cr
 F(*)&\mapright .&G(*).\cr}$$
So we need to show that the  natural transformation
$\beta$
 induces an equivalence on the derivatives at all $0$-connected spaces.
Since
$\tilde K(R;\tilde M[ *]) \simeq W(\u{\sc P_R};\u{\tilde M}[*])  \simeq *$, this will imply that for any $0$-connected $X$, $\beta$ is an equivalence.

So let $X$ be connected; we should consider spaces $Y\to X$ over $X$;
%
but for simplicity of writing, we would like to eliminate the spaces from our calculation.  In the following sections, we will prove
  \bigskip
\noindent
{\caps Technical Lemma 9.4:} Let $R$ be an associative ring with unit, and let $M, N$ be two simplicial
$R$-bimodules. If $N$ is $k$-connected, $\beta$ induces
a $2k$-connected map
$${\rm hofib} \bigl(\tilde K(R;B_.M\oplus B_.N) \to \tilde K(R;B_.M)\bigr)
\to 
{\rm hofib} \bigl(  W(\u{\sc P_R};\u{B_.M}\oplus\u{B_.N})\to  W(\u{\sc P_R};\u{B_. M})\bigr).$$

\bigskip

Having this lemma lets us conclude the proof of Theorem 9.2: because of the obvious equivalences
$$\tilde M[X\vee A]\cong \tilde M[X]\oplus \tilde M[A]$$
and, for connected spaces $Z$, $\tilde M[Z]\cong B_.(\Omega \tilde M[Z])$ (for $\Omega$ a simplicial
model of the loop space), Technical Lemma 9.4 would tell us that 
$${\rm hofib} \bigl(\tilde K(R;\tilde M[ X\vee S^k]) \to \tilde K(R;\tilde M[ X])\bigr)
\to 
{\rm hofib} \bigl(  W(\u{\sc P_R};\u{\tilde M}[X\vee S^k])\to  W(\u{\sc P_R};\u{\tilde M}[X])\bigr)$$
is $2(k-1)$-connected for all $k\geq 1$.  But by the definition of a derivative at a space  $X$ of a functor
 (at the basepoint of $X$) from [G1] (for functors to spaces) and [G2] (section 5), the spectrum
 $k\mapsto {\rm hofib} \bigl(\tilde K(R;\tilde M[ X\vee S^k]) \to \tilde K(R;\tilde M[ X])\bigr)
$
is equivalent to the derivative of the functor $\tilde K(R;\tilde M[ -])$ at $X$, and similarly the spectrum 
 $k\mapsto{\rm hofib} \bigl(  W(\u{\sc P_R};\u{\tilde M}[X\vee S^k])\to  W(\u{\sc P_R};\u{\tilde M}[X])\bigr)
 $
 is equivalent to the derivative of the functor $W(\u{\sc P_R};\u{\tilde M}[-])$ at $X$.  The $2(k-1)$-equivalences of the $k$'th spaces in these two spectra make the two of them equivalent, and thus the two derivatives agree, as we needed to show.
 
\bigskip
\noindent{\bf 10. \underbar{Reduction of Technical Lemma 9.4 to $THH$}}

\bigskip
In this section we are going to prove that for $k$-connected $N$,
the  two fibers  in  Technical Lemma
9.4 in are both $2k$--equivalent  to $THH(\u{R\semiprod M},\u{B.N})$. In section 11 we will prove that the
trace map $\beta$ induces the abstract equivalence obtained in this section.

\bigskip
\noindent
{\caps Lemma 10.1:} Let $R$ be an associative ring with unit, and let $M$ and $N$ be
two simplicial $R$-bimodules, then
$$K(R\semiprod M;B.N)\simeq {\rm hofib}(\tilde K(R;B.M\oplus B.N)\mapright.\tilde K(R;B.M)).$$

\bigskip\noindent
{\it Proof.} Using the isomorphisms  $B.M\oplus B.N\simeq B.(M\oplus N)$ and
$R\semiprod (M\oplus N)\cong (R\semiprod M)\semiprod N$
the result follows from the commuting diagram whose vertical maps are
equivalences by 4.1 of [DMc1] as well as the identification of
the homotopy fiber on the bottom row:
$$\matrix{{\rm hofib}&\mapright.&\tilde K(R;B.M\oplus B.N)&\mapright.&\tilde K(R;B.M)\cr
\mapdown\simeq.&&\mapdown\simeq.&&\mapdown\simeq.\cr
\tilde K(R\semiprod M;B.N)&\mapright.&\tilde K((R\semiprod M)\semiprod N)
&\mapright.&\tilde K(R\semiprod M)\cr}
\leqno{(1)}$$

\bigskip
We are interested in $\tilde K(R\semiprod M;B.N)$ in a $2k$--connected range when
$N$ is $k$-connected. 
By Theorem 3.4 of [DMc1], the trace map
$\tilde K(R\semiprod M;B.N)\mapright.THH(\u{R\semiprod M};\u{B.N})$ is
$2(k+1)$--connected if $N$ is $k$-connected. 
We let $\gamma$ be the composite
obtained using the natural splitting of the rows in (1):
$$\matrix{\tilde K(R,B.(M\oplus N))&\mapright\simeq. &\tilde K(R\semiprod (M\oplus N))
\cong \tilde K((R\semiprod M)\semiprod N) \cr
\mapdown\gamma.&&\mapdown.\cr
THH(\u{R\semiprod M};\u{B.N})&\mapleft 2(k+1).&\tilde K(R\semiprod M;B.N)
}\leqno{(2)}$$

\bigskip\noindent
{\caps Proposition 10.2:} For a ring $R$ and simplicial $R$-bimodules $M$ and $N$,
using the natural map $\gamma$ from (2), we obtain
$$\tilde K(R;B.M\oplus B.N)\mapright 2(k+1). THH(\u{R\semiprod M};\u{B.N})\times \tilde K(R;B.M).$$

\bigskip
In order to identify the homotopy fiber of the map
$$W(\u R;\u{B.M}\oplus \u{B.N})\mapright. W(\u R;\u{B.M})$$
we first observe that by the multi-linearity of $U^n(F;-)$, for any FSP $F$ and $F$-bimodules
$P_0$ and $P_1$, we have
a $C_n$-equivariant decomposition
$$U^n(F;P_0\oplus P_1) \mapright\simeq.
\prod_{\alpha\in {\rm Hom}_{\rm Sets}(\{1,2,\ldots,n\},\{0,1\})} U^n(F;P_{\alpha(1)},\ldots,P_{\alpha(n)}).$$
Thus, if $P_1$ is $k$-connected we have a
$2k$-connected $C_n$--equivariant map
$$U^n(F;P_0\oplus P_1)\mapright. U^n(F;P_0)\times (C_n)_+\wedge U^n(R;P_0,\dots,P_0,P_1).\leqno{(3)}$$
Taking $C_n$ fixed points, this gives a map
$$U^n(F;P_0\oplus P_1)^{C_n}\mapright. U^n(F;P_0)^{C_n} \times  U^n(R;P_0,\dots,P_0,P_1).$$
We obtain maps for all $k$ dividing $n$,
$$U^n(F;P_0\oplus P_1)^{C_n}\mapright res.
U^k(F;P_0\oplus P_1)^{C_k}\mapright. U^k(R;P_0\ldots,P_0,P_1)$$
and hence maps $\epsilon_n$:
$$U^n(F;P_0\oplus P_1)^{C_n}\mapright\epsilon_n. U^n(F;P_0)^{C_n}\times
\prod_{k|n} U^k(F;P_0,\ldots,P_0,P_1)$$
which take the restriction maps for the fixed points of the $U$'s applied
to $P_0\oplus P_1$ to the restriction maps of the $U$'s applied to $P_0$ and
the projections on the product.

\bigskip\noindent
{\caps Lemma 10.3:} The map of homotopy inverse limits produced by the
$\epsilon_n$'s:
$$\matrix{
W(\u R;\u{B.M}\oplus \u{B.N})&\!\!\mapright\epsilon.&\!\!\!
{\rm holim}_{{\bf N}^{\times}} \!\!\left(U^n(\u R;\u{B.M})^{C_n}\ \!\!\times \!\prod_{k|n} \!U^k(\u R;\u{B.M},\ldots,\u {B.M},\u{B.N})\right)\cr
&&\mapdown\cong.\cr
&&W(\u R;\u{B.M})\ \times\ \prod_{a=0}^{\infty} U^a(\u R;\u{B.M},\ldots,\u{B.M},\u{B.N})\cr
}$$
is $2(k+1)-1$ connected if $N$ is $k$-connected.

\bigskip\noindent
{\it Proof.} It is enough to show the map for each $W_n$ to
${\rm holim}_{\{\leq n\}}$ is 
$2(k+1)$ connected. By Corollary 5.7 it is enough to show that the
$C_n$-homotopy orbits of the maps in (3) are $2(k+1)$ connected,
which they are since homotopy orbits preserve connectivity. 

\bigskip
In order to relate 
$\prod_{a=0}^{\infty} U^a(\u R;\u{B.M},\ldots,\u{B.M},\u{B.N})$ to
$THH(\u{R\semiprod M},\u{B.N})$, a proposition motivated by the
result in [L] will be useful. 
For $R\rightarrow S$ a map of FSP's and
$M$ an $S$-bimodule
we can form a bi-simplicial spectrum
$$[p]\times [q]\ \mapsto\ \ 
U^{p+1}_q(R;S,\ldots,S,M)$$
where the $q$-direction has the usual (diagonal) simplicial structure of $U$ and 
the $p$ direction has the evident simplicial structure defined so that
$$U^{*+1}_0(R;S\ldots,S,M)\ \cong\ THH(S,M).$$

\bigskip
\noindent
{\caps Proposition 10.4:} If  $R\rightarrow S$ is a map of FSP's over the same underlying set and
$M$ is an $S$-bimodule,
then the natural map
$$U^1(S,M)\mapright. U^{*+1}(R,S^*,M)$$
given by the inclusion of the zero simplicial dimension 
$$U^{*+1}(R,S^*,M)_{[0]} = U^1(S,M)_{[*]}$$
is an equivalence.

\bigskip\noindent
{\it Proof.} Since both theories are linear in the bimodule $M$ variable
it is enough to show that the map is an equivalence for
$M = S\otimes X\otimes S$ for $X$ a spectrum (and the tensor product taken over the sphere
spectrum).  More general bimodules $M$ can
be resolved by bimodules of this form, using the functor $X\mapsto S\otimes X\otimes S$
from spectra to $S$--bimodules and its adjoint, the forgetful functor.
In this case one can``break''  the circles to get
$$U^{n+1}(R,S^n,M)_{[*]}\mapright\simeq. X\otimes
\overbrace{S\hat\otimes_RS\hat\otimes\cdots\hat\otimes_RS}^{n+2\;{\rm times}}.$$
The induced structure maps
(from the unused simplicial direction) are simply $X\otimes (\ )$
applied to the standard bimodule resolution of $S$ as an $R$-algebra
(i.e. multiplication on the``insides'', no twists) and hence since this
is homotopy equivalent to $S$ again (using the extra degeneracy map),
we get $X\otimes S$. 

\bigskip\noindent
{\caps Corollary 10.5:} 
For a ring $R$ and simplicial $R$-bimodules $M$ and $N$,
$$THH(\u{R\semiprod M};\u{B.N})\simeq\prod_{a=0}^{\infty} U^a(\u R;\u{B.M},\ldots,\u{B.M},\u{B.N}).$$

\bigskip\noindent
{\it Proof.} By Proposition 10.4, we have
$$THH(\u{R\semiprod M};\u{B.N})\mapright\simeq. U^{*+1}(\u R;\u{R\semiprod M}^{*},\u{B.N}).$$
Using the multi-linearity of the $U$'s, , if we let $S^k = \Delta^k/\partial$, we have by the calculations at the end of the next section an equivalence $\rho_A$
of bisimplicial spectra
$$U^{*+1}(\u R;\u{R\semiprod M}^{*},\u{B.N})\ \cong
\prod_a^{\infty} U^a(\u R;\u M,\ldots,\u M,\u{B.N})\otimes S^a.$$
Using the fact that $S^a\simeq \overbrace{S^1\wedge\ldots\wedge S^1}^{a\;times}$,
$\u{B.M}\simeq\u M\otimes S^1$ and the multilinearity of the $U$'s, we
have that each $U^a(\u R;\u M,\ldots,\u M,\u{B.N})\otimes S^a$ is equivalent
to $U^a(\u R;\u{B.M},\ldots,\u{B.M},\u{B.N})$ and hence the result. 

\bigskip\noindent
{\caps Proposition 10.6:} For a ring $R$ and simplicial $R$-bimodules $M$ and $N$, using the
map $\epsilon$ and the equivalences of Lemma 10.3 and Corollary 10.5
$$W(\u R;\u{B.M}\oplus\u{ B.N})\mapright 2k+1. THH(\u {R\semiprod M};\u{B.N})\times W(\u R;\u{B.M}).$$

\bigskip
\noindent{\bf 11. \underbar{{Proof that the trace induces the equivalence on section 10}}}

\bigskip
In section 10 we showed that the two fibers used in Technical Lemma 9.4
agree in a $2k+1$ range if the bimodule $N$ being considered was $k$--connected.
In this section we will show that the trace map defined in section 9 induces
the equivalence on the fibers in a $2k$ range as suggested by section 10. 
We now recall, from [Mc1], primarily page 218 and [DMc1], section 4 details that
allow us to construct an unstable model for the composite $\gamma$ used
in Proposition 10.2.
The following notation will be convenient for this purpose.

\bigskip\noindent
{\bf Notation.} 
If $\sc P$ is a category of diagrams of projective right $R$--modules and $\sc M$ is the
category of diagrams of the same form of right $R$--modules, we will for brevity write
$$R(P)=\Hom_\sc P(P,P)$$
and
$$M_P=\Hom_\sc M(P,P\otimes_R M)$$
for any $P\in\sc P$ and any $R$--bimodule $M$.

\bigskip
We recall, using this notation, that for a discrete ring $R$,
$$\tilde K(R;M) = {\rm stabilization\  w.r.t.\ Waldhausen's\ S\!\!-\!\!construction\ of}\ \ 
\sc P\mapsto\bigvee_{P\in\sc P}M_P$$
on $\sc P_R$.
As in the proof of Theorem 3.4 in [DMc1], the point here is that
$$K(R;M)=\Omega \vert\coprod_{\oc\in S_.\sc P_R} \Hom_{S_.\sc M_R}(\oc,\oc\otimes M)\vert$$
and 
$$K(R)=K(R;0)=\Omega\vert S_.\sc P_R\vert.$$
By [W], the S-construction stabilizing maps induce equivalences 
$$\Omega\vert S_.\sc P_R\vert \mapright\simeq.
\Omega^2\vert S_.^{(2)}\sc P_R\vert \mapright\simeq.
\Omega^3\vert S_.^{(3)}\sc P_R\vert \mapright\simeq.\cdots$$
and similarly 
$$\Omega \vert\coprod_{\oc\in S_.\sc P_R} \Hom_{S_.\sc M_R}(\oc,\oc\otimes M)\vert \mapright\simeq.
\Omega^2 \vert\coprod_{\oc\in S^{(2)}_.\sc P_R} \Hom_{S^{(2)}_.\sc M_R}(\oc,\oc\otimes M)\vert \mapright\simeq.
\cdots.$$
As is shown in the proof of Theorem 3.4 of [DMc1], for each $p$ the map
$$\eqalign{
{\rm hofib}\bigl( & \vert\coprod_{\oc\in S^{(p)}_.\sc P_R} \Hom_{S^{(p)}_.\sc M_R}(\oc,\oc\otimes M)\vert \to 
\vert S_.^{(p)}\sc P_R\vert \bigr)\cr & 
\to {\rm hocofib}\bigl(\vert S_.^{(p)}\sc P_R\vert \mapright 0.  \vert\coprod_{\oc\in S^{(p)}_.\sc P_R} \Hom_{S^{(p)}_.\sc M_R}(\oc,\oc\otimes M)\vert \bigr)\cr &
\simeq\vert\bigvee_{\oc\in S^{(p)}_.\sc P_R} \Hom_{S^{(p)}_.\sc M_R}(\oc,\oc\otimes M)\vert
}$$
is $(2p-3)$--connected (where $0$ is the map which sends every $\oc$ to the zero map $\oc\to\oc\otimes M$), and thus the homotopy cofiber spectrum is the same as the spectrum $\tilde K(R;M)$.

By Proposition 6.14 and  Corollary 6.17, we know that  
$$U^n(R,M)^{C_n}=  {\rm stabilization\  w.r.t.\ the\ S\!\!-\!\!construction\ of}\ \ \ 
\sc P\mapsto \bigvee_{P\in \sc P}U^n_0(\u\sc P\vert_P;
\u M\vert_P)^{C_n}$$
on $\sc P_R$.
We observe that
$${\rm Hom}_{\sc P_{R\semiprod M}}(P\otimes_R(R\semiprod M),
P\otimes_R(R\semiprod M)\otimes_{R\semiprod M}N)\ \cong\ 
{\rm Hom}_{\sc P_R}(P,P\otimes_RN)$$
In other words,
$$N_{P\otimes_RR\semiprod M}\ \cong\ N_P$$
Similarly, 
$$\eqalign{{\rm Hom}_{\sc P_{R\semiprod M}}&(P\otimes_R(R\semiprod M),
P\otimes_R(R\semiprod M))\cr
&  \cong \Hom_{{\sc P}_R}(P, P\otimes_R(R\semiprod M)) \cong
 {\rm Hom}_{\sc P_R}(P,P)\oplus {\rm Hom}_{\sc P_R}(P,P\otimes_RM)}$$
and taking into account the composition product, we can write this as
$$R(P\otimes_R(R\semiprod M))\cong R(P)\semiprod M_P.$$
Given $m\in M_P$, we let
$(1,m)$ be the obvious  isomorphism in 
$R(P)\semiprod M_P$.
We note that $(1,m)\circ (1,m') = (1, m+ m')$
and for $n\in N_P$,
$(1,m)\circ n = n = n\circ (1,m)$. Thus, we have a representation
$$M_P\mapright 1+\star.GL(R(P)\semiprod M_P)$$
whose image acts trivially on $N_P$.  Here $GL$ is used to indicate the invertible elements.
 
We obtain a natural simplicial map
$$N_P\times B_*M_P\mapright B_*(1+\star).  
THH_*(\u{R(P)\semiprod M_P};\u{N_P})$$
where $\u{R(P)\semiprod M_P}$, $\u{N_P}$ are viewed as an FSP  on the category consisting of one point  associated to the ring $R(P)\semiprod M_P$, and the bimodule on that FSP associated to the 
$R(P)\semiprod M_P$--bimodule $N_P$.
Note that
$$\u{R(P)\semiprod M_P}\cong\u{\sc P_{R\semiprod M}}\vert_{P\otimes_R(R\semiprod M)}$$
and
$$\u{N_P}\cong \u N\vert_{P\otimes_R(R\semiprod M)}.$$

\bigskip\noindent
{\caps Lemma 11.1:} The natural transformation 
$$\bigvee_{P\in\sc P_R}THH_*(\u{R(P)\semiprod M_P};\u{N_P})
\mapright \star\otimes_RR\semiprod M.
\bigvee_{P\in\sc P_{R\semiprod M}}THH_*(\u{R(P)\semiprod M_P};\u{N_P})$$
(obtained by tensoring the indexing category  $\star\otimes_RR\semiprod M$)
is an equivalence after stabilization.

\bigskip\noindent
{\it Proof.} The domain of this map can be written as
$$\bigvee_{P\in\sc P_R}U_*^1(\u{R(P)\semiprod M_P};\u{N_P})
=\bigvee_{P\in\sc P_R}U_*^1(\u{R\semiprod M}\vert_P,\u{N}\vert_P)
=\bigvee_{P\in\sc P_R}U_0^{*+1}(\u{\sc P_R}\vert_P;\u{R\semiprod M}\vert_P^*,\u{N}\vert_P)
,$$
and so stabilizes to $U^{*+1}(\u{\sc P_R};\u{R\semiprod M}^*,\u{N})$ by the non--equivariant analog
of Corollary 6.17 which allows different bimodules in the $*+1$ bimodule positions (and is proved in exactly the same way). 

The target of the map can be written as 
$$\bigvee_{P\in\sc P_{R\semiprod M}}\!\!\!\!\!\!U_*^1(\u{R(P)\semiprod M_P};\u{N_P})
=\!\!\!\!\!\!\bigvee_{P\in\sc P_{R\semiprod M}}\!\!\!\!\!\!\!\!U_*^1(\u{R\semiprod M}\vert_P,\u{N}\vert_P)
=\!\!\!\!\!\!\bigvee_{P\in\sc P_{R\semiprod M}}\!\!\!\!\!\!\!\!U_0^{*+1}(\u{\sc P_{R\semiprod M}}\vert_P;\u{R\semiprod M}\vert_P^*,\u{N}\vert_P)
,$$
and similarly stabilizes to $U^{*+1}(\u{\sc P_{R\semiprod M}};\u{R\semiprod M}^*,\u{N})$.
On $\sc P_{R\semiprod M}$, $\u{R\semiprod M}$ is the same as $\sc P_{R\semiprod M}$, but we write
it in this way to keep the parallel clear.  If we abbreviate the map $\star\otimes_RR\semiprod M$ to
$t:\sc P_R\to\sc P_{R\otimes M}$, the map in the statement of the lemma stabilized to the obvious map induced by $t$, noting that $\u{R\semiprod M}$ on $\sc P_R$
is the same as $t^*$ of $\u{R\semiprod M}$ on $\sc P_{R\semiprod M}$ and that $\u N$ on $\sc P_R$
is the same as $t^*$ of $\u{N}$ on $\sc P_{R\semiprod M}$.  Thus the map we are interested in is the
bottom horizontal map in the commutative diagram
$$\matrix{
U^{*+1}(\u{\sc P_{R}}\vert_R;(t^*\u{R\semiprod M})\vert_R^*,(t^*\u{N})\vert_R)
&\mapright t_*.
&U^{*+1}(\u{\sc P_{R\semiprod M}}\vert_{R\semiprod M};\u{R\semiprod M}\vert_{R\semiprod M}^*,\u{N}\vert_{R\semiprod M})\cr
\mapdown.&&\mapdown.\cr
U^{*+1}(\u{\sc P_{R}};(t^*\u{R\semiprod M})^*,t^*\u{N})
&\mapright t_*.
&U^{*+1}(\u{\sc P_{R\semiprod M}};\u{R\semiprod M}^*,\u{N})\cr
}$$
where the vertical maps are induced by the inclusion.   By the non-equivariant analog
of Proposition 6.13 which allows different bimodules in the $*+1$ bimodule positions (and is proved in exactly the same way), these vertical maps are equivalences,  so to show that the bottom horizontal map is an equivalence, it suffices to show that the top horizontal map is.  The top horizontal map is 
a map of $U$'s of FSP's and bimodules over a one point set, and can also be written as the map
$$U^{*+1}(\u{R};\u{R\semiprod M}^*,\u{N})\mapright t_*.U^{*+1}(\u{R\semiprod M};\u{R\semiprod M}^*,\u{N})$$
induced by $t$ on the FSP's.  It is an equivalence because of the diagram
$$\matrix{
U^{*+1}_0(\u{R};\u{R\semiprod M}^*,\u{N})
&=
&U^{*+1}_0(\u{R\semiprod M};\u{R\semiprod M}^*,\u{N})\cr
\mapdown.&&\mapdown.\cr
U^{*+1}(\u{R};\u{R\semiprod M}^*,\u{N})
&\mapright t_*.
&U^{*+1}(\u{R\semiprod M};\u{R\semiprod M}^*,\u{N})\cr
}$$
where the vertical maps are equivalences by the non-equivariant analog
of Proposition 10.4 which allows different bimodules in the $*+1$ bimodule positions (and is proved in exactly the same way).

\bigskip
\noindent
{\caps11.2. A Model for $\gamma$:}
The stabilization of the composite
$$\bigvee_{P\in\sc P}N_P\times B_*M_P\mapright 1+\star. \!\!
\bigvee_{P\in\sc P}THH_*(R(P)\semiprod M_P;N_P)
\mapright \star\otimes_RR\semiprod M. \!\!\!\!\!
\bigvee_{P\in\sc P_{R\semiprod M}}THH_*(R(P)\semiprod M_P;N_P)$$
is a natural transformation from $K(R;N\oplus B.M)$ to 
$K(R\semiprod M;N)$. When $N$ is replaced by $B.N$, this 
(by [DMc1] and [Mc1]) is a model for the natural transformation $\gamma$
used in Proposition 10.2. 

\bigskip
By the discussion in the beginning of section 9, for any $R$-bimodule $M$, the map of spaces (not spectra) from
$M_P$ to $U^n_0(\u{\sc P_R},\u M)^{C_n}$ 
$$\eqalign{
m\ \mapsto &\ m\wedge m\wedge\cdots\wedge m\cr
\in &[\ \u M\vert_P(S^0)\wedge\cdots\wedge \u M\vert_P(S^0)]^{C_n}\cr
\mapright. &\ 
[{\rm holim}_{I^{k+1}}
\Omega^{X_0\sqcup\ldots\sqcup X_k}
(\u M\vert_P(S^{X_0})\wedge \u M\vert_P(S^{X_1})\wedge\ldots\wedge
\u M\vert_P(S^{X_k})]^{C_n}\cr
}$$
stabilizes to the map $\tilde K(R;M)\mapright.U^n(\u {\sc P_R};\u M)^{C_n}$ used in the
construction of the trace map $\beta$. 
Using this and  the definition of the maps $\epsilon_n$ in Lemma 10.3 we obtain
the following. 

\bigskip\noindent
{\caps11.3 A Model for the $W$ fiber map:} The composite
$$K(R;M\oplus N)\mapright \beta. W(\u{\sc P}_R;\u {M\oplus N})\mapright.\prod_{a=1}^{\infty} U^a(\u{\sc P_ R};\u M^{a-1},\u N)$$
is equivalent to the stabilization of the natural transformation determined by
the product of the maps $\beta_a$
$$\beta_a(m\times n) = \  \overbrace{m\wedge\cdots\wedge m}^{(a-1)\;times}\wedge n $$
$$(M\oplus N)_P \cong
M_P\oplus N_P\mapright\beta_a.
U^a_0(\u {\sc P_R}\vert _P;\u M\vert _P^{a-1},
\u N\vert_P)$$
for all $P\in\sc P_R$.

\bigskip
By the models in 11.2 and 11.3 and the equivalence in Corollary 10.5 it will suffice to to prove the following 
proposition to finish the proof of Technical Lemma 9.4.

\bigskip\noindent
{\caps Proposition 11.4:} For $R$ a ring and $M$ and $N$ $R$-bimodules, there are maps $\rho_A$
which make 
the following diagram commute:
$$\matrix{B.M_P\times N_P&\mapright 1+\star.&THH(\u{R(P)\semiprod M_P},\u{ N_P})\cr
\mapdown \prod_a \beta_a.&&\mapdown\simeq.\cr
\prod_{a=0}^{\infty} U^a_0(\u {R(P)};\u{B.M_P},\ldots,\u {B.M_P},\u{N_P} )&\mapright \prod_a \rho_a.&
\prod_{a=0}^{\infty} U(a)_0(\u {R(P)};\u{M_P},\u{N_P})\cr
\mapdown.&&\mapdown\simeq.\cr
\prod_{a=0}^{\infty} U^a_{\star}(\u {R(P)};\u{B.M_P},\ldots,\u {B.M_P},\u{N_P} )&
\mapright \prod_a \rho_a\ \simeq.
&
\prod_{a=0}^{\infty} U(a)_{\star}(\u {R(P)};\u{M_P},\u{N_P})\cr
}$$
where we define bisimplicial spectra $U(a)(\u {R(P)};\u{M_P},\u{N_P})$ by
$$U(a)^*_{\star}(\u {R(P)};\u{M_P},\u{N_P}) = 
 \prod_{1\leq j_1 <\cdots < j_a\leq *} U^{*+1}_{\star}(\u {R(P)};F_1,\ldots,F_*,\u {N_P})
$$
with
$$
  F_t =\cases{\u {M_P}&if $t\in\{j_1,\ldots,j_a\}$\cr
              \u {R(P)}&otherwise.\cr}
$$
In other words, if we break $\u{R(P)\semiprod M_P}$ as a bimodule down into $\u {R(P)}\oplus \u{M_P}$, then
 $U(a)(\u {R(P)};\u{M_P},\u{N_P})$ is the part of $THH(\u{R(P)\semiprod M_P},\u{ N_P})$
which contains exactly $a$ $\u{M_P}$'s.
\bigskip 
We recall from Proposition 10.4 the equivalence
$$THH_*(\u{R(P)\semiprod M_P},\u {N_P})\cong
U^{*+1}_0(\u{R(P)};\u{R(P)\semiprod M_P},\u {N_P})
\mapright\simeq.
U^{*+1}_{\star}(\u{R(P)};\u{R(P)\semiprod M_P},\u {N_P}).$$
By the multi-linear property of the $U$, we see that we have
a natural equivalence of bi-simplicial spectra
$$U^{*+1}_{\star}(\u{R(P)};\u{R(P)\semiprod M_P},\u {N_P})\mapright\simeq.\prod_{a=0}^{\infty}
U(a)^*_{\star}(\u {R(P)};\u{M_P},\u{N_P}) .$$

\bigskip
\noindent
{\caps Lemma 11.5:} For all $a$, $U(a)^*_0(\u {R(P)};\u{M_P},\u{N_P})\mapright\simeq.U(a)^*_{\star}(\u {R(P)};\u{M_P},\u{N_P})$.

\bigskip\noindent
{\it Proof.}
We have a commutative diagram:
$$\matrix{U^{*+1}_0(\u{R(P)};\u{R(P)\semiprod M_P},\u {N_P})&\mapright\simeq.&\prod_{a=0}^{\infty}U(a)^*_0(\u {R(P)};\u{M_P},\u{N_P})\cr
\mapdown\simeq.&&\mapdown\simeq.\cr
U^{*+1}_{\star}(\u{R(P)};\u{R(P)\semiprod M_P},\u {N_P})&\mapright\simeq.&\prod_{a=0}^{\infty}U(a)^*_{\star}(\u {R(P)};\u{M_P},\u{N_P})\cr
}$$
where the horizontal maps are equivalences by the decomposition into homogeneous
pieces as explained above and the left vertical map is an equivalence by  Proposition 10.4. Thus, the
right vertical map which makes the diagram commute, namely the inclusion, is an equivalence. However, this inclusion is  a product
of maps and hence each of them is an  equivalence as well.

\bigskip
Composing $1+\star$ with the projection equivalence from
$THH(\u{R(P)\semiprod M_P},\u {N_P})\mapright\simeq. \prod_{a=0}^{\infty} U(a)(\u {R(P)};\u{M_P},\u{N_P})$
we obtain simplicial maps
$$\eta^a:B_*M_P\times N_P\mapright. U(a)_0^*(\u {R(P)};\u{M_P},\u{N_P})\leqno{(1)}$$

In order to better express the maps $\eta^a$ we make the following
observations. 
Let ${\rm Surj}_{\Delta}([a],[k])$ be the set of
surjective monotone maps from $[a]$ to $[k]$. 
For $\sigma\in {\rm Surj}_{\Delta}([a],[k])$ and
$1\leq j\leq k$, we will write 
$\mu_j(\sigma) = {\rm min}\{\sigma^{-1}(j)\}$. 
We have an
isomorphism of 
sets from ${\rm Surj}_{\Delta}([k],[a])$ to
$\{1\leq j_1 <\ldots < j_k\leq a\}$ given
by sending $\sigma\in {\rm Surj}_{\Delta}([a],[k])$ to
$\{\mu _1(\sigma),\ldots, \mu _k(\sigma)\}$.
With these conventions, 
the maps $\eta^a$ are the composites:
$$(n,m _1\times\cdots\times m _k)$$ 
$$\eqalign{
&\ \mapsto\ 
\prod_{\sigma\in {\rm Surj}_{\Delta}([k],[a])}
\sigma^*(n\wedge m_{\mu _1(\sigma)}\wedge\cdots\wedge
m_{\mu _k(\sigma)})\cr
\in&\ \Omega^{0\sqcup\ldots\sqcup 0}
 \prod_{1\leq \mu _1(\sigma)<\ldots < \mu _k(\sigma)\leq a}
( \u{N_P}(S^{0})\wedge{F_1}[{S^{0}}]\wedge\ldots\wedge
{F_k}[{S^{0}}])\cr
\mapright.&\ 
{\rm holim}_{I^{k+1}}
\Omega^{x_0\sqcup\ldots\sqcup x_k}
\prod_{1\leq \mu _1(\sigma)<\ldots < \mu _k(\sigma)\leq a}
( \u{N_P}(S^{x_0})\wedge{F_1}[{S^{x_1}}]\wedge\ldots\wedge
{F_k}[{S^{x_k}}])\cr
=&\  U(a)^k_0(\u {R(P)};\u{M_P},\u{N_P})\cr
}$$
where $F_t$ are as, as before, $\u{M_P}$ if $t=\mu_j(\sigma)$ for some $j$ and $\u {R(P)}$ otherwise,
and where $\sigma^*(n\wedge m_{\mu _1(\sigma)}\wedge\cdots\wedge
m_{\mu _k(\sigma)})$ has $n$ in the $0$'th coordinate, $m_{\mu_j(\sigma)}$ in the $\mu_j(\sigma)$ coordinate, $1\leq j\leq a$, and $1_R$ in the others.

We now define, for all $a$,  bi-simplicial maps (which are equivalences)
$$U^a_{\star}(\u {R(P)}; (B_*\u{ M_P})^a,\u{N_P})\to U(a)^*_{\star}(\u {R(P)};\u{M_P},\u{N_P}).$$ 
In order to do this, we first make a few general remarks
about simplicial constructions and then apply them
to this specific application. 

Let $S^k = \Delta^k/\partial$ as a simplicial
set. That is, $\Delta^k_{[n]} = {\rm Hom}_{\Delta}([n],[k])$ and $\partial$ is the
usual subsimplicial set determined by the $k-1$ subskeleton.
We also have the simplicial set
$$\overbrace{S^1\wedge\cdots\wedge S^1}^{k\;  times}$$
by which we mean the diagonal of the obvious $k$-fold multisimplicial set. 
We know that after realization, this is homeomorphic to 
$S^k$ but there is not a simplicial map that realizes to a homeomorphism.
There are exactly $k!$ simplcial
maps from $\overbrace{S^1\wedge\cdots\wedge S^1}^{k\;  times}$ to
$S^k$ which realize to homotopy equivalences, one for each of the
non-degenerate cells of $\overbrace{S^1\wedge\cdots\wedge S^1}^{k\;  times}$
in dimension $k$. We will be wanting to use one of these simplicial maps.

We can first look at $S^k$. It has one point in every simplicial dimension less than $k$. 
In dimension $k$ it has one non-degenerate element corresponding to the
identity on $[k]$. In dimension n larger than k it has an element for every
ordered surjection from $[n]$ to $[k]$ and one other point, $*$. Thus,
we can write
$$S^k_{[n]} = {\rm Surj}_{\Delta}([n],[k])_+$$

Now we want to write down simplicial maps from
$\overbrace{S^1\wedge\cdots\wedge S^1}^{k\;  times}$ to $S^{k}$ which
are homotopy equivalences. The critical dimension is $k$. 
Given a sequence $(\tau_1,\ldots,\tau_k)$ of surjections $[n]\to [1]$, 
we can associate a sequence of non-zero positive integers
by looking at the cardinality of the inverse image of $0$, i.e.
$\{|\tau^{-1}_1(0)|,|\tau^{-1}_2(0)|,\ldots,|\tau^{-1}_k(0)|\}$.
However, we can also express this as
$\{\mu _1(\tau_1),\mu _1(\tau_2),\ldots,\mu _1(\tau_k)\}$
since $\mu _1(\tau) = {\rm min}(\tau^{-1}(1)) = |\tau^{-1}(0)|$. 
A sequence
$$(\tau_1,\ldots,\tau_k)\in \overbrace{S^1\wedge\cdots\wedge S^1_{[k]}}^{k\;  times}$$ 
is non-degenerate if and only if the sequence
$\{\mu _1(\tau_1),\mu _1(\tau_2),\ldots,\mu _1(\tau_k)\}$
has no repeated terms.

We define $\alpha$ to be the simplicial map from
$ \overbrace{S^1\wedge\cdots\wedge S^1}^{k\;  times}$ to
$S^k$ as follows:
On the $k-1$ skeleton it must be trivial,
in simplicial dimension $k$ it takes the non-degenerate
simplex whose sequence is $(1, 2, \ldots, k)$ to the identity in
$S^k_{[k]}$ and all others to the basepoint. This determines
a map of the $k$-skeletons and hence by extension of degeneracies
a map of the simplicial sets (as both have only degenerate simplicies in
dimensions greater than $k$).

We actually need to understand $\alpha$ explicitly in all simplicial dimensions. 
To do this, we simply need to see how degeneracies operate in both settings. 
When we do, we find that we can describe $\alpha$ as follows:
a simplex ${\tau_1,\ldots,\tau_k}\in 
 \overbrace{S^1\wedge\cdots\wedge S^1_{[n]}}^{k\;  times}$
 is sent to the basepoint if its associated sequence
 $\{m_1(\tau_1),m_1(\tau_2),\ldots,m_1(\tau_k)\}\ \in\{1,2,\ldots,n\}^k$
is not monotone increasing. If it is monotone increasing, then it
is sent to the monotone surjection $\alpha(\tau_1,\ldots,\tau_n)\in
{\rm Surj}_{\Delta}([n],[k])$ given by:
$$\mu _j(\alpha(\tau_1,\ldots,\tau_n)) = \mu _1(\tau_j).$$

We now return to constructing  bi-simplicial equivalences
$$
\rho_a: U^a_{\star}(\u{R(P)};B_*\u{ M_P},\ldots,B_*\u{ M_P},\u{ N_P})\ \mapright.\  U(a)^*_{\star}(\u{R(P)};\u{ M_P},\u{ N_P}).$$ 
We observe that because $U^n$ commutes with realizations in each of its bimodule variables, we have a simplicial map
which is an equivalence in each simplicial dimension:
$$U^a_{\star}(\u{R(P)};B_*\u{ M_P},\ldots,B_*\u{ M_P},\u{ N_P})\mapright\simeq.
U^a_{\star}(\u{R(P)};\u{ M_P},\ldots,\u{ M_P},\u{ N_P})\otimes
 \overbrace{S^1\wedge\cdots\wedge S^1}^{a\;  times}
.\leqno{(2)}$$
Applying the simplicial map $\alpha$ we obtain an equivalence
$$U^a_{\star}(\u{R(P)};\u{ M_P},\ldots,\u{ M_P},\u{ N_P})\otimes
 \overbrace{S^1\wedge\cdots\wedge S^1}^{a\;  times}
 \mapright\alpha\ \simeq.
U^a_{\star}(\u{R(P)};\u{ M_P},\ldots,\u{ M_P},\u{ N_P})\otimes S^a.$$
We have a bi-simplicial map which is a homeomorphism 
$$ U^a_{\star}(\u{R(P)};\u{ M_P},\ldots,\u{ M_P},\u{ N_P})\otimes S^a_*\mapright.
U(a)^*_{\star}(\u{R(P)};\u{M_P},\u{N_P})$$
$$\eqalign{ &\prod_{\sigma\in{\rm Surj}_{\Delta}([k],[a])} U^a_{\star}(\u{R(P)};\u{ M_P},\ldots,\u{ M_P},\u{ N_P})
\cr &  \mapsto\ 
\prod_{\sigma\in{\rm Surj}_{\Delta}([k],[a])} \sigma^*(U^a_{\star}(\u{R(P)};\u{ M_P},\ldots,\u{ M_P},\u{ N_P}))
= U^k_{\star}(\u{R(P)};\u{ M_P},\u{ N_P}).
}\leqno{(3)}$$
We define $\rho_a$ to be the composite of these equivalences.

\bigskip\noindent
{\it Proof of Proposition 11.4.}\ 
It is enough to show $\eta^a = \rho_a\circ \beta_a$ for all $a$
where $\eta^a$ is the piece of $1+\star$ as defined in (1) and 
$\beta_a$ was defined in 11.3.
Let 
$(m_1,\ldots,m_k,n)\in B_kM_P\times N_P$. Then
$\beta_a (m_1,\ldots,m_k,n)$ is the image in the homotopy colimit
of the element 
$n\wedge (m_1\times\cdots\times m_k)^{\wedge a}$. Via the
equivalence in (2), this element is mapped to the element
in the homotopy colimit repesented by
$$\prod_{\gamma_1,\ldots,
\gamma_a\in {\rm Surj}_{\Delta}([k],[1])}
(n\wedge m_{\mu _1(\gamma_1)}\wedge m_{\mu _1(\gamma_2)}\wedge\ldots\wedge
m_{\mu _1(\gamma_a)})
$$
The simplicial map $\alpha$ composed with the map in (3) will send this
element to the image of the element in the colimit of
$\prod_{\sigma\in {\rm Surj}_{\Delta}([k],[a])}
\sigma^*(n\wedge m_{\mu _1(\sigma)}\wedge\cdots\wedge
m_{\mu _k(\sigma)})$ which is the image of $\eta^a$.

\bigskip\noindent{\bf Bibliography}
\item{[B]}{\caps Marcel  B\"okstedt}, {\it Topological Hochschild homology}, Preprint, 1985. 

\item{[BHM]}{\caps Marcel  B\"okstedt, Wu-Chung Hsiang, and Ib Madsen}, {\it The cyclotomic trace and algebraic K-theory of spaces}, Invent. Math, {\bf 111} (1993), 465--539. 

\item{[BK]}{\caps A. K. Bousfield and D. M. Kan}, Homotopy limits, completions and localizations,
Lecture Notes in Mathematics, Vol. 304. Springer-Verlag, Berlin-New York, 1972.

\item{[D]} {\caps Bjorn Ian Dundas}, {\it Relative K-theory and topological cyclic homology.}, Acta Math., {\bf 179} (1997), no. 2, 223--242. 

\item{[DMc1]} {\caps Bjorn Ian Dundas and Randy McCarthy}, {\it Stable K-theory and
topological Hoch\-schild homology}, Annals of Mathematics, {\bf 140} (1994), 685--701.

\item{[DMc1--err]} {\caps Bjorn Ian Dundas and Randy McCarthy}, {\it Erratum Stable K-theory and
topological Hochschild homology}, Annals of Mathmatics, {\bf 142} (1995), 425--426.

\item{[DMc2]}{\caps Bjorn Ian Dundas and Randy McCarthy},{\it Topological Hochschild homology of ring functors and exact categories},  J. Pure Appl. Algebra  {\bf 109}  (1996),  no. 3, 231--294.

\item{[G]}{\caps Thomas G. Goodwillie}, {\it Notes on the cyclotomic trace}, notes from a lecture series given at MSRI during the spring of 1990.

\item{[G1]}{\caps Thomas G. Goodwillie}, {\it Calculus I. The first derivative of pseudoisotopy theory},  K-Theory  {\bf 4}  (1990),  no. 1, 1--27.
\item{[G2]}{\caps Thomas G. Goodwillie}, {\it  Calculus II. Analytic functors},  $K$-Theory  {\bf 5}  (1991/92),  no. 4, 295--332.
\item{[G3]}{\caps Thomas G. Goodwillie}, {\it Calculus III. Taylor series},  Geom. Topol.  {\bf 7}  (2003), 645--711.

\item{[L]}{\caps Ayelet Lindenstrauss}, {\it  A relative spectral sequence for topological Hochschild homology of spectra},  J. Pure Appl. Algebra  {\bf 148}  (2000),  no. 1, 77--88. 

\item{[M]}{\caps Ib Madsen}, {\it Algebraic $K$-theory and traces},  Current developments in mathematics, 1995 (Cambridge, MA),  191--321, Int. Press, Cambridge, MA, 1994.

\item{[Mc1]} {\caps Randy McCarthy}, {\it Relative algebraic K-theory and 
topological cyclic homology}, Acta Math., 179 (1997), 197--222.

\item{[Mc2]} {\caps Randy McCarthy}, {\it A Chain Complex for the Spectrum Homology of the
Algebraic K-theory of an Exact Category}, Fields Institute Communications, Vol. {\bf 16}, 1997, 199--220.

\item{[PW]} {\caps Teimuraz Pirashvili and Friedhelm Waldhausen}, {\it Mac Lane homology and topological Hochschild homology},  J. Pure Appl. Algebra,  {\bf 82}, (1992), no.1, 81--98.

\item{[W]} {\caps Friedhelm Waldhausen}, {\it Algebraic K-theory of spaces}, Springer Lecture Notes in Math.
{\bf 1126} (1985), 318--419.

\bye